\documentclass[a4paper]{amsart}
\usepackage{graphicx, amsmath, amssymb, amsfonts}
\usepackage{longtable}
\usepackage[all]{xy}
\numberwithin{equation}{section}
\theoremstyle{plain}
 \newtheorem{thm}{Theorem}[section]
 \newtheorem{cor}[thm]{Corollary}
 \newtheorem{lem}[thm]{Lemma}
 \newtheorem{prop}[thm]{Proposition}
\theoremstyle{definition}
 \newtheorem{defn}[thm]{Definition}
 \newtheorem{exmp}[thm]{Example}
\theoremstyle{remark}
 \newtheorem{rem}[thm]{Remark}

\DeclareMathOperator{\Hom}{Hom}
\DeclareMathOperator{\Aut}{Aut}
\DeclareMathOperator{\Out}{Out}
\DeclareMathOperator{\rank}{rank}

\def\ccirc{\text{\raisebox{1pt}{\tiny$\circledcirc$}}}
\def\dcirc{\text{\raisebox{1pt}{\tiny$\odot$}}}
\def\vin{\text{\rotatebox{90}{$\in$}}}

\begin{document}
\title{A classification of subsystems of a root system}%
\author{Toshio Oshima}

\address{Graduate School of Mathematical Sciences,
University of Tokyo, 3-8-1, Komaba, Meguro-ku, Tokyo 153-8914, Japan}
\email{oshima@ms.u-tokyo.ac.jp}
\begin{abstract}
We classify isomorphic classes of the homomorphisms of a root system
$\Xi$ to a root system $\Sigma$ which do not change Cartan integers.
We examine several types of isomorphic classes defined by the Weyl group of
$\Sigma$, that of $\Xi$ and the automorphisms of $\Sigma$ or $\Xi$
etc.
We also distinguish the subsystem generated by a subset of a
fundamental system.
We introduce the concept of the dual pair for root systems which helps
to study the action of the outer automorphism of $\Xi$ on
the homomorphisms.
\end{abstract}
%%%%%%%%%%%%%%%%%%%%%%%%%%%%%%%%%%%%%%%%%%%%%%%%%%%%%%%%%%%%%%%%%%%%%
\maketitle
\tableofcontents
\section{Introduction}\label{sec:intro}
Root systems were introduced by W.~Killing and E.~Cartan for the study of 
semisimple Lie algebras and now they are basic in several fields of 
mathematics.
In this note a subsystem of a root system means a subset of a root system
which is stable under the reflections with respect to the roots in the subset.
The purpose of this note is to study subsystems of a root system.
It is not difficult to classify the subsystems if the root system is 
of the classical type but we do it in a unified way.
The method used here will be useful in particular when the root system
is of the exceptional type.

Let $\Xi$ and $\Xi'$ be subsystems of a root system $\Sigma$.
We define that $\Xi'$ is equivalent to $\Xi$ by $\Sigma$ and 
we write $\Xi\underset{\Sigma}{\sim}\Xi'$
if $w(\Xi)=\Xi'$ with an element $w$ of the Weyl group $W_\Sigma$ 
of $\Sigma$.
By the classification in this note we will get complete answers 
to the following fundamental questions 
(cf.~Remark~\ref{rem:ans} for the answers).

\medskip
Q1. What kinds of subsystems of $\Sigma$ exist as 
abstract root systems?

Q2. 
Suppose $\Xi'$ is isomorphic to $\Xi$ as abstract root systems,
which is denoted by $\Xi'\simeq\Xi$.
How do we know $\Xi'\underset{\Sigma}\sim\Xi$?

Q3. How many subsystems of $\Sigma$ exist which are equivalent to
$\Xi$?

Q4. Does the outer automorphism of $\Xi$ come from $W_\Sigma$?

Q5. Suppose $\sigma$ is an outer automorphism of $\Xi$ which stabilizes
   every irreducible component of $\Xi$.
   Is $\sigma$ realized by an element of $W_\Sigma$?

Q6. Suppose that $\Xi$ is transformed to $\Xi'$ by an outer 
  automorphism of $\Sigma$.  Is $\Xi\underset{\Sigma}\sim\Xi'$ valid?

Q7. Is $\Xi$ equivalent to a subsystem $\langle\Theta\rangle$ generated
by a subset $\Theta$ of a fundamental system $\Psi$ of $\Sigma$?
How many elements exist in 
$\{\Theta\subset\Psi\,;\,\langle\Theta\rangle\underset{\Sigma}{\sim}\Xi\}$?

For example, Q4 may be interesting if $\Xi$ has irreducible 
components which are mutually isomorphic to each other.  
An orthogonal system is its typical example (cf.~Remark~\ref{rem:orth}).

The first question of Q7 is studied by \cite{AL} and the answer is given there when $\Xi$ is irreducible (cf.~Remark~\ref{rem:fund} iii)).

To answer these questions we will study subsystems as follows.

Let $\Xi$ and $\Sigma$ be reduced root systems
and let $\Hom(\Xi,\Sigma)$ denote the set of maps of $\Xi$ to $\Sigma$
which keep the Cartan integers
$2\frac{(\alpha|\beta)}{(\beta|\beta)}$ invariant for the roots $\alpha$ and $\beta$.
Since the map is injective and its image is a root system,
the image is a subsystem of $\Sigma$ isomorphic to $\Xi$.

Let $W_\Xi$ and $W_\Sigma$ denote the Weyl groups of $\Xi$
and $\Sigma$ respectively and put
$\Aut(\Xi)=\Hom(\Xi,\Xi)$ and $\Aut(\Sigma)=\Hom(\Sigma,\Sigma)$.
We will first study the most refined classification, that is, 
$W_\Sigma\backslash\!\Hom(\Xi,\Sigma)$ after the review of the standard 
materials for root systems in \S\ref{sec:not}.
In \S\ref{sec:thm} we will give Theorem~\ref{thm} which reduces
% and Corollary~\ref{cor:thm} which reduce
the structure of $W_\Sigma\backslash\!\Hom(\Xi,\Sigma)$ to a simple 
graphic combinatorics in the extended Dynkin diagrams.
It is a generalization of the fact that an element of
$W_\Sigma\backslash\!\Aut(\Sigma)$ corresponds to a graph 
automorphism of the Dynkin diagram associated to $\Sigma$ 
(cf.~Example \ref{exmp:A})
and will be proved in \S\ref{sec:prf} after the preparation in 
\S\ref{sec:lem}.

In \S\ref{sec:dual} we define the dual pair of subsystems,
which helps us to study the action of $\Aut(\Xi)$ on $\Hom(\Xi,\Sigma)$.
In \S\ref{sec:table} we have the table of all the non-empty $\Hom(\Xi,\Sigma)$
with irreducible $\Sigma$.
The table gives the numbers of the elements of the cosets
\begin{gather*}
 W_\Sigma\backslash\!\Hom(\Xi,\Sigma),\  
 \Aut(\Sigma)\backslash\!\Hom(\Xi,\Sigma),\\
 W_\Sigma\backslash\!\Hom(\Xi,\Sigma)/\!\Aut(\Xi),\ 
 W_\Sigma\backslash\!\Hom(\Xi,\Sigma)/\!\Aut'(\Xi)
\end{gather*}
and the number of the subsystems generated 
by subsets of a fundamental system of $\Sigma$ which correspond
to a coset.
Here $\Aut'(\Xi)$ is the subgroup of $\Aut(\Xi)$ defined by the direct
product of the automorphisms of the irreducible components of $\Xi$.
The table also determines certain closures of $\Xi$ (cf.~Definition~
\ref{def:dual}, \ref{defn:S-closed}).

In many cases $\#\bigl(W_\Sigma\backslash\!\Hom(\Xi,\Sigma)/\!\Aut(\Xi)\bigr)=1$, 
which is equivalent to say that the subsystems of $\Sigma$ which are 
isomorphic to $\Xi$ form a single $W_\Sigma$-orbit.  
We will also distinguish the orbits when the number is larger than one.

In \S\ref{sec:rem} we give some remarks obtained by our study.
For example, Q4 will be examined for the orthogonal systems
of the root systems of type $E_7$ and $E_8$ .

In \S\ref{sec:Dynkin} we give the extended Dynkin diagrams and roots
of the irreducible root systems following the notation in \cite{Bo},
which is for the reader's convenience and will be constantly used in 
this note.
A proof of the classification of the root systems is also given for the
completeness (cf.~Proposition~\ref{prop:clasify} and Remark~\ref{rem:cls} iv)).

Dynkin \cite{Dy} classified regular subalgebras of a simple Lie algebra
in his study of semisimple subalgebras. 
The classification is given by Table 9 and Table 11 in \cite{Dy}.
In Table 11, $A_6+A_2$ and the second one of $A_7+A_1$ should be replaced by 
$E_6+A_2$ and $E_7+A_1$, respectively.
These tables describe the classification of 
$\Aut(\Sigma)\backslash\!\Hom(\Xi,\Sigma)/\!\Aut(\Xi)$ for $S$-closed 
subsystems (cf.~Definition~\ref{defn:S-closed})
in our table in \S\ref{sec:table} (cf.~Remark~\ref{rem:last}~ii))
and were obtained from diagrams given by successive
procedures removing roots from extended Dynkin diagrams. 
The procedure is the way to classify maximal $S$-closed subsystems used by \cite{BS}
(cf.~Remark~\ref{rem:maximal}).
The maximal $S$-closed subsystems are also classified by \cite{W}.
Our classification based on Theorem~\ref{thm} gives
a more refined classification of $W_\Sigma\backslash\!\Hom(\Xi,\Sigma)$. 
In fact we give a simple algorithm to give the complete representatives
of the coset $W_\Sigma\backslash\!\Hom(\Xi,\Sigma)$.

The author would like to thank E.~Opdam for pointing out
\eqref{eq:LF} and related errors in the table in \S10.
%%%%%%%
\section{Notation}\label{sec:not}
In this section we review the root systems and fix the notation 
related to them.
All the materials in this section are elementary and found in
\cite{Bo}.

Fix a standard inner product $(\ |\ )$ of $\mathbb R^n$ and
the orthonormal basis $\{\epsilon_1,\ldots,\epsilon_n\}$ of $\mathbb R^n$.
For $\alpha\in\mathbb R^n\setminus\{0\}$ the \textsl{reflection $s_\alpha$ with respect to $\alpha$} is defined by
\begin{equation}
 \begin{matrix}
  s_\alpha: & \mathbb R^n & \to  & \mathbb R^n\\
            & \vin        & & \vin\\
            & x           & \mapsto & s_\alpha(x) := 
     x-2\frac{(\alpha|x)}{(\alpha|\alpha)}\alpha
 \end{matrix}
\end{equation}

\begin{defn}
A \textsl{reduced} root system of rank $n$ is a finite subset $\Sigma$ 
of $\mathbb R^n\setminus\{0\}$ which satisfies
\begin{gather}
 \mathbb R^n=\textstyle\sum_{\alpha\in\Sigma}\mathbb R\alpha,\\
 s_\alpha(\Sigma)=\Sigma\qquad(\forall\alpha\in\Sigma),\\
 2\frac{(\alpha|\beta)}{(\alpha|\alpha)}
  \in\mathbb Z
   \qquad(\forall\alpha,\,\beta\in\Sigma),\\
 \mathbb R \alpha\cap\Sigma 
  = \{\pm\alpha\}\qquad(\forall\alpha\in\Sigma).\label{eq:reduced}
\end{gather}
In general the rank of a root system $\Sigma$ is denoted by $\rank \Sigma$.
\end{defn}
%%%
\begin{rem}
{\rm i)}
In this note any non-reduced root system, which doesn't satisfy 
\eqref{eq:reduced}, doesn't appear except in \S\ref{sec:Dynkin}
and hereafter for simplicity a root system always means a \textsl{reduced} root 
system.

{\rm ii)}
We use the notation $\mathbb N$ for the set $\{0,1,2,\ldots\}$ of 
non-negative integers.
\end{rem}
%%%

\begin{defn}\label{defn:fund}
A \textsl{fundamental system} $\Psi$ of the root system $\Sigma$ of rank
$n$ is a finite subset $\{\alpha_1,\ldots,\alpha_n\}$ of $\Sigma$ which
satisfies
\begin{gather}
  \mathbb R^n =\mathbb R\alpha_1+\mathbb R\alpha_2+\cdots+\mathbb R\alpha_n,\\
  \alpha =\sum_{j=1}^n m_j(\alpha)\alpha_j\in\Sigma\, \Rightarrow\, 
 \bigl(m_1(\alpha),\ldots,\ldots,m_n(\alpha)\bigr)\in\mathbb N^n 
 \text{ or }-\mathbb N^n.\label{eq:coeff}
\end{gather}
The fundamental system $\Psi$ exists for any root system $\Sigma$ and 
the root $\alpha\in\Sigma$ is \textsl{positive} (with respect to $\Psi$)
if $m_j(\alpha)\ge0$ for $j=1,\ldots,n$, 
which is denoted by $\alpha>0$.
\end{defn}
\begin{defn}\label{def:perp}
Let $\Theta$ be a finite subset of $\Sigma$ and put
\begin{align}
  W_\Theta&:=\langle s_\alpha\,;\,\alpha\in\Theta\rangle =
   \text{the group generated by }
  \{s_\alpha\,;\,\alpha\in\Theta\},\\
  W&:=W_\Sigma=W_\Psi,\\
  \langle \Theta \rangle &:= W_\Theta\Theta,\\
  \Theta^\perp &\,= \Theta^\perp\cap\Sigma
    :=\{\alpha\in\Sigma\,;\,(\alpha|\beta) =
  0\quad(\forall\beta\in\Theta)\}.
\end{align}
The reflection group $W$ generated by $s_\alpha$ for $\alpha\in\Sigma$
is called the \textsl{Weyl group} of $\Sigma$.
A subset $\Xi$ of $\Sigma$ is called a \textsl{subsystem} of $\Sigma$ if
$s_\alpha(\Xi)=\Xi$ for any $\alpha\in\Xi$. Then $\Xi$ is a root system
with $\rank\Xi=\dim\sum_{\alpha\in\Xi}\mathbb R\alpha$.
Note that $\langle\Theta\rangle$ and $\Theta^\perp$ are subsystems of $\Sigma$
and
\begin{equation}
 \rank\langle\Theta\rangle+\rank\Theta^\perp\le \rank\Sigma.
\end{equation}

\end{defn}
\begin{defn}
A map $\iota$ of a root system $\Xi$ to a root system $\Sigma$ 
is a \textsl{homomorphism} if $\iota$ keeps the Cartan integers:
\begin{equation}
 2\frac{(\iota(\alpha)|\iota(\beta))}
      {(\iota(\alpha)|\iota(\alpha))} = 
 2\frac{(\alpha|\beta)}{(\alpha|\alpha)}
 \quad(\forall\alpha,\,\beta\in\Xi).
\end{equation}
In this case $\iota$ is injective and $\iota(\Xi)$ is a subsystem
of $\Sigma$.

The set of all homomorphisms of $\Xi$ to $\Sigma$ is denoted
by $\Hom(\Xi,\Sigma)$ and define
\begin{equation}
  \Aut(\Sigma) := \Hom(\Sigma,\Sigma).
\end{equation}
Note that $W_\Sigma$ and $W_\Xi$ naturally act on $\Hom(\Xi,\Sigma)$
and
\begin{align}
  \iota\circ s_\alpha = s_{\iota(\alpha)}\circ\iota
  \quad(\iota\in\Hom(\Xi,\Sigma),\,\alpha\in\Xi).
\end{align}

Two homomorphisms $\iota$ and $\iota'$ of $\Xi$ to $\Sigma$ 
are \textsl{isomorphic} if
\begin{equation}
  \iota' = w\circ \iota
\end{equation}
for a suitable $w\in W_\Sigma$ and we define 
\begin{align}
 \overline\Hom(\Xi,\Sigma)&:=W_\Sigma\backslash\!\Hom(\Xi,\Sigma)
 \simeq W_\Sigma\backslash\!\Hom(\Xi,\Sigma)/W_{\Xi},\\
 \Out(\Sigma)&:= W_{\Sigma}\backslash\!\Aut(\Sigma)
                =\overline\Hom(\Sigma,\Sigma)
                \simeq\Aut(\Sigma)/W_{\Sigma}\\
             &\overset{\sim}{\leftarrow}\{g\in\Aut(\Sigma)\,;\,g(\Psi)=\Psi\}.\label{eq:OAut}
\end{align}

The root system $\Xi$ is \textsl{isomorphic to} $\Sigma$,
which is denoted by $\Xi\simeq\Sigma$,
if there exists a surjective homomorphism of $\Xi$ onto $\Sigma$.

Suppose $\Sigma_1$ and $\Sigma_2$ are the subsystems of $\Sigma$ 
such that $\Sigma=\Sigma_1\cup\Sigma_2$ and $\Sigma_1\perp\Sigma_2$.
Then we say that $\Sigma$ is a \textsl{direct sum} of $\Sigma_1$ and 
$\Sigma_2$, which is denoted by $\Sigma=\Sigma_1+\Sigma_2$.
The \textsl{irreducible root system} is a root system 
which has no non-trivial direct sum decomposition.
Note that every root system is decomposed into 
a direct sum of irreducible root systems and
\begin{equation}
 \Aut(\Sigma)\simeq\{g\in O(n)\,;\,g(\Sigma)=\Sigma\}
\end{equation}
if $\Sigma$ is an irreducible root system of rank $n$.
Here $O(n)$ is the orthogonal group of $\mathbb R^n$ with respect to $(\ |\ )$.

For root systems $\Sigma_1$ and $\Sigma_2$ there exists a root system
$\Sigma=\Sigma'_1+\Sigma'_2$ such that $\Sigma_j\simeq\Sigma'_j$
for $j=1$ and $2$. 
This root system $\Sigma$ is determined modulo isomorphisms and 
hence we simply write $\Sigma=\Sigma_1+\Sigma_2$.
When $\Sigma_1=\Sigma_2$, we sometimes write $2\Sigma_1$ in place of
$\Sigma_1+\Sigma_2$.
\end{defn}

For any two elements $\alpha$ and $\alpha'$ in $\Psi$, 
there exists an isomorphism $\iota$ of $\langle\alpha,\alpha'\rangle$ to 
one of the following four root systems with the fundamental system 
$\{\beta,\beta'\}$ such that
$\iota(\alpha)=\beta$ and $\iota(\alpha')=\beta'$:

\noindent
$A_1+A_1=2A_1$:
$(\beta,\beta')=(\epsilon_1,\epsilon_2)$
\hfill
$2\frac{(\beta|\beta')}{(\beta|\beta)} =0,\phantom{-}\quad
 2\frac{(\beta|\beta')}{(\beta'|\beta')}=0\phantom{-}$
\quad
\begin{xy}
  *+!D{\beta} *\cir<3pt>{} ;
  (10,0)  *+!D{\beta'} *\cir<3pt>{}
\end{xy}

\noindent
$A_2$: 
$(\beta,\beta')=(\epsilon_1-\epsilon_2,\epsilon_2-\epsilon_3)$
\hfill
$2\frac{(\beta|\beta')}{(\beta|\beta)}=-1,\quad
 2\frac{(\beta|\beta')}{(\beta'|\beta')}=-1$
\quad
\begin{xy}
 \ar@{-}  *+!D{\beta} *\cir<3pt>{} ;
  (10,0)  *+!D{\beta'} *\cir<3pt>{}
\end{xy}

\noindent
$B_2$: 
$(\beta,\beta')=(\epsilon_1-\epsilon_2,\epsilon_2)$
\hfill
$2\frac{(\beta|\beta')}{(\beta|\beta)}=-1,\quad
 2\frac{(\beta|\beta')}{(\beta'|\beta')}=-2$
\quad
\begin{xy}
 \ar@{=>}  *+!D{\beta} *\cir<3pt>{} ;
  (10,0)  *+!D{\beta'} *\cir<3pt>{}
\end{xy}

\noindent
$G_2$: 
$(\beta,\beta')=(-2\epsilon_1+\epsilon_2+\epsilon_3,\epsilon_1-\epsilon_2)$
\hfill
$2\frac{(\beta|\beta')}{(\beta|\beta)}=-1,\quad
 2\frac{(\beta|\beta')}{(\beta'|\beta')}=-3$
\quad
\begin{xy}
 \ar@3{->}  *+!D{\beta} *\cir<3pt>{} ;
  (10,0)  *+!D{\beta'} *\cir<3pt>{}
\end{xy}

\smallskip
The \textsl{Dynkin diagram} $G(\Psi)$ of a root system $\Sigma$ with the 
fundamental system $\Psi$ is the graph which consists of vertices 
expressed by small circles and edges expressed by some lines or arrows
such that the vertices are associated to the elements of $\Psi$.
The lines or arrows connecting two vertices represent 
the isomorphic classes of the corresponding two roots 
in $\Psi$ according to the above Dynkin diagram of rank 2.
Here the number of lines which link $\beta$ to $\beta'$ in the diagram equals
$\frac{-2(\beta|\beta')}{\min\{(\beta|\beta),(\beta'|\beta')\}}$.
The arrow points toward a shorter root.

\begin{defn}\label{def:extended}
A root $\alpha$ of an irreducible root system $\Sigma$ is called
\textsl{maximal} and denoted by $\alpha_{max}$ if every number 
$m_j(\alpha)$ for $j=1,\ldots,n$ in Definition~\ref{defn:fund} is
 maximal among the roots of $\Sigma$.
It is known that the maximal root uniquely exists.

Let $\Psi=\{\alpha_1,\ldots,\alpha_n\}$ be a fundamental system of $\Sigma$.
Define
\begin{align}
 \alpha_0 &:= %\begin{cases}
                -\alpha_{max},%&(n>1)
\\
              %  \alpha_1&(n=1)
%              \end{cases},\\
 \tilde\Psi&:=\Psi\cup\{\alpha_0\}.\label{eq:extFund}
\end{align}
The \textsl{extended Dynkin diagram} of $\Sigma$ in this note is 
the graph $G(\tilde\Psi)$ associated to $\tilde\Psi$ which is defined in 
the same way as $G(\Psi)$ associated to $\Psi$.
We call $\tilde\Psi$ the \textsl{extended fundamental system} of $\Sigma$.
A \textsl{subdiagram} of $G(\tilde\Psi)$ is the Dynkin diagram $G(\Theta)$
associated to a certain subset $\Theta\subset\tilde\Psi$.
\end{defn}
In \S\ref{sec:Dynkin} the extended Dynkin diagrams of all 
the irreducible root systems is listed, which are based on the notation in
\cite{Bo}.
The vertex expressed by a circled circle in the diagram corresponds to 
the special root $\alpha_0$. If the vertex and the lines starting from 
it are removed from the diagram, 
we get the corresponding Dynkin diagram of the irreducible root system.
The numbers below vertices $\alpha_j$ in the diagram in 
\S\ref{sec:Dynkin} are the numbers $m_j(\alpha_{max})$ given 
by \eqref{eq:coeff}. We define $m_0(\alpha_{\max})=1$ and then
\begin{equation}\label{eq:zerosum}
\sum_{\alpha_j\in\Tilde\Psi} m_j(\alpha_{\max})\alpha_j=0.
\end{equation}
%%%
\begin{rem} \label{rem:graphauto}
{\rm i)}
There is a bijection between the isomorphic classes of root systems
and the Dynkin diagrams.

The irreducible decomposition of a root system $\Sigma$ corresponds to
the decomposition of its Dynkin diagram $G(\Psi)$ into the connected 
components $G(\Psi_j)$.
It also induces the decomposition of the fundamental system
$\Psi = \Psi_1\amalg\cdots\amalg \Psi_m$ such that
$\Sigma=\langle\Psi_1\rangle+\cdots+\langle\Psi_m\rangle$ is the
decomposition into irreducible root systems.
Then we call each $\Psi_j$ an \textsl{irreducible component} 
of $\Psi$.

The irreducible root systems are classified as follows 
(cf.~\S\ref{sec:Dynkin}):
\begin{equation}
  A_n (n\ge1),\ B_n(n\ge2),\ C_n(n\ge3),\ 
 D_n(n\ge 4),\ E_6,\ E_7,\ E_8,\ F_4,\ G_2.
\end{equation}
We will also use this notation $A_n,\ldots$ for a root system or a 
fundamental system.  For example, $A_2+2B_3$ means a root system 
isomorphic to the direct sum of the root system of type $A_2$ and 
two copies of the root system of type $B_3$ or it means its 
fundamental system.

{\rm ii)} $\Out(\Sigma)$ is naturally isomorphic to the group of graph 
automorphisms of the Dynkin diagram associated to $\Sigma$.
If $\Sigma$ is irreducible, it also corresponds to the graph automorphisms of 
the extended Dynkin diagram which fix the vertex corresponding to $\alpha_0$.
Here we give the list of irreducible root systems $\Sigma$ with non-trivial
$\Out(\Sigma)$:
\begin{equation}
 \begin{cases}
   \Out(A_n) \simeq \mathbb Z/2\mathbb Z&(n\ge 2),\\
   \Out(D_4) \simeq \mathfrak S_ 3\, 
     &(=\text{the symmetric group of degree 3}),\\
   \Out(D_n) \simeq \mathbb Z/2\mathbb Z&(n\ge 5),\\
   \Out(E_6) \,\simeq \mathbb Z/2\mathbb Z.
 \end{cases}
\end{equation}

{\rm iii)} The graph automorphism $\sigma$ of 
the extended Dynkin diagram $G(\tilde \Psi)$ with the following property
corresponds to a transformation by an element of $W_\Sigma$.
\begin{equation}
 \begin{cases}
   \text{A rotation of }G(\tilde\Psi) &(\Sigma=A_n,\ E_6),\\
   \text{Any automorphism} &(\Sigma=B_n,\ C_n,\ E_7),\\
   \sigma\bigl((\alpha_0,\alpha_1,\alpha_{n-1},\alpha_n)\bigr)
   =(\alpha_1,\alpha_0,\alpha_n,\alpha_{n-1}) &(\Sigma=D_n),\\
   \sigma\bigl((\alpha_0,\alpha_1,\alpha_{n-1},\alpha_n)\bigr)
   =(\alpha_n,\alpha_{n-1},\alpha_1,\alpha_0) 
        &(\Sigma=D_n,\ n:\text{even}\ge4),\\
   \sigma\bigl((\alpha_0,\alpha_1,\alpha_{n-1},\alpha_n)\bigr)
   =(\alpha_n,\alpha_{n-1},\alpha_0,\alpha_1) 
        &(\Sigma=D_n,\ n:\text{odd}).
 \end{cases}
\end{equation}
When $\Sigma$ is irreducible, we have the bijection:
\begin{equation}\label{eq:mult1}
\begin{matrix}
 \{w\in W_\Sigma\,;\,w(\tilde\Psi)=\tilde\Psi\} &\overset{\sim}{\to}
 &\{\alpha_j\in\tilde\Psi\,;\,m_j(\alpha_{max})=1\}&\\
  \vin & & \vin\\
 \sigma &\mapsto &\sigma(\alpha_0)
\end{matrix}
\end{equation}
\end{rem}
To classify subsystems contained in a root system we prepare
more definitions.
\begin{defn}\label{def:weq}
We put
\begin{align}
 \Aut'(\Xi)&:=\Aut(\Xi_1)\times\cdots\times\Aut(\Xi_m)
   \subset\Aut(\Xi),\\
 \Out'(\Xi)&:=\Aut'(\Xi)/W_\Xi
\end{align}
for a root system $\Xi$ with an irreducible decomposition
$\Xi=\Xi_1+\cdots+\Xi_m$.
\end{defn}

\begin{defn}
Let $\Xi$, $\Xi'$ and $\Theta$ be subsystems of $\Sigma$.
\begin{align}
 \Xi\underset{\Theta}\sim\Xi'\ &\Leftrightarrow\ 
 \exists w\in W_\Theta\text{ such that }
 \Xi'=w(\Xi),\\
 \Xi\underset{\Theta}{\overset{w}\sim}\Xi'\ &\Leftrightarrow\ 
 \exists g\in\Aut(\Theta)\text{ such that }
 \Xi'=g(\Xi).
\end{align}
If $\Xi\underset{\Theta}\sim\Xi'$
(resp.~$\Xi\underset{\Theta}{\overset{w}\sim}\Xi'$),
we say that $\Xi'$ is \textsl{equivalent }
(resp.~\textsl{weakly equivalent}) \textsl{to $\Xi$ by $\Theta$}.
Since $\Aut(\Xi)\simeq\{\iota\in\Hom(\Xi,\Sigma)\,;\,\iota(\Xi)=\Xi\}$,
we have
\begin{align}
 \begin{split}
 &\{\Xi'\subset \Sigma\,;\, 
    s_\alpha(\Xi')=\Xi'\ (\forall\alpha\in\Xi')\text{ and }
  \Xi'\simeq\Xi\}/\!\underset{\Sigma}{\sim}\\
 &\qquad\simeq W_{\Sigma}\backslash\!\Hom(\Xi,\Sigma)/\!\Aut(\Xi)\\
 &\qquad\simeq \overline\Hom(\Xi,\Sigma)/\!\Out(\Xi),
 \end{split}\\
 \begin{split}
 &\{\Xi'\subset \Sigma\,;\,  
    s_\alpha(\Xi')=\Xi'\ (\forall\alpha\in\Xi')\text{ and }
  \Xi'\simeq\Xi\}/\!\underset{\Sigma}{\overset{w}\sim}\\
 &\qquad\simeq\Aut(\Sigma)\backslash\!\Hom(\Xi,\Sigma)/\!\Aut(\Xi)\\
 &\qquad\simeq \Out(\Sigma)\backslash\overline\Hom(\Xi,\Sigma)/\!\Out(\Xi),
 \end{split}\\
 \begin{split}
 &W_{\Sigma}\backslash\!\Hom(\Xi,\Sigma)/\!\Aut'(\Xi)\simeq
  \overline\Hom(\Xi,\Sigma)/\!\Out'(\Xi).
 \end{split}
\end{align}
\end{defn}
%%%%%%%%%%
\begin{defn}[fundamental subsystems]
A subsystem $\Xi$ of $\Sigma$ is called \textsl{fundamental} if 
there exists $\Theta\subset\Psi$ such that
$\Xi\underset{\Sigma}{\sim}\langle\Theta\rangle$.
\end{defn}
\begin{rem}
Suppose $\Sigma$ is of type $A_n$.
Then it is clear that
%For subsystems $\Xi$ and $\Xi'$ of $\Sigma$ it is clear that
\begin{align}
&\text{any subsystem of $\Sigma$ is fundamental},\\
&\bigl(\Xi\underset{\Sigma}{\sim}\Xi'
\ \Leftrightarrow\ 
   \Xi\simeq \Xi'\bigr)
 \text{ \ for subsystems $\Xi$ and $\Xi'$ of $\Sigma$}.
\end{align}
\end{rem}
%%%%%%%%%%
Our aim in this note is to clarify the structure of
\begin{gather*}
 \overline\Hom(\Xi,\Sigma),\ 
 \Out(\Sigma)\backslash\overline\Hom(\Xi,\Sigma),\ 
 \overline\Hom(\Xi,\Sigma)/\!\Out(\Xi),\ 
  \overline\Hom(\Xi,\Sigma)/\!\Out'(\Xi),\\ 
 \Out(\Sigma)\backslash\overline\Hom(\Xi,\Sigma)/\!\Out(\Xi)
 \text{ and fundamental subsystems of $\Sigma$}.
\end{gather*}
For this purpose we prepare the following definition.
\begin{defn}\label{def:endroot}
 {\rm i)}
A root $\alpha\in\Psi$ (resp.\ $\tilde\Psi$) is an \textsl{end root} 
of $\Psi$ (resp.\ $\tilde\Psi$) if 
\begin{equation}
\#\{\beta\in\Psi\text{ (resp.~$\tilde\Psi$)}\,;\,(\beta,\alpha)<0\}\le 1.
\end{equation}
A root $\alpha\in\Psi$ (resp.~$\tilde\Psi$) is called 
a \textsl{branching root} of $\Psi$ (resp.~$\tilde\Psi$) if
\begin{equation}
 \#\{\beta\in\Psi\text{ (resp.~$\tilde\Psi$)}\,;\,(\beta,\alpha)<0\}
 \ge 3.
\end{equation}
The corresponding vertex in the (extended) Dynkin diagram is also
called an end vertex or a branching vertex, respectively.

{\rm ii)}
When $\Sigma$ is irreducible, we put
\begin{equation}
  \Sigma^L := \{\alpha\in\Sigma\,;\,
   |\alpha|=|\alpha_{max}|\}
\end{equation}
and denote its fundamental system by $\tilde\Psi^L$.
Then $\Sigma^L$ is a subsystem of $\Sigma$ and
\begin{equation}\label{eq:SL}
 \begin{gathered}
  A_n^L = A_n,\ B_n^L = D_n\ (n\ge2),\ C_n^L = nA_1\ (n\ge 3),\ D_n^L = D_n\ (n\ge 4),\\
  E_6^L=E_6,\ E_7^L = E_7,\ E_8^L = E_8,\ F_4^L = D_4,\ G_2^L = A_2.
 \end{gathered}
\end{equation}
A root system whose Dynkin diagram contains no arrow is called
\textsl{simply laced}. % if any irreducible component $\Xi$ satisfies $\Xi^L=\Xi$.
\end{defn}
%%%%%%%%%%%%%%%%%%%%%%%%%%%%%%%%%%%%%%%%%%%%%%%%%%%
%%%%%%%%%%%%%%%%%% Theorem %%%%%%%%%%%%%%%%%%%%%%%%
\section{A theorem}\label{sec:thm}
In this section we will give a simple procedure to clarify
the set
$ %\begin{equation}
 \overline\Hom(\Xi,\Sigma):=W_{\Sigma}\backslash\!\Hom(\Xi,\Sigma)
$ %\end{equation}
for root systems $\Xi$ and $\Sigma$.
\begin{rem}\label{rem:red}
{\rm i) } Note that
\begin{align}
\overline{\Hom}\bigl(\Xi,\Sigma_1+\Sigma_2\bigr)
 &\simeq \coprod_{\Xi'\subset\Xi:\ \text{component}}\Bigl(
      \overline{\Hom}\bigl(\Xi',\Sigma_1\bigr),
      \overline{\Hom}\bigl((\Xi')^\perp,\Sigma_2\bigr)
    \Bigr),\label{eq:SigDiv}\\
\overline{\Hom}\bigl(\Xi_1+\Xi_2,\Sigma\bigr)
 &\simeq \coprod_{\bar\iota\in \overline{\Hom}(\Xi_1,\Sigma)}
   \Bigl(\bar\iota, \overline{\Hom}\bigl(\Xi_2,\iota(\Xi_1)^\perp\bigr)
   \Bigr).\label{eq:XiDiv}
\end{align}
Here $\iota$ is a representative of $\bar\iota$ and 
the \textsl{component} $\Xi'$ of $\Xi$ is the subsystem of $\Xi$ such 
that $\Xi=\Xi'+(\Xi')^\perp$.  
The empty set and $\Xi$ are also components of $\Xi$.

The identification \eqref{eq:XiDiv} follows from
\begin{equation}
 \{w\in W_\Sigma\,;\,w|_{\iota(\Xi)}=id\} = 
  W_{\iota(\Xi)^\perp}\subset W_\Sigma
\end{equation}
for any $\iota\in\Hom(\Xi,\Sigma)$ (cf.~\cite{Bo}).

{\rm ii) }
The identifications \eqref{eq:SigDiv} and \eqref{eq:XiDiv} assure
that we may assume $\Xi$ and $\Sigma$ are irreducible.
In fact, the study of the structure of $\overline{\Hom}(\Xi,\Sigma)$ 
is reduced to the study of $\bar\iota\in\overline{\Hom}(\Xi,\Sigma)$ 
and $\iota(\Xi)^\perp$ for irreducible $\Xi$ and $\Sigma$.

{\rm iii) }
We may moreover assume
$\iota(\Xi)\cap\Sigma^L\ne\emptyset$ 
by considering the dual root systems 
   $\Xi^{\vee}:=\bigl\{\frac{2\alpha}{(\alpha|\alpha)}\,;\,\alpha\in\Xi\bigr\}$
     and
   $\Sigma^{\vee}:=\bigl\{\frac{2\alpha}{(\alpha|\alpha)}\,;\,\alpha\in\Sigma\bigr\}$
in place of $\Xi$ and $\Sigma$, respectively.
\end{rem}

\begin{defn}
When $G(\Phi)$ is isomorphic to a subdiagram $G(\Theta)$ of 
$G(\tilde \Psi)$ with a map $\bar\iota:\Phi\to\Theta\subset \tilde\Psi$,
it is clear that $\bar\iota$ defines an element
of $\Hom(\Xi,\Sigma)$.  
In this case we say that $\bar\iota$ is an \textsl{imbedding} of 
$G(\Phi)$ into $G(\tilde\Phi)$.
\end{defn}

Recalling Definition~\ref{def:perp}, \ref{def:extended} and 
\ref{def:endroot}, 
we now state a main lemma in this note.
%, which is a generalization 
%of the description of $\Out(\Sigma)=\overline\Hom(\Sigma,\Sigma)$ 
%in Remark~\ref{rem:graphauto} ii).
%%%%%%%%%%%%%%%%%  Theorem %%%%%%
\begin{lem}\label{lem:thm}
Let $\Xi$ and $\Sigma$ be irreducible root systems and
let $\Phi$ and $\Psi$ be their fundamental systems, respectively.
Denoting
\begin{align}
  \Hom'(\Xi,\Sigma)&:=\bigl\{\iota\in \Hom(\Xi,\Sigma)\,;\,\iota(\Xi)\cap
  \Sigma^L\ne\emptyset\bigr\},\\
 \overline\Hom'(\Xi,\Sigma)&:=W_\Sigma\backslash\!\Hom'(\Xi,\Sigma),
\end{align}
 we have the following claims according to the type of $\Xi$:
%%%%%

{\rm 1)} \ $\Xi$ is of type $A_m$. 
\begin{gather*}
\overline\Hom'(\Xi,\Sigma) %\\
\overset{\sim}{\leftarrow}
\bigl\{\text{Imbeddings $\bar\iota$ of $G(\Phi)$ into $G(\tilde\Psi)$
with the end vertex $\alpha_0$%
}\bigr\}.
\end{gather*}
Let $\bar\iota$ be this graph imbedding corresponding to 
$\iota\in\Hom(\Xi,\Sigma)$.
Then
\begin{equation}\label{eq:orthA}
 \iota(\Xi)^\perp
 \simeq \bigl\langle\alpha\in\tilde\Psi\,;\,
   \alpha\perp\bar\iota(\Phi)\bigr\rangle.
\end{equation}
In the case $\#\overline\Hom'(\Xi,\Sigma)>1$, we have
$\#\overline\Hom'(\Xi,\Sigma)=3$ if $(\Xi,\Sigma)$ is of type
$(A_3,D_4)$ and $2$ if otherwise.
Moreover for 
$\bar\iota,\,\bar\iota'\in\overline\Hom'(\Xi,\Sigma)$
\begin{equation}
 \begin{gathered}
  \text{``\,$\bar\iota$ and $\bar\iota'$ are conjugate under 
  an element of $\Out(\Sigma)$ or $\Out(\Xi)$"}\\
  \text{$\Leftrightarrow$
  $\ \iota(\Xi)^\perp\simeq\iota'(\Xi)^\perp$}.
 \end{gathered}
\end{equation}
%%%%%%%%%

{\rm 2)} \ 
$\Xi$ is of type $D_m$ {\rm ($m\ge 4$)}.\\
Let $\Phi_m=\{\beta_0,\ldots,\beta_{m-1}\}$ be 
a fundamental system of $\Xi$ with the Dynkin diagram
\begin{xy}
 \ar@{-}  *+!D{\beta_0} *{\circ} ;
  (5,0)  *+!D{\beta_1} *{\circ}="B"
 \ar@{-} "B";(10,0) *+!D{\beta_2} *{\circ}="C"
% \ar@{-} "C";(15,0) *+!D{\beta_3} *{\bullet}="D"
% \ar@{-} "D";(20,0) *+!D{\alpha_5} *{\bullet}="E"
 \ar@{-} "C";(14,0) \ar@{.} (14,0);(19,0)^*!U{\!\!\!\cdots}
 \ar@{-} (19,0);(23,0) *+<5pt>!D{\beta_{m-3}} *{\circ}="H"
 \ar@{-} "H";(28,0) *+<5pt>!DL{\!\beta_{m-2}} *{\circ}
 \ar@{-} "H";(23,-5) *+!L{\beta_{m-1}} *{\circ}
\end{xy}
and $m_\Sigma$ denote the 
maximal integer $m$ such that there is an imbedding 
$\bar\iota_m$ of $G(\Phi_m)$ into $G(\tilde\Psi^L)$. 
We put $m_\Sigma=0$ if such imbedding doesn't exist.
Then
\[
  m_\Sigma =
 \begin{cases}
   0  &(\Sigma\text{ is of type } A_n,\ C_n,\ G_2),\\
   \rank \Sigma &(\Sigma\text{ is of type } B_n,\ D_n,\ E_8,\ F_4),\\
   5  &(\Sigma\text{ is of type } E_6),\\
   6  &(\Sigma\text{ is of type } E_7)
 \end{cases}
\]
and
\begin{align*}
  \Hom'(D_m,\Sigma)\ne\emptyset\ 
  &\Leftrightarrow\ 
   (4\le)m\le m_\Sigma\\
  &\Leftrightarrow\ 
   \#\bigl(\overline\Hom'(\Xi,\Sigma)/\!\Out(\Xi)\bigr) = 1.
\end{align*}

\underline{$\Sigma$ is of type $E_6$, $E_7$ or $E_8$}.
\begin{align*}
 &\#\overline\Hom(D_m,\Sigma)=
 \begin{cases}
    2 &(m=m_\Sigma),\\
    1 &(4\le m<m_\Sigma),
  \end{cases}\\
 &\iota(\Xi)^\perp\simeq
     D_{m_\Sigma-m}+
  \begin{cases}
      A_1 &(n=7), \\
      \emptyset &(n=6,\,8).
  \end{cases}
\end{align*}

%%%%
\underline{$\Sigma$ is of type $D_n$, $B_n$ or $F_4$ \rm ($m\le n$)}.
\begin{align*}
 \#\overline\Hom'(\Xi,\Sigma)=
 &\begin{cases}
  6 &(\Sigma:D_4\ (m=n=4)),\\
  3 &(\Sigma:B_n\text{ and }D_n\ (m=4<n)),\\
  2 &(\Sigma:D_n\ (4<m=n)),\\
  1 &(\Sigma:F_4\ (4=m),\ B_n\ (4 < m \le n),\  D_n\ (4<m<n)),
\end{cases}\\
 \iota(\Xi)^\perp\simeq
  &\begin{cases}
  D_{n-m} & (\iota\in\Hom(D_m,D_n)),\\
  B_{n-m} & (\iota\in\Hom(D_m,B_n)),\\
  \emptyset  & (\iota\in\Hom(D_4,F_4)).
 \end{cases}
\end{align*}

{\rm 3)} \ $\Xi$ is of type $B_m$ {\rm ($m\ge 2$)}.
\begin{gather*}
\Hom(\Xi,\Sigma)\ne\emptyset\ \Leftrightarrow\ 
 \#\overline\Hom(\Xi,\Sigma)=1\text{ and }
 \begin{cases}
  \Sigma \text{ is of type } B_n \text{ with } m\le n,\\
  \Sigma \text{ is of type } C_n \text{ with } m=2,\\
  \Sigma \text{ is of type } F_4 \text{ with } m\le 4,
 \end{cases}\\
\iota(\Xi)^\perp\cap T_n
 \simeq T_{n-m}\quad (T=B,\,C,\,F,\ F_2=B_2,\,
  \,F_1=A_1\text{ and }\iota(B_3)^\perp\cap F_4\not\subset \Sigma^L).
\end{gather*}
%%%

{\rm 4)} \ 
$\Xi$ is of type $C_m$ {\rm ($m\ge3$)}.
\begin{gather*}
 \Hom(\Xi,\Sigma)\ne\emptyset\ \Leftrightarrow\ 
 \#\overline\Hom(\Xi,\Sigma)=1\text{ and }
 \begin{cases}
  \Sigma \text{ is of type } C_n \text{ with }m\le n,\\
  \Sigma \text{ is of type } F_4 \text{ with } m\le 4,
\end{cases}\\
\iota(\Xi)^\perp\cap T_n
 \simeq T_{n-m}\quad (T=C,\,F,\ F_2=C_2,\,F_1=A_1\text{ and }\iota(C_3)^\perp\cap F_4\subset \Sigma^L).
\end{gather*}

%%%

{\rm 5)} \ 
$\Xi$ is of type $E_m$ {\rm($m=6$, $7$ and $8$)}.
\begin{align*}
\overline\Hom(\Xi,\Sigma)&\overset{\sim}{\leftarrow}
\bigl\{\text{Imbeddings $\bar\iota$ of $G(\Phi)$ into $G(\tilde\Psi)$}\bigr\}\big/\!\!\sim,\\
 \iota(\Xi)^\perp&
 \simeq\bigl\langle\alpha\in\tilde\Psi\,;\, \alpha\perp\bar\iota(\Phi)\bigr\rangle.
\end{align*}
Here $/\!\!\sim$ is interpreted that all the imbeddings of $G(\Phi)$
are considered to be isomorphic except for $(\Xi,\Sigma)\simeq (E_6,E_6)$.
Namely
$\#\overline\Hom(\Xi,\Sigma)\le 1$ if $(\Xi,\Sigma)\not\simeq
(E_6,E_6)$.

{\rm 6)} \ \
$\Xi$ is of type $G_2$ or $F_4$.
\[
 \Hom(\Xi,\Sigma)\ne\emptyset\ \Leftrightarrow\ 
 \#\overline\Hom(\Xi,\Sigma)=1\text{ and }\Xi\simeq\Sigma.
\]
\end{lem}
%%%
\begin{rem}\label{rem:thmauto}
{\rm i)}
In the proof of Lemma~\ref{lem:thm}~2) we will have
\[
 \iota(\Xi)^\perp\simeq
  \begin{cases}
   \bigl\langle
     \bar\iota_{m_\Sigma}(\Phi_{m_\Sigma})^\perp\cap\Psi
   \bigr\rangle
     & (m_\Sigma-1\le m \le m_\Sigma),\\
   \bigl\langle
     \bar\iota_{m_\Sigma}(\Phi_{m_\Sigma})^\perp\cap\Psi,\ 
     \bar\iota_{m_\Sigma}(\beta_m),\ldots
     \bar\iota_{m_\Sigma}(\beta_{m_\Sigma-1})
   \bigr\rangle
     &(4\le m\le m_\Sigma-2)
  \end{cases}
\]
for the imbedding $\bar\iota_{m_\Sigma}$ with 
$\bar\iota(\beta_{m_\Sigma})=\alpha_0$ if $\Sigma$ is of type 
$D_n$, $E_6$, $E_7$ or $E_8$.

Let $\Theta_m$ be a subset of $\tilde\Psi$ such that
$\langle\Theta_m\rangle\simeq D_m$.
If $\Sigma$ is of type $B_n$ or $D_n$, we may assume that
$\iota_m\in\overline\Hom(D_m,\Sigma)$ satisfies
$\iota_m(\Xi)=\langle\Theta_m\rangle$ and then
\begin{gather}
  \iota_m(\Phi_m)^\perp = 
   \bigl\langle\Theta_m^\perp\cap\tilde\Psi\bigr\rangle.
\end{gather}
Suppose $\Sigma$ is of type $E_6$, $E_7$ or $E_8$.
Let $\tilde\alpha_{max}$ be the maximal root of 
$\langle \Theta_{m_\Sigma}\rangle$. 
Put $\tilde\alpha_0=-\tilde\alpha_{max}$ and 
$\tilde\Theta_{m_\Sigma}=\Theta_{m_\Sigma}\cup\{\tilde\alpha_0\}$.
We may assume 
$\iota_m\in\overline\Hom(D_m,\Sigma)$ satisfies
$\iota_m(\Xi)=\langle\Theta_m\rangle$ and 
$\Theta_m\subset\tilde\Theta_{m_\Sigma}$.
Then
\begin{gather}
  \iota_m(\Phi_m)^\perp = 
    \bigl\langle\Theta_m^\perp\cap\tilde\Theta_{m_\Sigma},\,
  \Theta_{m_\Sigma}^\perp\cap\tilde\Psi\bigr\rangle.
\end{gather}
Note that $G(\tilde\Theta_{m_\Sigma})$ is the extended Dynkin diagram of
$\langle\Theta_{m_\Sigma}\rangle\simeq D_{m_\Sigma}$.
See Example~\ref{exmp:A} viii) and ix).

{\rm ii)}
Using a graph automorphism of $G(\tilde\Psi)$ corresponding to
a suitable element of $W_\Sigma$, we may replace $\alpha_0$ by another 
element $\alpha_j$ of $\Psi$ with $m_j(\alpha_{max})=1$
in Theorem~\ref{thm} and in the remark above
(cf.~Remark~\ref{rem:graphauto} ii)).

{\rm iii)}
The image $\iota(\Xi)$ corresponding to the graph automorphism 
$\bar\iota$ in Lemma~\ref{lem:thm} is obtained by Proposition~\ref{prop:subsystem}.
\iffalse
{\rm iv)}
Suppose $\Sigma$ is of type $R$ with
$R=A_n$, $B_n$, $C_n$ or $D_n$.
Then for any imbedding $\Xi\subset\Sigma$ of the irreducible root system 
$\Xi$, we can find a graphic imbedding of $G(\Phi)$ into a subdiagram
$G_\Phi$ of $\tilde R$ or $\tilde R'$ containing $\alpha_0$
such the subsystem $\Sigma_\Phi$ corresponding to $G_\Phi$ is equivalent 
to $\Xi$ by $\Sigma$.
Note that $\Sigma_\Phi^\perp$
%corresponding to $G_\Phi$ 
is generated by the roots in $\tilde R$ or 
$\tilde R'$ orthogonal to $G_\Phi$.
Here $\tilde R=G(\Psi)$ and $\tilde R'$ is an affine Dynkin diagram
given in Proposition~\ref{prop:clasify} (cf.~Remark~\ref{rem:cls} iii)) if
$R$ is not simply laced
and we put $\tilde R'=\tilde R$ if otherwise.
%false
If $\Xi$ is not of type $D_4$, we have
\begin{equation}\begin{split}
 \overline\Hom(\Xi,\Sigma)\simeq
 \{
 &\text{Imbeddings $\bar\iota$ of $G(\Phi)$ to $G(\tilde\Psi)$ or 
 $G(\tilde\Psi')$
 such that}\\
 &\text{a certain fixed end vertex of $G(\Phi)$ corresponds to
 $\alpha_0$ or $\alpha_0'$}
 \}.
\end{split}\end{equation}
\fi
\end{rem}
%%%
Lemma~\ref{lem:thm} can be summarized in the following form.
\begin{thm}\label{thm}
Let $\Sigma$ and $\Xi$ be irreducible root systems and
let $\Psi$ and $\Phi$ be their fundamental systems, respectively.
Retain the notation given in 
Definition~\ref{def:perp}--\ref{def:extended} and \ref{def:endroot}.
If $\Sigma$ is not simply laced, 
we denote the maximal root in $\Sigma\setminus\Sigma^L$ by 
$\alpha_{max}'$ and the Dynkin diagram of $\tilde\Psi'$ by
$G(\tilde\Psi')$.  Here we put $\alpha_0'=-\alpha_{max}'$ and
$\tilde\Psi'=\Psi\cup\{\alpha_0'\}$.

{\rm i)}
Suppose $\Sigma$ is of the classical type or $\Xi\simeq A_m$ 
with $m\ge1$.

When $\Xi\not\simeq D_4$ or $(\Sigma,\Xi)\simeq(D_4,D_4)$,
\begin{equation}\label{eq:GA}
\begin{split}
 \overline\Hom(\Xi,\Sigma)\overset{\sim}{\leftarrow}
 \{
 &\text{Imbeddings $\bar\iota$ of $G(\Phi)$ to $G(\tilde\Psi)$ or 
 $G(\tilde\Psi')$}\\
 &\text{such that $\beta_0$ corresponds to
 $\alpha_0$ or $\alpha_0'$ by $\bar\iota$}
 \}
\end{split}\end{equation}
for a suitable root $\beta_0\in\Phi$.
Here we delete $G(\tilde\Psi')$ and $\alpha_0'$ in the above
if $\Sigma$ is simply laced.
Moreover $\beta_0$ is any root in $\Phi$ such that the 
right hand side of \eqref{eq:GA} is not empty and 
if such $\beta_0$ doesn't exit, $\overline{\Hom}(\Xi,\Sigma)=\emptyset$.

When $\Xi\simeq D_4$, 
\begin{equation}\label{eq:Out1}
  \#\bigl(\overline\Hom(\Xi,\Sigma)/\!\Out(\Xi)\bigr)\le1
\end{equation}
and the representative of %the coset 
$\overline{\Hom}(\Xi,\Sigma)/\!\Out(\Xi)$
is given by the above imbedding $\bar\iota$ and
\begin{equation}
 \#\overline{\Hom}(D_4,B_n)=\#\overline{\Hom}(D_4,C_n)=
 \#\overline{\Hom}(D_4,D_{n+1})=3\quad(n\ge4).
\end{equation}

For $\iota\in\Hom(\Xi,\Sigma)$ corresponding to this 
imbedding $\bar\iota$ of $G(\Phi)$ we have
\begin{equation}\label{eq:orthA0}
 \iota(\Xi)^\perp
 = \bigl\langle\alpha\in\Psi\,;\,
   \alpha\perp\bar\iota(\Phi)\bigr\rangle.
\end{equation}
Moreover for 
$\bar\iota,\,\bar\iota'\in\overline\Hom(\Xi,\Sigma)$
\begin{equation}
 \begin{gathered}
  \text{``\,$\bar\iota$ and $\bar\iota'$ are conjugate under 
  an element of $\Out(\Sigma)$ or $\Out(\Xi)$"}\\
  \text{$\Leftrightarrow$
  $\ \iota(\Xi)\cap\Sigma^L\simeq\iota'(\Xi)\cap\Sigma^L$ and 
  $\iota(\Xi)^\perp\simeq\iota'(\Xi)^\perp$}.
 \end{gathered}
\end{equation}

{\rm ii)}
Suppose $\Sigma$ is of the exceptional type and 
$\Xi=R_m$ with $R=B$, $C$, $D$, $E$, $F$ and $G$.
Put
$m_0^R=2$, $3$, $4$, $6$, $4$ and $2$ according to $R=B$, $C$, $D$, $E$, 
$F$ and $G$, respectively, and moreover suppose $m\ge m_0^R$.
%Then we have \eqref{eq:Out1}.
%
Let $m_\Sigma^R$ be the maximal number $m$ such that
the Dynkin diagram $G(R_m)$ of the root system $R_m$ is
a subdiagram of $G(\tilde\Psi)$ or $G(\tilde\Psi')$.
Thus for a subset $\Phi_\Sigma^R$ of $\tilde\Psi$ or $\tilde\Psi'$
we identify $G(R_{m_\Sigma^R})$ with the subdiagram
$G(\Phi_\Sigma^R)$.
Put $m_\Sigma^R=0$ if such number $m$ with $m\ge m_0^R$ does not exists.

When $(\Sigma,R_m)\not\simeq(F_4,D_4)$, we have \eqref{eq:Out1} and
\begin{gather}
 \#\overline\Hom(R_m, \Sigma) =
  \begin{cases}
   0 &(m > m_\Sigma^R),\\ 
   \#\Out(R_{m_\Sigma^R}) &(m=m_\Sigma^R),\\ % \text{ and }
%      (\Sigma,R_m)\ne(F_4,D_4)),\\
%   2 &(m=m_\Sigma^R \text{ and } (\Sigma,R_m)=(F_4,D_4),\\
    1 &(m_0^R\le m<m_\Sigma^R),% \text{ and } (\Sigma,R_m)\ne(D_n,D_4))%,\\
%   3 &(m_0^R\le m<m_\Sigma^R \text{ and } (\Sigma,R_m)=(D_n,D_4))
  \end{cases}\\
  R_m^\perp\cap\Sigma = (R_m^\perp\cap R_{m_\Sigma^R})+
  \bigl\langle(\Phi_\Sigma^R)^\perp\cap\tilde\Psi
  \ (\text{or }\tilde\Psi')\bigr\rangle
  \qquad(m_0^R\le m\le m_\Sigma^R)
\end{gather}
through the natural map $G(R_m)\subset G(R_{m_\Sigma^R})\simeq G(\Phi_\Sigma^R)
\subset G(\tilde\Psi)$ {\rm(}or $G(\tilde\Psi')${\rm)} and
$R_m^\perp\cap R_{m^R_\Sigma}$ is given by i) or Lemma~\ref{lem:thm}~5).
The coset $\overline\Hom(D_4,F_4)$ consists of the two elements corresponding
to the identifications 
$D_4\simeq F_4^L$ and $D_4\simeq F_4\setminus F_4^L$.
\end{thm}
\begin{proof}
When $\Sigma$ is of type $R$ with $R=B_n$, $C_n$, $F_4$ or $G_2$,
$G(\tilde \Psi')$ is the affine Dynkin diagram $\tilde R'$ given by
Proposition~\ref{prop:clasify}. 
%Replacing the arrows of the extended Dynkin diagram of the dual 
%root system by those with the opposite directions,
%we get $G(\tilde\Psi')$ and 
This theorem follows from Lemma~\ref{lem:thm}, Remark~\ref{rem:red}~iv), 
Remark~\ref{rem:cls}~iii) and Remark~\ref{rem:orth0}.
\end{proof}
%%%             EXAMPLE    %%%%%%%%%%%%%
\begin{exmp}$\bigl(\overline{\Hom}(\Xi,\Sigma)$ and 
$\Xi^\perp\bigr)$\label{exmp:A}

{\rm i)} $\#\overline\Hom(A_2,A_n)=2$
and $A_2^\perp\cap A_n\simeq A_{n-3}\ (n\ge2)$.\\
Two elements of $\#\overline\Hom(A_2,A_n)$ are defined by 
$(\alpha_1,\alpha_2)\mapsto(\alpha_0,\alpha_1)$ and
$(\alpha_1,\alpha_2)\mapsto(\alpha_0,\alpha_n)$, respectively.
They are isomorphic to each other under $\Out(A_2)$.
Note that the rotation of the extended Dynkin diagram 
corresponds to an element of $W_{A_n}$.
\[
\begin{xy}
 \ar@{.}  *+!D{\alpha_1} *{\circ}="A";
  (5,0)  *+!D{\alpha_2} *{\ast}="B"
 \ar@{.} "B";(10,0) *+!D{\alpha_3} *{\bullet}="C"
 \ar@{-} "C";(12.5,0) \ar@{.} (12.5,0);(17.5,0)^*!U{\cdots}
 \ar@{-} (17.5,0);(20,0) %*+<5pt>!D{\alpha_{n-2}} 
  *{\bullet}="H"
 \ar@{-} "H";(25,0) *+<5pt>!D{\alpha_{n-1}} *{\bullet}="I"
 \ar@{.} "I";(30,0) *+!LD{\!\!\alpha_{n}}  *{\ast}="J"
 \ar@{-} "A";(15,-5) *+!LU{\alpha_0} *{\ccirc}="H"
 \ar@{.} "H";"J"
\end{xy}\qquad
\begin{xy}
 \ar@{.}  *+!D{\alpha_1} *{\ast}="A";
  (5,0)  *+!D{\alpha_2} *{\bullet}="B"
 \ar@{-} "B";(10,0) *+!D{\alpha_3} *{\bullet}="C"
 \ar@{-} "C";(12.5,0) \ar@{.} (12.5,0);(17.5,0)^*!U{\cdots}
 \ar@{-} (17.5,0);(20,0) %*+<5pt>!D{\alpha_{n-2}} 
  *{\bullet}="H"
 \ar@{.} "H";(25,0) *+<5pt>!D{\alpha_{n-1}} *{\ast}="I"
 \ar@{.} "I";(30,0) *+!LD{\!\!\alpha_{n}}  *{\circ}="J"
 \ar@{.} "A";(15,-5) *+!LU{\alpha_0} *{\ccirc}="H"
 \ar@{-} "H";"J"
\end{xy}
\]

{\rm ii)}
$\#\overline\Hom(A_3,D_4)=3$, 
$\#\bigl(\Out(D_4)\backslash\overline\Hom(A_3,D_4)\bigr)=1$
and $A_3^\perp\cap D_4=\emptyset$.\\
The group $\Out(D_4)\simeq\mathfrak S_3$ corresponds to the graph 
automorphisms of the extended Dynkin diagram which fix $\alpha_0$.
\[
\begin{xy}
 \ar@{-}  *+!D{\alpha_1} *{\circ} ;
  (5,0)  *+!LD{\alpha_2} *{\circ}="B"
 \ar@{.} "B";(10,0) *+!L{\alpha_3} *{\ast}="C"
 \ar@{.} "B";(5,5) *+<3pt>!LD{\alpha_4} *{\ast}
 \ar@{-} "B";(5,-5) *+!L{\alpha_0} *{\ccirc};
\end{xy}\qquad
\begin{xy}
 \ar@{.}  *+!D{\alpha_1} *{\ast} ;
  (5,0)  *+!LD{\alpha_2} *{\circ}="B"
 \ar@{.} "B";(10,0) *+!L{\alpha_3} *{\ast}="C"
 \ar@{-} "B";(5,5) *+<3pt>!LD{\alpha_4} *{\circ}
 \ar@{-} "B";(5,-5) *+!L{\alpha_0} *{\ccirc};
\end{xy}\qquad
\begin{xy}
 \ar@{.}  *+!D{\alpha_1} *{\ast} ;
  (5,0)  *+!LD{\alpha_2} *{\circ}="B"
 \ar@{-} "B";(10,0) *+!L{\alpha_3} *{\circ}="C"
 \ar@{.} "B";(5,5) *+<3pt>!LD{\alpha_4} *{\ast}
 \ar@{-} "B";(5,-5) *+!L{\alpha_0} *{\ccirc};
\end{xy}
\]

{\rm iii)}
$\#\overline\Hom(A_3,D_n) = \#\bigl(\Out(D_n)\backslash
\overline\Hom(A_3,D_n)/\!\Out(A_3)\bigr) = 2$ for $n>4$.
\[
\begin{xy}
 \ar@{-}  *+!D{\alpha_1} *{\circ} ;
  (5,0)  *+!D{\alpha_2} *{\circ}="B"
 \ar@{.} "B";(10,0) *+!D{\alpha_3} *{\ast}="C"
 \ar@{.} "C";(15,0) *+!D{\alpha_4} *{\bullet}="D"
 \ar@{-} "D";(20,0) *+!D{\alpha_5} *{\bullet}="E"
 \ar@{-} "E";(23,0) \ar@{.} (23,0);(28,0)^*!U{\!\!\!\cdots}
 \ar@{-} (28,0);(31,0) *+<5pt>!D{\alpha_{n-2}} *{\bullet}="H"
 \ar@{-} "H";(36,0) *+<5pt>!DL{\!\alpha_{n-1}} *{\bullet}
 \ar@{-} "H";(31,-5) *+!L{\alpha_n} *{\bullet}
 \ar@{-} "B";(5,-5) *+!L{\alpha_0} *{\ccirc};
\end{xy}\qquad
\begin{xy}
 \ar@{.}  *+!D{\alpha_1} *{\ast} ;
  (5,0)  *+!D{\alpha_2} *{\circ}="B"
 \ar@{-} "B";(10,0) *+!D{\alpha_3} *{\circ}="C"
 \ar@{.} "C";(15,0) *+!D{\alpha_4} *{\ast}="D"
 \ar@{.} "D";(20,0) *+!D{\alpha_5} *{\bullet}="E"
 \ar@{-} "E";(23,0) \ar@{.} (23,0);(28,0)^*!U{\!\!\!\cdots}
 \ar@{-} (28,0);(31,0) *+<5pt>!D{\alpha_{n-2}} *{\bullet}="H"
 \ar@{-} "H";(36,0) *+<5pt>!DL{\!\alpha_{n-1}} *{\bullet}
 \ar@{-} "H";(31,-5) *+!L{\alpha_n} *{\bullet}
 \ar@{-} "B";(5,-5) *+!L{\alpha_0} *{\ccirc};
\end{xy}
\]
\qquad\qquad
$A_3^\perp\cap D_n \simeq D_{n-3}$ or $D_{n-4}$ according to the imbeddings
$A_3\subset D_n$.
\smallskip

{\rm iv)}
$\#\overline\Hom(A_2,E_6)=1$ and $A_2^\perp\cap E_6\simeq 2A_2$.
Then $3A_2\subset E_6$ and
$\#\overline\Hom(3A_2,E_6)=\#\overline\Hom(2A_2,2A_2)=8$ (cf.~\S\ref{sec:4A2E8}) .
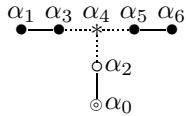
\[
\begin{xy}
 \ar@{-}  *+!D{\alpha_1} *{\bullet} ;
  (5,0)  *+!D{\alpha_3} *{\bullet}="C"
 \ar@{.} "C";(10,0) *+!D{\alpha_4} *{\ast}="D"
 \ar@{.} "D";(15,0) *+!D{\alpha_5} *{\bullet}="E"
 \ar@{-} "E";(20,0) *+!D{\alpha_6} *{\bullet},
 \ar@{.} "D";(10,-5) *+!L{\alpha_2} *{\circ}="X"
 \ar@{-} "X";(10,-10) *+!L{\alpha_0} *{\ccirc};
\end{xy}
\]

{\rm v)}
$\#\overline\Hom(A_4,E_6)=2$, 
$\#\bigl(\Out(E_6)\backslash\overline\Hom(A_4,E_6)\bigr)=1$
and $A_4^\perp\cap E_6\simeq A_1$.
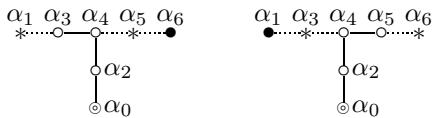
\[
\begin{xy}
 \ar@{.}  *+!D{\alpha_1} *{\ast} ;
  (5,0)  *+!D{\alpha_3} *{\circ}="C"
 \ar@{-} "C";(10,0) *+!D{\alpha_4} *{\circ}="D"
 \ar@{.} "D";(15,0) *+!D{\alpha_5} *{\ast}="E"
 \ar@{.} "E";(20,0) *+!D{\alpha_6} *{\bullet},
 \ar@{-} "D";(10,-5) *+!L{\alpha_2} *{\circ}="X"
 \ar@{-} "X";(10,-10) *+!L{\alpha_0} *{\ccirc};
\end{xy}\qquad
\begin{xy}
 \ar@{.}  *+!D{\alpha_1} *{\bullet} ;
  (5,0)  *+!D{\alpha_3} *{\ast}="C"
 \ar@{.} "C";(10,0) *+!D{\alpha_4} *{\circ}="D"
 \ar@{-} "D";(15,0) *+!D{\alpha_5} *{\circ}="E"
 \ar@{.} "E";(20,0) *+!D{\alpha_6} *{\ast},
 \ar@{-} "D";(10,-5) *+!L{\alpha_2} *{\circ}="X"
 \ar@{-} "X";(10,-10) *+!L{\alpha_0} *{\ccirc};
\end{xy}
\]

%$\bullet$ %$\circ$ $\ast$ $\star$ $\ccirc$
{\rm vi)}
$\#\overline\Hom(A_5,E_7) = 
\#\bigl(\Out(E_7)\backslash\overline\Hom(A_5,E_7)/\!\Out(A_5)\bigr)=2$.
\[
\begin{xy}
 \ar@{-}  *+!D{\alpha_0} *{\ccirc};
  (5,0)  *+!D{\alpha_1} *{\circ}="A"
 \ar@{-} "A";(10,0) *+!D{\alpha_3} *{\circ}="C"
 \ar@{-} "C";(15,0) *+!D{\alpha_4} *{\circ}="D"
 \ar@{.} "D";(20,0) *+!D{\alpha_5} *{\ast}="E"
 \ar@{.} "E";(25,0) *+!D{\alpha_6} *{\bullet}="F"
 \ar@{-} "F";(30,0) *+!D{\alpha_7} *{\bullet},
 \ar@{-} "D";(15,-5) *+!L{\alpha_2} *{\circ}
\end{xy}\qquad
\begin{xy}
 \ar@{-}  *+!D{\alpha_0} *{\ccirc};
  (5,0)  *+!D{\alpha_1} *{\circ}="A"
 \ar@{-} "A";(10,0) *+!D{\alpha_3} *{\circ}="C"
 \ar@{-} "C";(15,0) *+!D{\alpha_4} *{\circ}="D"
 \ar@{-} "D";(20,0) *+!D{\alpha_5} *{\circ}="E"
 \ar@{.} "E";(25,0) *+!D{\alpha_6} *{\ast}="F"
 \ar@{.} "F";(30,0) *+!D{\alpha_7} *{\bullet}
 \ar@{.} "D";(15,-5) *+!L{\alpha_2} *{\ast}
\end{xy}
\]
\qquad\qquad
$A_5^\perp\cap E_7 \simeq A_2$ or $A_1$ according to the imbeddings $A_5\subset E_7$.

\smallskip
{\rm vii)}
$\#\overline\Hom(4A_1,D_4)=6$ and 
$\#\bigl(\Out(D_4)\backslash\overline\Hom(4A_1,D_4)\bigr)=1$.
\[
\begin{xy}
 \ar@{.}  *+!D{\alpha_1} *{\bullet} ;
  (5,0)  *+!LD{\alpha_2} *{\ast}="B"
 \ar@{.} "B";(10,0) *+!L{\alpha_3} *{\bullet}="C"
 \ar@{.} "B";(5,5) *+<3pt>!LD{\alpha_4} *{\bullet}
 \ar@{.} "B";(5,-5) *+!L{\alpha_0} *{\ccirc}
 \ar@{}  (40,1) *{\text{$\alpha_0^\perp=\{\pm\alpha_1,\pm\alpha_3,\pm\alpha_4\}$}}
\end{xy}
\]

{\rm viii)}
$\#\overline\Hom(D_4,E_n)=1$ (cf.~Remark~\ref{rem:thmauto}~i))
\[
\begin{xy}
 \ar@{}  *+!D{\alpha_1} *{\ast} ;
  (5,0)  *+!D{\alpha_3} *{\ast}="C"
 \ar@{.} "C";(10,0) *+!D{\alpha_4} *{\circ}="D"
 \ar@{.} "D";(15,0) *+!D{\alpha_5} *{\ast}="E"
 \ar@{} "E";(20,0) *+!D{\alpha_6} *{\ast},
 \ar@{-} "D";(10,-5) *+!L{\alpha_2} *{\circ}="X"
 \ar@{-} "X";(10,-10) *+!L{\alpha_0} *{\ccirc}
 \ar@{-} "X";(5,-5) *+!R{\tilde\alpha_0} *{\dcirc}
 \ar@{}  (10,-14) *{\text{$D_4^\perp\cap E_6=\emptyset$}}
\end{xy}\qquad
\begin{xy}
 \ar@{-}  *+!D{\alpha_0} *{\ccirc};
  (5,0)  *+!D{\alpha_1} *{\circ}="A"
 \ar@{-} "A";(10,0) *+!D{\alpha_3} *{\circ}="C"
 \ar@{.} "C";(15,0) *+!D{\alpha_4} *{\circ}="D"
 \ar@{.} "D";(20,0) *+!D{\alpha_5} *{\bullet}="E"
 \ar@{}  "E";(25,0) *+!D{\alpha_6} *{\ast}="F"
 \ar@{}  "F";(30,0) *+!D{\alpha_7} *{\bullet},
 \ar@{.} "D";(15,-5) *+!L{\alpha_2} *{\bullet}
 \ar@{-} "A";(5,-5) *+!L{\tilde\alpha_0} *{\dcirc}
 \ar@{}  (15,-14) *{\text{$D_4^\perp\cap E_7\simeq3A_1$}}
\end{xy}\qquad
\begin{xy}
 \ar@{}  *+!D{\alpha_1} *{\ast};
  (5,0)  *+!D{\alpha_3} *{\bullet}="A"
 \ar@{-} "A";(10,0) *+!D{\alpha_4} *{\bullet}="C"
 \ar@{-} "C";(15,0) *+!D{\alpha_5} *{\bullet}="D"
 \ar@{.} "D";(20,0) *+!D{\alpha_6} *{\ast}="E"
 \ar@{.} "E";(25,0) *+!D{\alpha_7} *{\circ}="F"
 \ar@{-} "F";(30,0) *+!D{\alpha_8} *{\circ}="G"
 \ar@{-} "G";(35,0) *+!D{\alpha_0} *{\ccirc}
 \ar@{-} "C";(10,-5) *+!L{\alpha_2} *{\bullet}
 \ar@{-} "G";(30,-5) *+!L{\tilde\alpha_0} *{\dcirc}
 \ar@{}  (20,-14) *{\text{$D_4^\perp\cap E_8\simeq D_4$}}
\end{xy}
\]

{\rm ix)} $\#\overline\Hom(D_5,D_9)=\#\overline\Hom(D_5,B_9)=1$
(cf.~Remark~\ref{rem:thmauto}~i)).
\[
\begin{xy}
 \ar@{-}  *+!D{\alpha_1} *{\circ};
  (5,0)  *+!D{\alpha_2} *{\circ}="A"
 \ar@{-} "A";(10,0) *+!D{\alpha_3} *{\circ}="C"
 \ar@{-} "C";(15,0) *+!D{\alpha_4} *{\circ}="D"
 \ar@{.} "D";(20,0) *+!D{\alpha_5} *{\ast}="E"
 \ar@{.} "E";(25,0) *+!D{\alpha_6} *{\bullet}="F"
 \ar@{-} "F";(30,0) *+!D{\alpha_7} *{\bullet}="G"
 \ar@{-} "G";(35,0) *+!D{\alpha_8} *{\bullet}
 \ar@{-} "A";(5,-5) *+!L{\alpha_0} *{\ccirc}
 \ar@{-} "G";(30,-5) *+!L{\alpha_9} *{\bullet}
 \ar@{}  (20,-10) *{\text{$D_5^\perp\cap D_9\simeq D_4$}}
\end{xy}\qquad
\begin{xy}
 \ar@{-}  *+!D{\alpha_1} *{\circ};
  (5,0)  *+!D{\alpha_2} *{\circ}="A"
 \ar@{-} "A";(10,0) *+!D{\alpha_3} *{\circ}="C"
 \ar@{-} "C";(15,0) *+!D{\alpha_4} *{\circ}="D"
 \ar@{.} "D";(20,0) *+!D{\alpha_5} *{\ast}="E"
 \ar@{.} "E";(25,0) *+!D{\alpha_6} *{\bullet}="F"
 \ar@{-} "F";(30,0) *+!D{\alpha_7} *{\bullet}="G"
 \ar@{-} "G";(35,0) *+!D{\alpha_8} *{\bullet}="H"
 \ar@{-} "A";(5,-5) *+!L{\alpha_0} *{\ccirc}
 \ar@{=>} "H";(40,0) *+!D{\alpha_9} *{\bullet}
 \ar@{}  (20,-10) *{\text{$D_5^\perp\cap B_9\simeq B_4$}}
\end{xy}
\]

{\rm x)}
$\#\overline\Hom(A_2,F_4)=2$, $(A_2^L)^\perp\cap F_4\simeq A_2^S$, 
$(A_2^S)^\perp\cap F_4\simeq A_2^L$.
\begin{align*}
\begin{xy}
 \ar@{-} *+!D{\alpha_0} *{\ccirc} ;
   (5,0)  *+!D{\alpha_1} *{\circ}="C"
 \ar@{.} "C";(10,0) *+!D{\alpha_2} *{\ast} ="D"
 \ar@{:>} "D";(15,0) *+!D{\alpha_3} *{\bullet} ="E"
 \ar@{-} "E";(20,0) *+!D{\alpha_4} *{\bullet}="F"
\end{xy}
\qquad
\begin{xy}
 \ar@{-} *+!D{\alpha_1} *{\bullet} ;
   (5,0)  *+!D{\alpha_2} *{\bullet}="C"
 \ar@{:>} "C";(10,0) *+!D{\alpha_3} *{\ast} ="D"
 \ar@{.} "D";(15,0) *+!D{\alpha_4} *{\circ} ="E"
 \ar@{-} "E";(20,0) *+!D{\alpha_0'} *{\ccirc}="F"
\end{xy}
\end{align*}

{\rm xi)}
$\#\overline\Hom(C_3,F_4)=1$, $G(C_4)\subset 
G(\tilde F_4')$, $G(C_3)\subset G(\tilde C_4)$ 
and $C_3^\perp\cap F_4=A_1^L$.
\begin{align*}
\begin{xy}
 \ar@{} *+!D{\alpha_1} *{\ast} ;
   (5,0)  *+!D{\alpha_2} *{\bullet}="C"
 \ar@{:>} "C";(10,0) *+!D{\alpha_3} *{\ast} ="D"
 \ar@{.} "D";(15,0) *+!D{\alpha_4} *{\circ} ="E"
 \ar@{-} "E";(20,0) *+!D{\alpha_0'} *{\ccirc}="F"
 \ar@{<=} "F";(25,0) *+!D{\tilde\alpha_0} *{\dcirc}="F"
\end{xy}\qquad
C_4=\{\alpha_2,\alpha_3,\alpha_4,\alpha_0'\},\ 
C_3=\{\alpha_4,\alpha_0',\tilde\alpha_0\}.
%\ C_3^\perp\cap F_4=\langle\alpha_2\rangle
\end{align*}

{\rm xii)} $\#\overline\Hom(A_4+A_2,E_8)=2$ and
$\#\bigl(\overline\Hom(A_4+A_2,E_8)/\!\Out(A_4+A_2)\bigr)=1$.\\
The two reductions $\overline\Hom(A_2,A_4^\perp)$ and 
$\overline\Hom(A_4,A_2^\perp)$ given below lead to the same result
(cf.~Remark~\ref{rem:red}).
In particular $(A_4+A_2)^\perp\cap E_8\simeq A_1$.
Note that $(A_4,A_4)$ and $(A_2,E_6)$ are special dual pairs
in $E_8$ (cf.~Definition~\ref{def:dual}).
\begin{align*}
\begin{xy}
 \ar@{-}  *+!D{\alpha_1} *{\bullet};
  (5,0)  *+!D{\alpha_3} *{\bullet}="A"
 \ar@{-} "A";(10,0) *+!D{\alpha_4} *{\bullet}="C"
 \ar@{.} "C";(15,0) *+!D{\alpha_5} *{\ast}="D"
 \ar@{.} "D";(20,0) *+!D{\alpha_6} *{\circ}="E"
 \ar@{-} "E";(25,0) *+!D{\alpha_7} *{\circ}="F"
 \ar@{-} "F";(30,0) *+!D{\alpha_8} *{\circ}="G"
 \ar@{-} "G";(35,0) *+!D{\alpha_0} *{\ccirc}
 \ar@{-} "C";(10,-5) *+!L{\alpha_2} *{\bullet}
% \ar@{-} "G";(30,-5) *+!L{\alpha'_0} *{\bullet}
\end{xy}
&\ \to\ &
\begin{xy}
 \ar@{.}  *+!D{\alpha_1} *{\circ}="O";
  (5,0)  *+!D{\alpha_3} *{\ast}="A"
 \ar@{.} "A";(10,0) *+!D{\alpha_4} *{\bullet}="C"
 \ar@{.} "C";(15,0) *+!D{\alpha_2} *{\ast}="D"
 \ar@{-} "O";(7.5,-5) *{\dcirc}="E"
 \ar@{.} "E";"D"
\end{xy}
&\text{ \ or \ }
\begin{xy}
 \ar@{.}  *+!D{\alpha_1} *{\ast}="O";
  (5,0)  *+!D{\alpha_3} *{\bullet}="A"
 \ar@{.} "A";(10,0) *+!D{\alpha_4} *{\ast}="C"
 \ar@{.} "C";(15,0) *+!D{\alpha_2} *{\circ}="D"
 \ar@{.} "O";(7.5,-5) *{\dcirc}="E"
 \ar@{-} "E";"D"
\end{xy}\\
\begin{xy}
 \ar@{-}  *+!D{\alpha_1} *{\bullet};
  (5,0)  *+!D{\alpha_3} *{\bullet}="A"
 \ar@{-} "A";(10,0) *+!D{\alpha_4} *{\bullet}="C"
 \ar@{-} "C";(15,0) *+!D{\alpha_5} *{\bullet}="D"
 \ar@{-} "D";(20,0) *+!D{\alpha_6} *{\bullet}="E"
 \ar@{.} "E";(25,0) *+!D{\alpha_7} *{\ast}="F"
 \ar@{.} "F";(30,0) *+!D{\alpha_8} *{\circ}="G"
 \ar@{-} "G";(35,0) *+!D{\alpha_0} *{\ccirc}
 \ar@{-} "C";(10,-5) *+!L{\alpha_2} *{\bullet}
\end{xy}
&\ \to\ &
\begin{xy}
 \ar@{.}  *+!D{\alpha_1} *{\ast} ;
  (5,0)  *+!D{\alpha_3} *{\circ}="C"
 \ar@{-} "C";(10,0) *+!D{\alpha_4} *{\circ}="D"
 \ar@{.} "D";(15,0) *+!D{\alpha_5} *{\ast}="E"
 \ar@{.} "E";(20,0) *+!D{\alpha_6} *{\bullet},
 \ar@{-} "D";(10,-5) *+!L{\alpha_2} *{\circ}="X"
 \ar@{-} "X";(10,-10)               *{\dcirc};
\end{xy}&\text{ \ or \ }
\begin{xy}
 \ar@{.}  *+!D{\alpha_1} *{\bullet} ;
  (5,0)  *+!D{\alpha_3} *{\ast}="C"
 \ar@{.} "C";(10,0) *+!D{\alpha_4} *{\circ}="D"
 \ar@{-} "D";(15,0) *+!D{\alpha_5} *{\circ}="E"
 \ar@{.} "E";(20,0) *+!D{\alpha_6} *{\ast},
 \ar@{-} "D";(10,-5) *+!L{\alpha_2} *{\circ}="X"
 \ar@{-} "X";(10,-10)               *{\dcirc};
\end{xy}
\end{align*}
\end{exmp}
%%%
\begin{cor} \label{cor:A}
{\rm i)}
Suppose $\Sigma$ is not of type $A$.
Let $G(\{\alpha_0,\alpha_{j_1},\ldots,\alpha_{j_{m-1}}\})$ be a maximal 
subdiagram of $G(\tilde\Psi)$ isomorphic to $G(A_{m})$
such that $\alpha_0$ and $\alpha_{j_{m-1}}$ are the end vertices of 
the subdiagram and $\alpha_{j_\nu}$ are not the branching vertex of 
$G(\tilde\Psi)$ for $\nu=1,\ldots,m-2$.
Then
\begin{align*}
 \#\overline\Hom'(A_k,\Sigma)&
   =1\qquad(k=1,\ldots,m),\\
 \#\overline\Hom'(A_{m+1},\Sigma)
 &\begin{cases}
  = 0&
 \bigl(\alpha_{j_{m-1}}\text{ is not a branching vertex of }G(\tilde\Psi)
 \bigr)\\
  > 1&
 \bigl(\alpha_{j_{m-1}}\text{ is a branching vertex of }G(\tilde\Psi)
 \bigr)
 \end{cases}
\end{align*}
with
\[
 m = 
  \begin{cases}
    2    &(\Sigma=B_n,\ n\ge3),\qquad
    1\quad(\Sigma=B_2,\ C_n),\\
    2    &(\Sigma=D_n,\ n\ge4),\\
    3    &(\Sigma=E_6),\qquad
    4\quad(\Sigma=E_7),\qquad
    6\quad(\Sigma=E_8),\\
    3    &(\Sigma=F_4),\,\qquad
    2\quad(\Sigma=G_2).
  \end{cases}
\]
Here $\alpha_{j_{m-1}}$ is the branching vertex if $\Sigma=B_n$ {\rm($n\ge3$)},
$D_n$ {\rm($n\ge4$)}, $E_6$, $E_7$ or $E_8$.

{\rm ii)}
We consider the following procedure for a Dynkin diagram $X$:
%%%
\begin{quote}
If $X$ is connected, we replace it by the subdiagram $X'$ of the 
extended Dynkin diagram $\tilde X$ of $X$ where the vertices of $X'$ 
correspond to the roots orthogonal to the maximal root of $\tilde X$.
If an irreducible component of $X'$ has no root with the length of 
the maximal root, we remove the component.

\noindent
If $X$ is not connected, we choose one of the connected
component of $X$ and change the component by the above procedure.
\end{quote}

Then $\overline\Hom'(rA_1,\Sigma)$ corresponds to the $r$ steps
of the above procedures starting from $G(\Psi)$.
The existence of these steps implies $\Hom'(rA_1,\Sigma)\ne\emptyset$
and in this case $\overline\Hom'(rA_1,\Sigma)=1$ if and only if
any non-connected Dynkin diagram does not appear except for the final
step.
In particular, we have the following:

Let $r(\Sigma)$ be the maximal integer $r$ satisfying
$\Hom'(rA_1,\Sigma)\ne\emptyset$.
Then
\begin{align}
  r(\Sigma) &= 1 + \sum_j r(\Sigma'_j).
\end{align}
Here $\{\Sigma'_j\}$ 
is the set of irreducible components
of $\alpha_0^\perp$ such that $\Sigma'_j\cap\Sigma^L\ne\emptyset$ and
\begin{align*}
 r(A_n) &= 1+r(A_{n-2}) = [\tfrac{n+1}2]&(n\ge2),\quad r(A_1)&=1,\ r(A_0)=0,\\
 r(B_n) &= 2+r(B_{n-2}) = 2[\tfrac n2]&(n\ge4),\quad r(B_3)&= r(B_2)= 2,\\
 r(C_n) &= 1+r(C_{n-1}) = n&(n\ge3),\quad r(C_2)&= 2,\\
 r(D_n) &= 2+r(D_{n-2}) = 2[\tfrac n2]&(n\ge4),\quad r(D_3)&= r(D_2)= 2,\\
 r(E_6) &= 1+r(A_5) = 4,&
 r(E_7) &= 1+r(D_6) = 7,\\
 r(E_8) &= 1+r(E_7) = 8,\\
 r(F_4) &= 1+r(C_3) = 4,&
 r(G_2) &= 1.
\end{align*}
\end{cor}
\begin{rem} \label{rem:A1}
{\rm i)}
If $\Sigma$ is of type $A$, $D$ or $E$, 
then $\overline\Hom(rA_1,\Sigma)$ is figured as follows
according to the procedures in Corollary~\ref{cor:A}~ii) and
the notation in \S\ref{sec:Dynkin}.
\[
\begin{xy}
 \ar@{->} *+{A_n};(15,0)
  *+{A_{n-2}}="A"
 \ar@{->}
  "A";(32,0) *+{A_{n-4}}="B"
 \ar@{->}  "B";(49,0) *+{A_{n-6}}="C"
 \ar@{->} "C";(64,0) *+{\cdots}
\end{xy}
\qquad
\begin{xy} \ar@{->} *+{D_n};(20,0)
  *+{D_{n-2}+A_1}="A"
 \ar@{->}
  "A";(39,0) *+{\cdots}
\end{xy}
\]
\[
\begin{xy}
 \ar@{->}^*{\alpha_0} *+{E_6};(15,0)
  *+{A_5}="A"
 \ar@{->}^*{\frac{\epsilon_1+\cdots+\epsilon_4-\epsilon_5+\epsilon_6+\epsilon_7-\epsilon_8}2}
  "A";(40,0) *+{A_3}="B"
 \ar@{->}^*{\epsilon_4-\epsilon_1}  "B";(60,0) *+{A_1}="C"
 \ar@{->}^*{\epsilon_3-\epsilon_2}  "C";(80,0)  *+{\emptyset}
\end{xy}
\qquad\qquad
\begin{xy}
  \ar@{->}^*{-\epsilon_7-\epsilon_8} *+{E_8};(18,0)
  *+{E_7}="A"
\end{xy}
\]
\[
\begin{xy}
 \ar@{->}^*{\epsilon_7-\epsilon_8} *+{E_7};(15,0)
  *+{D_6}="A"
 \ar@{->}^*{\!\!\!\!\!\!\!\!\!-\epsilon_5-\epsilon_6} "A";(40,0) *+{D_4+A_1}="B"
 \ar@{->}^<(0.5)*{\epsilon_6-\epsilon_5}  "B";(60,-8) *+{D_4}="C"
 \ar@{->}^<(0.2)*{\!\!\!\!\!\!\!-\epsilon_3-\epsilon_4}  "B";(60,8)  *+{4A_1}="D"
 \ar@{->} "C";(75,-8) *+{3A_1}="E"
 \ar@{->} "D";(75,8)  *+{3A_1}="F"
 \ar@3{->} "D";"E"
 \ar@3{->} "E";(90,0) *+{2A_1}="G"
 \ar@3{->} "F";"G"
 \ar@{=>} "G";(105,0) *+{A_1}="H"
 \ar@{->} "H";(120,0)  *+{\emptyset}
\end{xy}
\]
There appears $3A_1$ twice in the above. They are distinguished 
by the structure of $(3A_1)^\perp\cap E_7$ but isomorphic
in $E_8$ (cf.~\S\ref{sec:TE7}, \S\ref{sec:TE8} and \S\ref{sec:8A1E8}). 

For example, it follows from the procedures shown above that
\begin{equation}
 \begin{split}
 \#\overline\Hom(5A_1,E_7)
 &= \#\overline\Hom(4A_1,D_6)
 = \#\overline\Hom(3A_1,D_4+A_1)\\
 &= \#\overline\Hom(2A_1,D_4)+\#\overline\Hom(2A_1,4A_1)\\
 &= \#\overline\Hom(A_1,3A_1)+4\#\overline\Hom(A_1,3A_1)
 = 3+4\cdot 3 = 15.
 \end{split}
\end{equation}

{\rm ii)} For an irreducible root system $\Sigma$, we can easily
calculate $\#\overline\Hom(\Xi,\Sigma)$ and $\Xi^\perp\cap\Sigma$
for any root system $\Xi$ in virtue of Theorem~\ref{thm} together with 
Remark~\ref{rem:red} (cf.~Example \ref{rem:thmauto} x)).
The complete list for non-trivial $\overline\Hom(\Xi,\Sigma)$ is
given in \S\ref{sec:table}.
More refined structures related to the actions of
$\Out(\Sigma)$ and $\Out(\Xi)$ etc.\ are also given in \S\ref{sec:table}, 
which will be studied in later sections.
\end{rem}
%%%%%%%%%%%%%%%%%%%%%%%%%  Lemmas  %%%%%%%%%%%%%%%%%%%%%%%%%%
\section{Lemmas}\label{sec:lem}
In this section we always assume that $\Psi$ is a fundamental system
of an \textsl{irreducible} root system $\Sigma$ and $\tilde\Psi$ is 
the corresponding extended fundamental system. 

First note that for $\alpha\in\tilde\Psi\cap\Sigma^L$
we have
\begin{equation}\label{eq:inprod}
 2\frac{(\alpha|\beta)}{(\alpha|\alpha)}
 \in \begin{cases}
    \{0,-1\} &(\forall\beta\in\langle\Psi\setminus\{\alpha\}\rangle
       \text{ and }\beta>0),\\
    \{0,1\}  &(\forall\beta\in\langle\Psi\setminus\{\alpha\}\rangle
       \text{ and }\beta<0).
  \end{cases}
\end{equation}
Here we put $\Psi\setminus\{\alpha\}=\Psi$ if $\alpha\notin\Psi$.
%and the above is clear from the root system generated by $\alpha$ and $\beta$.

\begin{lem}\label{lem:orth}
If a subset $\Theta$ of $\tilde\Psi$ contains $\alpha_0$ and 
the diagram $G(\Theta)$ is connected,
\begin{equation}
 \Theta^\perp=
  \bigl\langle\alpha\in\tilde\Psi\,;\,\alpha\perp\Theta\bigr\rangle.
\end{equation}
\end{lem}
\begin{proof}
Note that $(\alpha_i|\alpha_j)\le 0$ for $0\le i<j\le n$.\\
We will prove the lemma by the induction on $\#\Theta$.

We may put $\Theta=\{\alpha_0,\alpha_1,\ldots,\alpha_m\}$ and we may assume
$\Theta':=\{\alpha_0,\alpha_1,\ldots,\alpha_{m-1}\}$ is empty or forms
a connected subdiagram.

Let $\alpha=\sum_{j=1}^n m_j(\alpha)\alpha_j\in\Theta^\perp$
with $m_j(\alpha)\ge0$.
Then the induction hypothesis for $\Theta'$ implies $m_j(\alpha)=0$ for
$j\le m$ and
\[
 0=(\alpha_m|\alpha)
  =\textstyle\sum_{j={m+1}}^n m_j(\alpha)(\alpha_m|\alpha_j). 
\]
Hence $(\alpha_m|\alpha_j)\ne 0$ means $m_j(\alpha)=0$.
\end{proof}
%%%
\begin{rem}\label{rem:orth0}
In Lemma~\ref{lem:orth} we may replace $\alpha_0$ by any element 
$\alpha_0'$ satisfying $(\alpha_0'|\alpha)\le0$ for all $\alpha\in\Psi$.
Then the Dynkin diagram $G(\tilde\Psi')$ of 
$\tilde\Psi'=\Psi\cup\{\alpha_0'\}$ is an affine Dynkin diagram
in Proposition~\ref{prop:clasify}.
\end{rem}
%%%
\begin{lem}\label{lem:fibre}
Fix $\Theta\subset \Psi$ and 
$\mathbf m\in\mathbb Z^{\#\Theta}\setminus\{0\}$.
Define the map
\[
  \begin{matrix}
   p_\Theta\,: &\Sigma & \to & \mathbb Z^{\#\Theta}\\
          &\vin && \vin\\
   & \beta=\sum_{\alpha_i\in\Psi}m_i(\beta)\alpha_i & \mapsto &
   \bigl(m_i(\beta)\bigr)_{\alpha_i\in\Theta}
  \end{matrix}.
\]
Then $p_\Theta^{-1}(\mathbf m)\cap\Sigma^L$ is empty or 
a single $W_{\Psi\setminus\Theta}$-orbit.
Moreover $p_\Theta^{-1}(\mathbf m)\setminus\Sigma^L$ is also empty or 
a single $W_{\Psi\setminus\Theta}$-orbit.
\end{lem}
%%%
\begin{proof}
Fix $0\ne \mathbf m=(m_i)_{\alpha\in\Theta}$ 
in the image of $p_\Theta$.

Let $\mathfrak g$ be the complex simple Lie algebra with the root
system $\Sigma$ and let $X_\alpha\in\mathfrak g$ be a root vector
for $\alpha\in\Sigma$.
We denote by $\mathfrak g_{\Psi\setminus\Theta}$ the semisimple Lie algebra
generated by $\{X_\alpha\,;\,\alpha\in\Psi\setminus\Theta\}$.
Then the space 
\[
 V_{\mathbf m}:=\textstyle\sum_{\alpha\in p_\Theta^{-1}(\mathbf m)}\mathbb CX_\alpha
 \subset\mathfrak g
\]
is a $\mathfrak g_{\Psi\setminus\Theta}$-stable subset
under the adjoint representation of $\mathfrak g$,
which is an irreducible representation of 
$\mathfrak g_{\Psi\setminus\Theta}$ as is shown in \cite[Proposition~2.39 ii)]{OO}.

Let $\pi_\Theta$ be the orthogonal projection of
$\mathbb R^n$ onto $\sum_{\alpha\in\Psi\setminus\Theta}\mathbb R\alpha$
with respect to $(\ |\ )$.
Put $v_{\mathbf m} = \sum_{\alpha_i\in\Theta}m_i\alpha_i - 
  \pi_\Theta(\sum_{\alpha_i\in\Theta}m_i\alpha_i)$.
Then
\[
 \pi_\Theta(\alpha) = \alpha - v_{\mathbf m},\ 
 (\pi_\Theta(\alpha)|\pi_\Theta(\alpha))
  =(\alpha|\alpha) - (v_{\mathbf m}|v_{\mathbf m})
 \quad(\forall\alpha\in p_\Theta^{-1}(\mathbf m)).
\]
The set of the weights of the irreducible representation 
$(\mathfrak g_{\Psi\setminus\Theta}, V_{\mathbf m})$
is $\pi_\Theta(p_\Theta^{-1}(\mathbf m))$ and the set of the weights
with the longest length %(i.e.\ extreme weights) 
is $\pi_\Theta
  (p_\Theta^{-1}(\mathbf m)\cap\Sigma^L)$.
Hence $p_\Theta^{-1}(\mathbf m)\cap\Sigma^L$ is a single
$W_{\Psi\setminus\Theta}$-orbit.

When $p_\Theta^{-1}(\mathbf m)\not\subset\Sigma^L$, 
we have the last statement in the lemma by 
combining the above argument with \cite[Proposition~2.37 ii)]{OO}.
\end{proof}
\begin{prop}\label{prop:subsystem}
For a proper subset $\Theta$ of the extended fundamental system 
$\tilde\Psi$ of $\Sigma$
we have
\begin{equation}
 \langle\Theta\rangle
 = \begin{cases}
    p_{\Psi\setminus\Theta}^{-1}\bigl(0\bigr) &(\alpha_0\notin\Theta),\\
    p_{\Psi\setminus\Theta}^{-1}\bigl(\{0,\pm p_{\Psi\setminus\Theta}(\alpha_0)\}\bigr)
    &(\alpha_0\in\Theta)
   \end{cases}
\end{equation}
under the notation in Lemma~\ref{lem:fibre}.
\end{prop}
\begin{proof}
Note that
$\langle\Theta\rangle\supset p_{\Psi\setminus\Theta}^{-1}\bigl(0\bigr)$.
We assume $\alpha_0\in\Psi$ because the claim is clear when 
$\alpha_0\notin\Theta$.
Then \eqref{eq:inprod} implies
$\langle \Theta\rangle\subset p_{\Psi\setminus\Theta}^{-1}
\bigl(\{0,\pm p_{\Psi\setminus\Theta}(\alpha_0)\}\bigr)$.
Let $\Theta_0$ be the irreducible component of $\Theta$ containing 
$\alpha_0$.
Since $\langle\Theta\rangle$ is $W_{\Theta\setminus\{\alpha_0\}}$-invariant,
Lemma~\ref{lem:fibre} implies that
\[
 \langle\Theta\rangle\setminus
 p_{\Psi\setminus\Theta}^{-1}\bigl(0\bigr)=
  p_{\Psi\setminus\Theta}^{-1}
  \bigl(\{\pm p_{\Psi\setminus\Theta}(\alpha_0)\}\bigr)%\\
  \text{ or }
  p_{\Psi\setminus\Theta}^{-1}
  \bigl(\{\pm p_{\Psi\setminus\Theta}(\alpha_0)\}\bigr)
  \cap\Sigma^L. %\Bigr).
\]
Hence the proposition is clear if
$\Sigma$ is simply laced %of type $A_n$, $D_n$ or $E_n$
or if $\Theta_0$ is of type $B_n$ or $C_n$.
It is also easy to check $p_{\Psi\setminus\Theta}^{-1}
  \bigl(p_{\Psi\setminus\Theta}(\alpha_0)\bigr)
  \subset \Sigma^L$
in any other case when
$(\Psi,\Theta_0)=(B_n,A_{m-1})$, $(B_n,D_m)$, $(C_n,A_1)$ or $(F_4,A_k)$
with $m\le n$ and $k\le 3$.
\end{proof}
\begin{lem}[roots orthogonal to the end root]\label{lem:orthe}
Suppose $\alpha_1$ is an end root of $\Psi$ with
$\alpha_1\in\Sigma^L$.
Then the set
\begin{equation}
 Q=\bigl\{\alpha=\alpha_1+m_2(\alpha)\alpha_2+m_3(\alpha)\alpha_3+
   \cdots+m_n(\alpha)\alpha_n\in\Sigma^L\,;\,
  (\alpha|\alpha_1)=0\bigr\}
\end{equation}
is empty if $\Psi$ is of type $A$ and it is a single 
$W_{\Psi\cap\alpha_1^\perp}$-orbit if otherwise.
\end{lem}
%%%
\begin{proof}
We may assume $\#\Psi>1$.  Then there is a unique
$\beta\in\Psi$ with $(\alpha_1|\beta)<0$.
We may assume $\beta=\alpha_2$ and we have
\[
 Q=\bigl\{\alpha=\alpha_1+2\alpha_2
 +m_3(\alpha)\alpha_3+\cdots+m_n(\alpha)\alpha_n
 \in\Sigma^L\bigr\}.
\]
Then $Q=\emptyset$ if and only if $\Psi$ is of type $A$.
If $\Psi$ is not of type $A$, Lemma~\ref{lem:fibre} assures that
$Q$ is a single $W_{\Psi\setminus\{\alpha_1,\alpha_2\}}$-orbit.
Note that 
$\Psi\cap \alpha_1^\perp = \Psi\setminus\{\alpha_1,\alpha_2\}$.
\end{proof}
%%%
\begin{lem}[special imbeddings of $A_2$ and $A_3$]\label{lem:A2}
Let $\Psi'\subset\Psi$.
If $\Psi'\ne\Psi$, we assume that we can choose
$\alpha'\in\Psi\cap\Sigma^L$ with $\alpha'\notin\Psi'$.
If $\Psi'=\Psi$, we put $\alpha'=\alpha_0$.
Define 
\begin{align*}
 Q_1&:=\bigl\{\beta\in\langle\Psi'\rangle\cap\Sigma^L\,;\,
   (\beta|\alpha')<0\bigr\},\\
 Q_2&:=\bigl\{(\beta_1,\beta_2)\in\langle(\Psi'\rangle\cap\Sigma^L)
          \times(\langle\Psi'\rangle\cap\Sigma^L)\,;\,\\
  &
  \qquad (\beta_1,\alpha')=(\beta_2|\alpha')<0
  \text{ and } (\beta_1|\beta_2)=0\bigr\},\\
 \Theta&:=\bigl\{\alpha\in\Psi'\,;\,(\alpha|\alpha')<0\bigr\},\\
 \Theta^L&:= \Theta\cap\Sigma^L.
\end{align*}
Then $\Theta^L$ is the set of complete representatives of 
$Q_1/W_{\Psi'\setminus\Theta}$. 
Moreover if $\Psi'\ne\Psi$, 
\begin{align*}
    Q_2/\#W_{\Psi'\setminus\Theta}
  &= \#\Theta^L\bigl(\#\Theta^L - 1\bigr) \\
  &\quad + \bigl\{\alpha\in\Theta^L\,;\,G(\Psi'_\alpha)\text{ is not of type }A
   \text{ or not an end root of }G(\Psi'_\alpha)\bigr\}.
\end{align*}
Here $\Psi'_\alpha$ is the irreducible component of 
$\Psi'$ containing $\alpha\in\Theta$.
\end{lem}
%%%
\begin{proof}
Let $\beta\in Q_1$.
It follows from \eqref{eq:inprod} that there exists 
$\alpha_m\in\Psi'$ satisfying
\begin{gather}
  \beta = \alpha_m + \sum_{\alpha_j\in\Psi'\setminus\Theta}
 m_j(\beta)\alpha_j,\label{eq:A2}\\
  (\alpha_m|\beta)<0.
\end{gather}
If $\alpha_m\notin\Sigma^L$, $\Psi'_{\alpha_m}$ 
is of type $A$ or $C$ and therefore $\beta$ of the 
form \eqref{eq:A2} does not belong to $\Sigma^L$.
Hence $\alpha_m\in\Theta^L$ and $\alpha_m\in W_{\Psi'\setminus\Theta}\beta$ by
lemma~\ref{lem:fibre}.

Let $\alpha_m$, $\alpha_{m'}\in\Theta^L$ with $m\ne m'$.
We have $\alpha_{m'}\not\in W_{\Psi'\setminus\Theta}\alpha_m$
and therefore $\Theta^L$ is the set of complete representatives of 
$Q_1/W_{\Psi'\setminus\Theta}$. 

Let $(\beta_1,\beta)\in Q_2$.
We may assume $\beta_1=\alpha_k\in\Theta^L$ by the argument above
and $\beta$ is of the form \eqref{eq:A2} with $\alpha_m\in\Theta^L$.

If $k\ne m$, we may similarly assume $\beta=\alpha_m$ and
$(\alpha_k,\alpha_m)\in Q_2$.
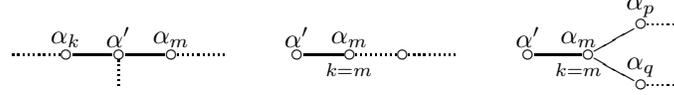
\[\begin{xy}
 \ar@{.}  *{};
  (7,0)  *+!D{\alpha_k} *{\circ}="B"
 \ar@{-} "B";(14,0) *+!D{\alpha'} *{\circ}="C"
 \ar@{.} "C";(14,-4) *{}
 \ar@{-} "C";(21,0) *+!D{\alpha_m} *{\circ} = "D"
 \ar@{.} "D";(28,0) *{}
\end{xy}
\qquad
\begin{xy}
 \ar@{-}  *+!D{\alpha'} *{\circ};
  (7,0)   *+!D{\alpha_m} *+!U{{}^{k=m}} *{\circ}="B"
 \ar@{.} "B";(14,0) *{\circ}="C"
 \ar@{.} "C";(21,0) *{}
\end{xy}
\qquad
\begin{xy}
 \ar@{-}  *+!D{\alpha'} *{\circ};
  (8,0)   *+!D{\!\!\!\!\alpha_m} *+!U{\!\!\!\!{}^{k=m}} *{\circ}="B"
 \ar@{-} "B";(15,4) *+!D{\alpha_p} *{\circ}="D"
 \ar@{.} "D";(20,4) *{}
 \ar@{-} "B";(15,-4)  *+!D{\alpha_q} *{\circ}="E"
 \ar@{.} "E";(20,-4)  *{}
\end{xy}
\]
Suppose $k=m$.
If $\alpha_m$ is the end root of $\Psi_{\alpha_m}$, it follows from
Lemma~\ref{lem:orthe} that 
$\Psi'_{\alpha_m}$ is not of type $A$ and 
$(\alpha_k,\beta)$ corresponds to a unique element of 
$Q_2/W_{\Psi'\setminus\Theta}$.

If $\alpha_m$ is not the end root of $\Psi'_{\alpha_m}$, 
$\Psi'_{\alpha_m}$ is of type $A$ and it is easy to see that 
$(\alpha_k,\beta)$ also corresponds to a unique class in 
$Q_2/W_{\psi'\setminus\Theta}$.
In fact, we may put
$\{\alpha\in\Psi_{\alpha_m}\,;\,(\alpha|\alpha_m)<0\}=
\{\alpha_p,\alpha_q\}$ and 
\[
  \beta = \alpha_m+\alpha_p+\alpha_q
   + \sum_{\alpha_j\in\Psi'\setminus\{\alpha_m,\alpha_p,\alpha_q\}}
%{\substack{\alpha_j\in \Psi'\\ j\ne m,\,p,\,q}}
    m_j(\beta)\beta\in\Sigma.
\]
Note that the roots $\beta$ with this expression are in
a single $W_{\Psi'\setminus\{\alpha_m,\alpha_p,\alpha_q\}}$-orbit.
Thus we have the lemma.
\end{proof}

%%%%%%%%%%%%%%%%%%%%%%%%%%%%  Proof **************************
\section{Proof of the main Lemma}\label{sec:prf}
Retain the notation in Lemma~\ref{lem:thm} to prove it.

1) Let $\Xi$ be of type $A_{m+1}$ with the fundamental system
$\Phi=\{\beta_0,\ldots,\beta_m\}$ and the Dynkin diagram
\begin{xy}
 \ar@{-} *+!D{\beta_0} *{\circ}="Z";
  (7,0) *+!D{\beta_1} *{\circ}="A"
 \ar@{-} "A";(10,0) \ar@{.} (10,0);(15,0)^*!U{\cdots}
 \ar@{-} (15,0);(18,0) *{\circ}="C"
 \ar@{-} "C";(25,0) *+<5pt>!D{\beta_{m-1}} *{\circ}="B"
 \ar@{-} "B";(32,0) *+!D{\beta_m} *{\circ}
\end{xy}.

First note that $\bar\iota$ naturally corresponds to an element
of $\Hom'(\Xi,\Sigma)$ and then \eqref{eq:orthA} follows from lemma~\ref{lem:orth}.
We will prove the lemma by the induction on $m$.

Let $\iota\in\Hom'(\Xi,\Sigma)$.
Since 
$\{\alpha\in\Sigma\,;\,|\alpha|=|\alpha_{max}|\}=W_\Sigma\alpha_{max}$,
the lemma is clear when $m=0$.
Suppose $m\ge1$.
By the induction hypothesis we may assume that there exists a
unique sequence $(\alpha_0,\ldots,\alpha_{m-1})$ of 
element of $\tilde \Psi$ and an element $w\in W_\Sigma$ such that
$w\circ\iota(\beta_j)=\alpha_j$ for $j=0,\ldots,m-1$.
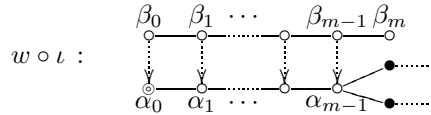
\[
 w\circ\iota\,:\quad
\raisebox{10pt}{
 \begin{xy}
 \ar@{-} *+!D{\beta_0} *{\circ}="Z";
  (7,0) *+!D{\beta_1} *{\circ}="A"
 \ar@{-} "A";(10,0) \ar@{.} (10,0);(15,0)^*!U{\cdots}
 \ar@{-} (15,0);(18,0) *{\circ}="C"
 \ar@{-} "C";(25,0) *+<5pt>!D{\beta_{m-1}} *{\circ}="B"
 \ar@{-} "B";(32,0) *+!D{\beta_m} *{\circ}
 \ar@{.>} "Z";(0,-7) *+!U{\alpha_0} *{\ccirc}="P"
 \ar@{-} "P";(7,-7) *+!U{\alpha_1} *{\circ}="Q"
 \ar@{-} "Q";(10,-7) \ar@{.} (10,-7);(15,-7)_*{\cdots}
 \ar@{-} (15,-7);(18,-7) *{\circ}="U"
 \ar@{-} "U";(25,-7) *+<5pt>!U{\alpha_{m-1}} *{\circ}="R"
 \ar@{-} "R";(32,-4) *{\bullet}="S"
 \ar@{-} "R";(32,-9) *{\bullet}="T"
 \ar@{.} "S";(37,-4)
 \ar@{.} "T";(37,-9)
 \ar@{.>} "A";"Q"
 \ar@{.>} "B";"R"
 \ar@{.>} "C";"U"
 \end{xy}
}
\]
Put $\alpha'_m=w\circ\iota(\beta_m)$ and
\begin{align*}
   \Psi'&= \{\alpha\in\tilde\Psi\,;\,
          (\alpha|\alpha_j)=0
          \quad(j=0,\ldots,m-2)\},\\
   \Theta &= \{\alpha\in\Psi'\,;\,(\alpha|\alpha_{m-1})<0\}.
\end{align*}
Since $(\alpha'_m|\alpha_j)=0$ for $j=0,\ldots,m-2$, 
$\alpha'_m\in\langle \Psi'\rangle$. 
Applying Lemma~\ref{lem:A2} to $\alpha':=\alpha_{m-1}$, we have
$\alpha_m\in\Theta\cap\Sigma^L$ and 
$w'\in W_{\Psi'\setminus\Theta}$ such that 
$w'(\alpha'_m)=\alpha_m$.
Hence $w'w\circ\iota$ corresponds to a required imbedding of $G(\Phi)$
into $G(\Psi)$.

The uniqueness of $\alpha_m\in\Theta\cap\Sigma^L$ is proved as follows.
Suppose there exists $w\in W_\Sigma$ such that
\[
  w\alpha_j = \alpha_j\quad\text{for }j=0,\ldots,m-1
  \quad\text{and}\quad
  w\alpha_m\in\Theta\cap\Sigma^L.
\]
Then $w\in W_{\Psi'\setminus\Theta}$ and Lemma~\ref{lem:A2} assures
$w\alpha_m=\alpha_m$.

Thus we have proved the first claim and then Lemma~\ref{lem:orth} assures
\eqref{eq:orthA}.  The last claim is easily obtained by applying the claims
we have proved to the extended Dynkin diagrams in \S\ref{sec:Dynkin}.

%%%%%%%%%%%%
\medskip
2) Let $\Xi$ is of type $D_m$ with $m\ge4$.
We may assume that $\Sigma$ is of type $B_n$, $D_n$, $E_n$ or $F_n$.
Let $\iota\in\Hom'(\Xi,\Sigma)$.
Lemma~\ref{lem:thm}~1) assures that there exists a unique sequence
$\alpha_0,\alpha_{j_1},\ldots,\alpha_{j_{m-3}}$ in $\tilde\Psi$
and an element $w\in W_\Sigma$ such that
\begin{equation}\label{eq:AinD}
  w\circ\iota(\beta_\nu) = \alpha_{j_\nu}\quad(\nu=0,\ldots,m-3)
  \text{ with }j_0 = 0.  
\end{equation}
Putting
\[
 \begin{split}
  \Psi' &= \bigl\{\alpha\in\Psi,\;\,(\alpha|\alpha_{j_\nu})=0
  \quad(\nu=0,\ldots,m-4)\bigr\},\\
  \Theta &= \bigl\{\alpha\in\Psi'\,;\,(\alpha,\alpha_{j_{m-3}}))<0\bigr\},\\
  \alpha'&=\alpha_{j_{m-3}},\\
  (\beta,\beta') &= \bigl(w\circ \iota(\beta_{m-2}),w\circ\iota(\beta_{m-1})\bigr),
 \end{split}\quad
 \begin{xy}
 \ar@{.} (-3,0);(0,0) *{}
 \ar@{-} (0,0);(5,0) *+!D{\!\!\beta_{m-4}} *{\circ}="A"
 \ar@{-} "A";(13,0) *+!D{\beta_{m-3}} *{\circ}="B"
 \ar@{-} "B";(21,0) *+!D{\ \beta_{m-2}} *{\circ}
 \ar@{-} "B";(13,-7) *+!L{\beta_{m-1}} *{\circ}
 \end{xy}
\]
we have $\beta$, $\beta'\in\langle\Psi'\rangle$ and we can apply 
Lemma~\ref{lem:A2} as in the case when $\Xi$ is of type $A$.
Thus
\begin{align*}
  &\# W_\Sigma\backslash\!\bigl\{\iota\in\Hom'(\Xi,\Sigma)\,;\,
   \exists w\in W_\Sigma\text{ such that }\eqref{eq:AinD}
   \text{ is satisfied.}\bigr\}\\
  &= \bigl(\#(\Theta\cap\Sigma^L)\bigr)\bigl(\#(\Theta\cap\Sigma^L)-1\bigr)
  + \#\bigl\{\alpha\in\Theta\cap\Sigma^L\,:\, \text{the irreducible component
    of }\\
  &\quad \Psi' \text{ containing $\alpha$ 
    is not of type $A$ or $\alpha$ is not an end vertex of the component}
  \bigr\}.
\end{align*}
Hence $\Hom'(D_m,\Sigma)=\emptyset$ if $\Sigma$ is of type $A_n$, $C_n$ or
$G_2$ or $m>\rank\Sigma$.
Moreover we have $\#\overline\Hom(D_m,\Sigma)$ shown in the following
table under the notation in \S\ref{sec:table}.

\smallskip
\smallskip
%\centerline{
%\begin{longtable}{|l|c|l|c|}\hline
\begin{tabular}{|l|c|l|c|}\hline
$\Sigma$ & $\Xi$ & $\Psi'$ & $ \#$ \\ \hline
$D_4$     & $D_4$ & $\{\alpha_1,\alpha_3,\alpha_4\}\simeq 3A_1$
 & 6 \\ \hline
$D_5$     & $D_4$ & $\Psi\setminus\{\alpha_2\}\simeq A_1+A_3
 (\ni\alpha_3:\text{not an end root})$
 & 3\\ \hline
$D_{n\ (n\ge6)}$ & $D_4$ & $\Psi\setminus\{\alpha_2\}\simeq A_1+D_{n-2}$
 & 3 \\ \hline 
$B_4$          & $D_4$ & $\Psi\setminus\{\alpha_2\}\simeq A_1+B_2$
 & 3\\ \hline
$F_4$ & $D_4$ & $\Psi\setminus\{\alpha_1\}\simeq C_3$
 & 1\\ \hline
$D_{n\ (4<m=n)}$ & $D_m$ & $\{\alpha_{n-1},\alpha_n\}\simeq 2A_1$
 & 2\\ \hline 
$D_{n\ (4<m=n-1)}$ & $D_m$ & $\{\alpha_{n-2},\alpha_{n-1},\alpha_n\}
  \simeq A_3(\ni\alpha_{n-2}:\text{not an end root})$
 & 1\\ \hline
$D_{n\ (4<m\le n-2)}$ & $D_m$ & $\{\alpha_{m-1},\ldots,\alpha_n\}
 \simeq D_{n-m+1}$
 & 1\\ \hline
$B_{n\ (4<m\le n)}$ & $D_m$ & $\{\alpha_{m-1},\ldots,\alpha_n\}
 \simeq B_{n-m+1}$
 & 1\\ \hline
$E_6$ & $D_4$  & $\Psi\setminus\{\alpha_2\}\simeq A_5
 (\ni\alpha_4: \text{not an end root})$
 & 1\\ \cline{2-4}
      & $D_5$ & $\{\alpha_1,\alpha_3,\alpha_5,\alpha_6\} \simeq 2A_2$
 & 2\\ \hline
$E_7$ & $D_4$ & $\Psi\setminus\{\alpha_1\}\simeq D_6$
 & 1\\ \cline{2-4}
      & $D_5$ & $\{\alpha_2,\alpha_4,\ldots,\alpha_7\}
 \simeq A_5(\ni\alpha_4:\text{not an end root})$
 & 1\\ \cline{2-4}
      & $D_6$ & $\{\alpha_2,\alpha_5,\alpha_6,\alpha_7\}
 \simeq A_1+A_3(\ni\alpha_5:\text{an end root})$
 & 2\\ \hline
$E_8$ & $D_4$ & $\Psi\setminus\{\alpha_8\}\simeq E_7$
 & 1\\ \cline{2-4}
    & $D_5$ & $\{\alpha_1,\ldots,\alpha_6\}\simeq E_6$
 & 1\\ \cline{2-4}
    & $D_6$ & $\{\alpha_1,\ldots,\alpha_5\}\simeq D_5$
 & 1\\ \cline{2-4}
    & $D_7$ & $\{\alpha_1,\alpha_2,\alpha_3,\alpha_4\}\simeq A_4
 (\ni\alpha_4:\text{not an end root})$
 & 1\\ \cline{2-4}
    & $D_8$ & $\{\alpha_2,\alpha_1,\alpha_3\}\simeq A_1+A_2$
 & 2\\ \hline
\end{tabular}

\smallskip
Here $m_\Sigma$ is the rank of the maximal subdiagram of type $D_m$ contained 
in the extended Dynkin diagram of $\Sigma^L$ and then
\begin{equation}
  m_\Sigma = 
     \begin{cases}
      n&(\Sigma\text{ is of type $B_n$ or $D_n$}),\\
      5,\,6,\,8 &(\Sigma
       \text{ is of type $E_6$, $E_7$ or $E_8$, respectively}).
     \end{cases}
\end{equation}
%%%%%%%%%%%%%%%% D_n in \Sigma %%%%%%%%%%%%%%%
\begin{equation}\label{eq:DinE}
\begin{gathered}
\begin{xy}
 \ar@{.}  *+!D{\alpha_1} *{\circ};
  (5,0)  *+!D{\alpha_2} *{\circ}="A"
 \ar@{-} "A";(10,0) *+!D{\alpha_3} *{\circ}="C"
 \ar@{-} "C";(15,0) *+!D{\alpha_4} *{\circ}="D"
 \ar@{-} "D";(18,0)  \ar@{.} (18,0);(22,0)
 \ar@{-} (22,0);(25,0) *{\circ}="F"
 \ar@{-} "F";(30,0) *+!D{\alpha_{n-2}} *{\circ}="G"
 \ar@{-} "G";(35,0) *+!LD{\!\!\alpha_{n-1}} *{\circ}
 \ar@{-} "A";(5,-5) *+!L{\alpha_0} *{\ccirc}
 \ar@{-} "G";(30,-5) *+!L{\alpha_n} *{\circ}
\end{xy}\qquad
\begin{xy}
 \ar@{-}  *+!D{\alpha_1} *{\circ};
  (5,0)  *+!D{\alpha_2} *{\circ}="A"
 \ar@{-} "A";(10,0) *+!D{\alpha_3} *{\circ}="C"
 \ar@{-} "C";(15,0) *+!D{\alpha_4} *{\circ}="D"
 \ar@{-} "D";(18,0)  \ar@{.} (18,0);(22,0)
 \ar@{-} (22,0);(25,0) *{\circ}="F"
 \ar@{-} "F";(30,0) *{\circ}="G"
 \ar@{-} "G";(35,0) *+!D{\!\alpha_{n-1}} *{\circ}="H"
 \ar@{-} "A";(5,-5) *+!L{\alpha_0} *{\ccirc}
 \ar@2{.>} "H";(40,0) *+!LD{\!\!\alpha_n} *{\circ}
\end{xy}\\[5pt]
\begin{xy}
 \ar@{.}  *+!D{\alpha_1} *{\ast} ;
  (5,0)  *+!D{\alpha_3} *{\circ}="C"
 \ar@{-} "C";(10,0) *+!D{\alpha_4} *{\circ}="D"
 \ar@{-} "D";(15,0) *+!D{\alpha_5} *{\circ}="E"
 \ar@{.} "E";(20,0) *+!D{\alpha_6} *{\ast},
 \ar@{-} "D";(10,-5) *+!L{\alpha_2} *{\circ}="X"
 \ar@{-} "X";(10,-10) *+!L{\alpha_0} *{\ccirc};
\end{xy}\quad
\begin{xy}
 \ar@{-}  *+!D{\alpha_0} *{\ccirc};
  (5,0)  *+!D{\alpha_1} *{\circ}="A"
 \ar@{-} "A";(10,0) *+!D{\alpha_3} *{\circ}="C"
 \ar@{-} "C";(15,0) *+!D{\alpha_4} *{\circ}="D"
 \ar@{-} "D";(20,0) *+!D{\alpha_5} *{\circ}="E"
 \ar@{.} "E";(25,0) *+!D{\alpha_6} *{\ast}="F"
 \ar@{.} "F";(30,0) *+!D{\alpha_7} *{\bullet},
 \ar@{-} "D";(15,-5) *+!L{\alpha_2} *{\circ}
\end{xy}\quad
\begin{xy}
 \ar@{.}  *+!D{\alpha_1} *{\ast};
  (5,0)  *+!D{\alpha_3} *{\circ}="A"
 \ar@{-} "A";(10,0) *+!D{\alpha_4} *{\circ}="C"
 \ar@{-} "C";(15,0) *+!D{\alpha_5} *{\circ}="D"
 \ar@{-} "D";(20,0) *+!D{\alpha_6} *{\circ}="E"
 \ar@{-} "E";(25,0) *+!D{\alpha_7} *{\circ}="F"
 \ar@{-} "F";(30,0) *+!D{\alpha_8} *{\circ}="G"
 \ar@{-} "G";(35,0) *+!D{\alpha_0} *{\ccirc}
 \ar@{-} "C";(10,-5) *+!L{\alpha_2} *{\circ}
\end{xy}
\end{gathered}
\end{equation}
%%%%%%%%%%%%%%%%%%%%%%%%%%%%%%%%%%%%%%

\smallskip
Fix $\iota\in\Hom(D_4,F_4)$.
Since $\Hom(D_4,F_4)$ is a single $W_{F_4}$-orbit, for any $g\in\Aut(D_4)$ 
there exists $w_g\in W_{F_4}$ with $\iota\circ g=w_g\circ\iota$.
Here $w_g$ is uniquely determined by $g$ because $\rank F_4=\rank D_4$.
Hence we have
\begin{equation}\label{eq:DF4}
 \begin{split}
  \Aut(D_4)&\simeq W_{F_4}\supset \Aut(B_4)=W_{B_4} \supset W_{D_4},\\
  \Out(D_4)&\simeq\mathfrak S_3,\quad W_{B_4}/W_{D_4}
   \simeq\mathbb Z/2\mathbb Z.
 \end{split}
\end{equation}
Let $\iota:D_4\subset D_n$ $(\subset B_n)$ be the natural imbedding 
given by the realization in \S\ref{sec:Dynkin} and let $g\in\Aut(D_4)$
be a non-trivial rotation of $G(D_4)$.  Then it is easy to see $\iota\circ g\ne w\circ\iota$ for any $w\in W_{B_n}$.  Hence if $n\ge4$, we have
\begin{equation}
  \#\bigl(\overline\Hom(D_4,D_n)/\!\Out(D_4)\bigr) =
  \#\bigl(\overline\Hom(D_4,B_n)/\!\Out(D_4)\bigr) = 1
\end{equation}
because $\#\bigl(\overline\Hom(D_4,D_n)\bigr)
=\#\bigl(\overline\Hom(D_4,B_n)\bigr)
=3$ and moreover we have
\begin{equation}\label{eq:DDperp}
 D_m^\perp\cap D_n\simeq D_{n-m},\quad
 D_m^\perp\cap B_n\simeq B_{n-m}
\end{equation}
for $m=4$.  
Here $D_0=D_1=B_0=\emptyset$, $D_2\simeq 2A_1$, $B_1\simeq A_1^S$ and
$A_1^S$ is the root space of type $A_1$ such that 
$A_1^S\cap (B_n)^L=\emptyset$.

Note that
\begin{equation}
  \Aut(D_n)\simeq W_{B_n}\text{ and }
  \Out(D_n):=\Aut(D_n)/W_{D_n} \simeq \mathbb Z/2\mathbb Z
  \quad(n\ge 5)
\end{equation}
under the natural imbedding $D_n\subset B_n$ of root spaces.
Thus we have
\begin{equation}
 \#\bigl(\overline\Hom(D_m,\Sigma)/\!\Out(D_m)\bigr) = 0\text{ or }
 1\text{ if $m\ge4$ and $\Sigma$ is irreducible}
\end{equation}
and therefore $\{\iota(D_m)^\perp\,;\,\iota\in\Hom(D_m,\Sigma)\}$ 
is a single $W_{\Sigma}$-orbit if it is non-empty.

Thus we have \eqref{eq:DDperp} for $4\le m\le n$ since it does not
depend on the imbedding of $D_m$.

Let $n\in\{6,7,8\}$ and put $m=m_{E_n}$.
There exists $\iota\in\Hom(D_m, E_n)$ %for $m=m_{E_n}$ 
such that $\iota$ corresponds to the imbedding of 
$\Phi_{m_\Sigma}$ to $\tilde\Psi$ with $\iota(\beta_0)=\alpha_0$.
%Here we may assume $\beta_1$, $\beta_{m-1}$, $\beta_m$ is the end roots
%of $\Phi$ and $(\beta_{m-2},\beta_m)<0$.
Then we have
$D_5^\perp\cap E_6=\emptyset$, $D_6^\perp\cap E_7\simeq A_1$ and
$D_8^\perp\cap E_8=\emptyset$ from Lemma~\ref{lem:orth}.

Moreover there exists $\iota'\in\Hom(D_{m-1}, E_n)$ such that
\[
  \iota'(D_{m-1}) = \bigl\{\iota(\beta_0),\ldots,\iota(\beta_{m-3}),
  \iota(\beta_{m-3})+\iota(\beta_{m-2})+\iota(\beta_{m-1})\bigr\}
\]
and it is clear that $D_{m-1}^\perp\cap E_n\simeq D_m^\perp\cap E_n$.

Let $4\le k\le 6$.
Then $D_k^\perp\cap E_8\supset D_k^\perp\cap D_8\simeq D_{8-k}$
and we can conclude $D_k^\perp\cap E_8\simeq D_{8-k}$
because $\rank(D_k^\perp\cap E_8)\le 8-k$ and there is no root system
containing $D_{8-k}$ as a proper subsystem such that its roots have the same
length and its rank is not larger than $8-k$.

Since $D_6^\perp\cap E_7\simeq A_1$ and $D_4^\perp\cap D_6\simeq 2A_1$,
$D_4^\perp\cap E_7\supset 3A_1$ and we have $D_4^\perp\cap E_7\simeq 3A_1$
by the same argument as above.

Thus we have obtained the claims in the lemma and therefore 
Remark ~\ref{rem:thmauto} i) is also clear.

%%%%%  B_m %%%
\medskip
3) Suppose $\Xi$ is of type $B_m$ with $m\ge 2$.

Note that for any $\beta\in\Xi\setminus\Xi^L$, there exists
$\beta_1,\,\beta_2\in\Xi^L$ such that
$\beta=\frac12(\beta_1+\beta_2)$ and $(\beta_1|\beta_2)=0$.
Hence $\iota\in\Hom(\Xi,\Sigma)$ is determined by $\iota|_{\Xi^L}$.
Note that $\Xi^L$ is of type $D_m$ with 
$D_2\simeq 2A_1$ and $D_3\simeq A_3$.

Then $\Hom(\Xi,\Sigma)\ne\emptyset$ means $\Sigma$ is of type
$B_n\ (n\ge m)$ or $F_4$ if $m>2$.

If $m>2$ or if $\Sigma$ is of type $F_4$, $\#\overline\Hom(\Xi^L,\Sigma)=1$ 
and therefore $\#\overline\Hom(\Xi,\Sigma)=1$.
If $m=2$ and $\Sigma=B_n$ or $C_n$, it is easy to see that
\[
  \bigl\{\iota\in\Hom(2A_1,\Sigma)\,;\,
    \tfrac12\bigl(\iota(\beta_1)+\iota(\beta_2)\bigr)\in \Sigma
  \bigr\}
\]
is a single $W_{\Sigma}$-orbit and we have 
also $\#\overline\Hom(\Xi,\Sigma)=1$.
Here $\Xi^L = \langle\beta_1\rangle+\langle\beta_2\rangle\simeq 2A_1$.

\medskip
4) When $\Sigma$ is of type $C_n$, we have the lemma from the case 3) 
by considering the dual root systems $\Xi^\vee$ and $\Sigma^\vee$.

%%%%%  E_6, E_7 in E_8 %%%
\medskip
5) We first examine $\overline\Hom(E_6,E_8)$ and $\overline\Hom(E_7,E_8)$
under the notation in \S\ref{sec:Dynkin}.\\
Since $\#\overline\Hom(A_5,E_8)=\#\overline\Hom(A_6,E_8)=1$,
we may assume
\[
E_6^o\supset \Psi_{\!A_5}
  =\{\alpha_0=-\epsilon_7-\epsilon_8,\,\alpha_8=\epsilon_7-\epsilon_6,\,
  \alpha_7=\epsilon_6-\epsilon_5,\,\alpha_6=\epsilon_5-\epsilon_4,\,\alpha_5=\epsilon_4-\epsilon_3\}
\]
for the imbedding $E_6\simeq E_6^o\subset E_8$.
Let $\tilde\alpha\in\Phi\setminus\Psi_{\!A_5}$.
We have
\[
\tilde\alpha=\sum_{j=1}^8c_j\epsilon_j\in E_6^o\subset E_8:\ 
\langle\tilde\alpha,\alpha_j\rangle
=\begin{cases}
   0&(j=0,\,8,\,6,\,5),\\
   -1&(j=7).
 \end{cases}
\]
Thus 
\[
  \tilde\alpha = c_1\epsilon_1+c_2\epsilon_2+c(\epsilon_3+\epsilon_4+\epsilon_5) + (c-1)(\epsilon_6+\epsilon_7-\epsilon_8).
\]
Since $\tilde\alpha$ is a root of $E_8$,
we have $c=\frac12$ and hence
\[
  \tilde\alpha = \alpha_{\pm} := \tfrac12(\pm(\epsilon_1+\epsilon_2)+\epsilon_3+\epsilon_4+\epsilon_5-\epsilon_6-\epsilon_7+\epsilon_8).
\]
Since $\alpha_2=\epsilon_1+\epsilon_2$ is orthogonal to $\alpha_j$ ($j=0,5,6,7,8$) and
$s_{\alpha_2}\alpha_+=\alpha_-$, 
we have $\#\overline\Hom(E_6,E_8)=1$ and 
\[
 \begin{split}
  (E_6^o)^\perp\cap E_8 &
  \simeq \Psi_{\!A_5}^\perp\cap\alpha_+^\perp\cap E_8\\
  &= (\langle \alpha_1,\alpha_3\rangle + \langle\alpha_2\rangle)
     \cap\alpha_+^\perp\\
  &=\langle \alpha_1,\alpha_3\rangle\simeq A_2.
 \end{split}\qquad
\raisebox{10pt}{\begin{xy}
 \ar@{.}  *+!D{\alpha_1} *{\ast};
  (5,0)  *+!D{\alpha_3} *{\ast}="A"
 \ar@{.} "A";(10,0) *+!D{\alpha_4} *{\ast}="C"
 \ar@{.} "C";(15,0) *+!D{\alpha_5} *{\circ}="D"
 \ar@{-} "D";(20,0) *+!D{\alpha_6} *{\circ}="E"
 \ar@{-} "E";(25,0) *+!D{\alpha_7} *{\circ}="F"
 \ar@{-} "F";(30,0) *+!D{\alpha_8} *{\circ}="G"
 \ar@{-} "G";(35,0) *+!D{\alpha_0} *{\ccirc}
 \ar@{.} "C";(10,-5) *+!L{\alpha_2} *{\ast}
 \ar@{.} "F";(25,-5) *+!L{\tilde\alpha} *{\circ}
\end{xy}}
\]
Let $E_7\simeq E_7^o\subset E_8$.
Then we may moreover assume $\alpha_4=\epsilon_3-\epsilon_2\in E_7$ and the condition
$\tilde\alpha \perp \alpha_4$ implies $\tilde\alpha=\alpha_+$.
Hence $\#\overline\Hom(E_7,E_8)=1$ and 
\[
 (E_7^o)^\perp\cap E_8\simeq
 \langle \alpha_1,\alpha_3\rangle\cap\alpha_4^\perp
 =\langle\alpha_1\rangle\simeq A_1.
\]

Now we examine $\overline\Hom(E_6,E_7)$.
Since $A_5\subset E_6\simeq E_6^o\subset E_7$, 
the argument in 1) assures that we may assume
\begin{align*}
 E_6^o\supset\Psi_{\!A_5}'&:=\Psi_{\!A_4}\cup\{\alpha_2=\epsilon_1+\epsilon_2\}
 \text{ or }
 E_6^o\supset\Psi_{\!A_5}:=\Psi_{\!A_4}\cup\{\alpha_5=\epsilon_4-\epsilon_3\}\\
 \Psi_{\!A_4}&:=\{\alpha_0=\epsilon_8-\epsilon_7,\, 
 \alpha_1 = \tfrac12(\epsilon_1-\epsilon_2-\epsilon_3-\epsilon_4-\epsilon_5-\epsilon_6-\epsilon_7+\epsilon_8),\\ 
 &\quad\quad \alpha_3=\epsilon_2-\epsilon_1,\, \alpha_4 = \epsilon_3-\epsilon_2\}.
\end{align*}
Then there exists 
$\tilde\alpha=\sum_{j=1}^8c_j\epsilon_j\in E_6^o\subset E_7$ such that
\[
\begin{split}
 (\tilde\alpha|\alpha_j) &= 0\quad(j=0,\,1,\,4),\\
 (\tilde\alpha|\alpha_3) &=-1,\\
 (\tilde\alpha|\alpha_2)(\tilde\alpha|\alpha_5)&=0.
\end{split}
\qquad
\raisebox{10pt}{
\begin{xy}
 \ar@{-}  *+!D{\alpha_0} *{\ccirc};
  (5,0)  *+!D{\alpha_1} *{\circ}="A"
 \ar@{-} "A";(10,0) *+!D{\alpha_3} *{\circ}="C"
 \ar@{-} "C";(15,0) *+!D{\alpha_4} *{\circ}="D"
 \ar@{-} "D";(20,0) *+!D{\alpha_5} *{\circ}="E"
 \ar@{.}  "E";(25,0) *+!D{\alpha_6} *{\ast}="F"
 \ar@{.}  "F";(30,0) *+!D{\alpha_7} *{\ast},
 \ar@{.} "D";(15,-5) *+!L{\alpha_2} *{\circ}
 \ar@{.} "C";(10,-5) *+!L{\tilde\alpha} *{\circ}
\end{xy}}
\]
Then the condition $(\tilde\alpha|\alpha_0)=0$ implies $c_7=c_8=0$
and
\begin{gather*}
\tilde\alpha 
 = (c+1)\epsilon_1+c(\epsilon_2+\epsilon_3)+c_4\epsilon_4+c_5\epsilon_5+c_6\epsilon_6,\\
 1-c-c_4-c_5-c_6 =0,\\
 (2c+1)(c-c_4)=0.
\end{gather*}
Hence $c=0$, $E_6^0\supset\Psi_{\!A_5}$ and 
$\tilde\alpha=\epsilon_1+\epsilon_5$ or $\epsilon_1+\epsilon_6$.
Since
\[
 \Psi_{\!A_5}^\perp\cap E_7 = \langle\alpha_7\rangle = \langle \epsilon_6-\epsilon_5\rangle
\]
and $s_{\epsilon_6-\epsilon_5}(\epsilon_1+\epsilon_5)=\epsilon_1+\epsilon_6$, we have
$\#\overline\Hom(\Xi,\Sigma)=1$ and $(E_6^o)^\perp\cap E_7=\emptyset$.

If $\#\Hom(\Xi,\Sigma)=1$, any element of $\Hom(\Xi,\Sigma)$ is 
isomorphic to the imbedding $\iota$ corresponding to the 
graphic imbedding $\bar\iota$ given in the claim. 
Since $\iota(\Xi)^\perp\supset\langle\Psi\cap\bar\iota(\Xi)^\perp\rangle$
and $\iota(\Xi)^\perp\simeq\langle\Psi\cap\bar\iota(\Xi)^\perp\rangle$, 
we have $\iota(\Xi)^\perp=\langle\Psi\cap\bar\iota(\Xi)^\perp\rangle$.

6)
If $\Xi$ is of type $G_2$ or $F_4$, the lemma is clear 
and thus we have completed the proof of the lemma.
%%%%%%%%%%%%%%%%%%%%%%%%%%%%

\section{Dual pairs and closures}\label{sec:dual}
\begin{defn}\label{def:OutX}
For a subsystem $\Xi$ of a root system $\Sigma$ and a subgroup $G$ of
$\Aut(\Sigma)$ we put
\begin{align}
 N_G(\Xi) &:= \{g\in G\,;\,g(\Xi)=\Xi\},\quad 
 Z_G(\Xi) := \{g\in G\,;\,g|_\Xi=id\},\\
 \Aut_{\Sigma}(\Xi) &:= 
   N_{W_\Sigma}(\Xi)/
   Z_{W_\Sigma}(\Xi)
 \subset \Aut(\Xi),\\
 \Out_{\Sigma}(\Xi) &:=\Aut_{\Sigma}(\Xi)/W_\Xi
 \simeq N_{W_\Sigma}(\Xi)
 /(W_\Xi\times W_{\Xi^\perp})\subset \Out(\Xi).\label{eq:OutX}
\end{align}
Note that the isomorphism in \eqref{eq:OutX} follows from the equality
$Z_{W_\Sigma}(\Xi) = W_{\Xi^\perp}$.
\end{defn}
\begin{prop}\label{prop:dual}
Let $\Xi_1$ be a subsystem of $\Sigma$. 
Put $\Xi_2=\Xi_1^\perp$.
Then there is a homomorphism
\begin{align}
 \varpi\,:\,\Out_{\Sigma}(\Xi_1)
 &\to \Out_{\Sigma}(\Xi_2)
 \simeq\Out_\Sigma(\Xi_2^\perp)\label{eq:D-hom}\allowdisplaybreaks
\intertext{and}\allowdisplaybreaks
 \varpi\text{ is bijective}&\text{ if \ }\Xi_2^\perp = \Xi_1,
  \label{eq:D-inj}\allowdisplaybreaks\\
 \Out_{\Sigma}(\Xi_1)\overset{\sim}{\to}
 \Out(\Xi_1)&\text{ if \ }
   \#\bigl(\overline\Hom(\Xi_1,\Sigma)/\!\Out(\Xi_1)\bigr)
  =\#\overline\Hom(\Xi_1,\Sigma),\label{eq:D-all}\allowdisplaybreaks\\
 \Out_{\Sigma}(\Xi_2)\overset{\sim}{\to}
 \Out(\Xi_2)&\text{ if \ }
   \#\bigl(\overline\Hom(\Xi_2,\Sigma)/\!\Out(\Xi_2)\bigr)
  =\#\overline\Hom(\Xi_2,\Sigma).\label{eq:D-surj}
\end{align}
\end{prop}
\begin{proof}
Since $N_{W_\Sigma}(\Xi_1)\subset N_{W_\Sigma}(\Xi_2)$
and $\Xi_2^\perp\supset\Xi_1$, 
\eqref{eq:D-hom} is well-defined
and \eqref{eq:D-inj} is clear.
Suppose $\#\bigl(\overline\Hom(\Xi_1,\Sigma)/\!\Out(\Xi_1)\bigr)
  =\#\overline\Hom(\Xi_1,\Sigma)$.
Then for any $g\in\Aut(\Xi_1)$ there exists $w\in W_\Sigma$ with
$w|_{\Xi_1}=g|_{\Xi_1}$ and \eqref{eq:D-all} is clear.
We similarly have \eqref{eq:D-surj}.
The isomorphism in \eqref{eq:D-hom} follows from 
\eqref{eq:D-inj} and the relation $(\Xi_2^\perp)^\perp=\Xi_2$.
\end{proof}
\begin{defn}[dual pairs]\label{def:dual}
A pair $(\Xi_1,\Xi_2)$ of subsystems of a root system $\Sigma$ is called
a \textsl{dual pair} in $\Sigma$ if 
\begin{equation}
  \Xi_1^\perp=\Xi_2\quad\text{and}\quad\Xi_2^\perp=\Xi_1.
\end{equation}
If $(\Xi_1,\Xi_2)$ is a dual pair, the map $\varpi$ in 
Proposition~\ref{prop:dual} is an isomorphism.
The dual pair is called \textsl{special}
if the map $\varpi$ is the isomorphism
\begin{equation}\label{eq:special}
  \varpi\,:\,\Out(\Sigma_1)\overset{\sim}{\to}\Out(\Sigma_2).
\end{equation}
For a subsystem $\Xi$ of $\Sigma$, its 
$\perp$-\textsl{closure} $\overline\Xi$ is defined by 
$\overline\Xi:=(\Xi^\perp)^\perp$.
Then $(\Xi,\Xi^\perp)$ is a dual pair if and only if $\Xi$ is $\perp$-closed
(i.e.\ $(\Xi^\perp)^\perp=\Xi$)
and hence $(\overline\Xi,\Xi^\perp)$ is always a dual pair.
We say that $\Xi$ is \textsl{$\perp$-dense} in $\Sigma$ if 
$\Xi^\perp=\emptyset$.
\end{defn}
\begin{cor}
Let $(\Xi_1,\Xi_2)$ be a dual pair in $\Sigma$. Then
\begin{equation}\label{eq:D-bad}
 \Out(\Xi_1)\not\simeq\Out(\Xi_2)\ \Rightarrow\ 
  \begin{cases}
   \#\bigl(\Hom(\Xi_1,\Sigma)/\!\Out(\Xi_1)\bigr)
    <\#\overline\Hom(\Xi_1,\Sigma)\\
    \text{or}\\
   \#\bigl(\Hom(\Xi_2,\Sigma)/\!\Out(\Xi_2)\bigr)
    <\#\overline\Hom(\Xi_2,\Sigma).
  \end{cases}
\end{equation}
Suppose $\#\overline\Hom(\Xi_2,\Sigma)=1$.
Let $\iota\in\Hom(\Xi_1,\Sigma)$.
Then we have
\begin{gather}
 (\Xi_1,\Xi_2)\text{ is a special dual pair }\Leftrightarrow\ 
 \#\Out(\Xi_1)=\#\Out(\Xi_2),\label{eq:dual1}\\
    \text{$\exists w\in W_\Sigma$ such that }
    \iota(\Xi_1)=w(\Xi_1)\ \Leftrightarrow\  
    \iota(\Xi_1)^\perp\simeq\Xi_2.\label{eq:dual2}
\end{gather}
\end{cor}
\begin{proof}
Note that \eqref{eq:D-bad} is the direct consequence of 
Proposition~\ref{prop:dual}.

Suppose $\#\overline\Hom(\Xi_2,\Sigma)=1$.
Then Proposition~\ref{prop:dual} implies
\[
  \Out(\Xi_1)\supset\Out_\Sigma(\Xi_1)\overset{\sim}\to\Out(\Xi_2)
\]
and \eqref{eq:dual1} is clear.
%Fix $\Xi_2\subset\Sigma$.
Then if $\iota(\Xi_1)^\perp\simeq\Xi_2$, there exists $w\in W_\Sigma$
with $\iota(\Xi_1^\perp)=w(\Xi_2)$ and therefore $\iota(\Xi_1)=w(\Xi_1)$,
which implies the claim.
\end{proof}
\begin{exmp} {\rm i)} \ 
The followings are examples of the triplets $(\Sigma,\Xi_1,\Xi_2)$
such that $(\Xi_1,\Xi_2)$ are special dual pairs in $\Sigma$.
\begin{gather*}
 (D_{m+n}, D_m,D_n)\quad(m\ge2,\,n\ge2,\,m\ne4,\,n\ne4),\\
 (E_6, A_3, 2A_1), \ 
 (E_7,A_5, A_2),\ 
 (E_7, A_3+A_1, A_3),\ 
 (E_7, 3A_1, D_4), \allowdisplaybreaks\\
 (E_8, E_6, A_2),\ 
 (E_8, A_5, A_2+A_1),\  
 (E_8, A_4, A_4),\ 
 (E_8, D_6, 2A_1),\ 
 (E_8, D_5, A_3),\allowdisplaybreaks\\
 (E_8, D_4, D_4),\ 
 (E_8, D_4+A_1, 3A_1),\ 
 (E_8, 2A_2, 2A_2),\allowdisplaybreaks\\
 (E_8, A_3+A_1,A_3+A_1),\ (E_8,4A_1,4A_1),\ 
 (F_4,A_2,A_2).
\end{gather*}
In these examples except for
$(D_4,2A_1,2A_1)$ and $(E_8,4A_1,4A_1)$, 
$\#\overline\Hom(\Xi_2,\Sigma)=1$ and 
the triplet is uniquely determined by the data 
$(\Sigma,\Xi_1,\Xi_2)$ up to the automorphisms defined 
by $W_\Sigma$.  If the imbedding $4A_1\subset E_8$ satisfies
$(4A_1)^\perp \simeq 4A_1$, we have a special dual pair
$(4A_1,4A_1)$ in $E_8$, which is also uniquely defined.
The imbedding $2A_1\subset D_4$ is unique up to $\Aut(D_4)$.

{\rm ii)}
The isomorphism $\varpi\in\Out(2A_2)$ defined by the dual 
pair $(2A_2,2A_2)$ in $E_8$ satisfies
\begin{equation}
 \varpi\bigl(\Out(A_2)\times\Out(A_2)\bigr)
 \ne\Out(A_2)\times\Out(A_2)
\end{equation}
because
$\#\bigl(\overline\Hom(4A_2,E_8)/\Out'(4A_2)\bigr)=
\#\overline\Hom(2A_2,E_8) = 1$.
%If $\Out(2A_2)$ is identified with $W_{B_2}$, $\varpi$
%corresponds to an outer automorphism of $W_{B_2}$.
See \S\ref{sec:4A2E8}.

{\rm iii)}
It happens that any dual pair of $E_8$ is special.
But for example, if $(\Sigma,\Xi_1,\Xi_2)$ is $(E_6,A_5,A_1)$ or
$(E_7,D_6,A_1)$, $(\Xi_1,\Xi_2)$ is a dual pair in $\Sigma$
satisfying $\Out(\Xi_1)\not\simeq\Out(\Xi_2)$ and 
$\#\overline\Hom(\Xi_2,\Sigma)=1$, which implies
$\#\overline\Hom(\Xi_1,\Sigma)>1$.
\end{exmp}

\begin{defn}[$S$-closure and $L$-closure]\label{defn:S-closed}
Let $\Xi$ be a subsystem of $\Sigma$.
Then $\Xi$ is \textsl{$S$-closed} if and only if
\begin{equation}\label{eq:S-cl}
 \alpha,\,\beta\in\Xi\text{ \ and \ }\alpha+\beta\in\Sigma
 \ \Rightarrow \ \alpha+\beta\in\Xi
\end{equation}
and \textsl{$L$-closed} if and only if
\begin{equation}\label{eq:L-cl}
 \beta\in\Sigma\cap\sum_{\alpha\in\Xi}\mathbb R\alpha
 \ \Rightarrow \ \beta\in\Xi.
\end{equation}
The smallest $S$-closed (resp.~$L$-closed) subsystem of 
$\Sigma$ containing $\Xi$ is called the $S$-closure 
(resp.~$L$-closure) of $\Xi$.
\end{defn}
\begin{rem}\label{rem:closed}
i) We have the following relation 
for a subsystem $\Xi$ of $\Sigma$:
\begin{equation}
 \text{$\perp$-closed}\ \Rightarrow\ 
 \text{$L$-closed}\ \Rightarrow
 \text{$S$-closed}.
\end{equation}

ii) Let $\mathfrak g$ be a complex semisimple Lie algebra with the 
root system $\Sigma$ and let $X_\alpha$ be root vectors corresponding to
$\alpha\in\Sigma$. Then the root system of the semisimple Lie algebra
$\mathfrak g_\Xi$ 
generated by $\{X_\alpha\,;\,\alpha\in\Xi\}$ is the $S$-closure of
$\Xi$.

Let $\Xi_1$ and $\Xi_2$ be  $S$-closed subsystems of $\Sigma$.
Then
\begin{equation}
 [\mathfrak g_{\Xi_1},\mathfrak g_{\Xi_2}] = 0
 \ \Leftrightarrow\ 
 \Xi_1\perp\Xi_2.
\end{equation}
Hence if $(\Xi_1,\Xi_2)$ is a dual pair 
with $\rank\Xi_1+\rank\Xi_2=\rank\Sigma$,
the dual pair of root systems gives a dual pair
in semisimple Lie algebras (cf.~\cite{Ru}).

iii)
Suppose $\Sigma$ is irreducible and 
there exist $\alpha,\,\beta\in\Sigma$ with 
$\alpha+\beta\in\Sigma\setminus\Xi$.
Then $\langle\alpha,\beta,\alpha+\beta\rangle$ is of type
$B_2$ or of type $G_2$, which implies that
$\Sigma$ is not simply laced. 
For example, $D_n\subset C_n$ is not $S$-closed and the
$S$-closure of $D_n$ equals $C_n$ ($n\ge 2$).

iv) Let $\Xi$ be an $L$-closed subsystem of $\Sigma$.
Then for any subsystem $\Xi'$ of $\Sigma$ 
\begin{equation}\label{eq:W2}
  W_{\Xi}\cap W_{\Xi'}=W_{\Xi\cap\Xi'}.
\end{equation}
This is proved as follows.
Choose a generic element $v$ of the orthogonal complement
of $\sum_{\alpha\in\Xi}\mathbb R\alpha$ in 
$\sum_{\alpha\in\Sigma}\mathbb R\alpha$ 
so that $\{\alpha\in\Sigma\,;\,(\alpha|v)=0\}=\Xi$.
Since $W_{\Xi}=\{w\in W_\Sigma\,;\,wv=v\}$, 
$W_{\Xi}\cap W_{\Xi'} =
  \{w\in W_{\Xi'}\,;\,wv=v\}
 =W_{\{\alpha\in\Xi'\,;\,(\alpha|v)=0\}}
 =W_{\Xi\cap\Xi'}$.

v) Put $\Xi=\{\pm\epsilon_1\pm\epsilon_2,\pm\epsilon_3\pm\epsilon_4\}
\simeq 4A_1$ and 
$\Xi'=\{\pm\epsilon_1\pm\epsilon_3,\pm\epsilon_2\pm\epsilon_4\}
\simeq 4A_1$. 
Then the subsystems $\Xi$ and $\Xi'$ of $D_4$ under the notation 
in \S\ref{sec:Dynkin} do not satisfy \eqref{eq:W2}.
When $\Sigma=B_n$, $C_n$, $F_4$ or $G_2$ and 
$\Xi=\Sigma^L$ and $\Xi'=\Sigma\setminus\Xi$,
\eqref{eq:W2} is not valid.
\end{rem}
%%%%%%%%%%%%%%%%%%%%%%%%%%%%%%%%%%
\section{Making tables}
We are ready to answer the questions in the introduction
by completing the tables in \S\ref{sec:table}.
In this section we do it when the root system $\Sigma$ is of 
the exceptional type. 
Following the argument in \S\ref{sec:thm}, 
we easily get all $\Xi$ satisfying $\Hom(\Xi,\Sigma)\ne\emptyset$ together
with $\#\overline\Hom(\Xi,\Sigma)$ and $\Xi^\perp$ by Theorem~\ref{thm}.
In fact, we start from the irreducible $\Xi$ and then examine other $\Xi$ 
by using \eqref{eq:XiDiv} in a suitable lexicographic order (as in the tables)
to avoid confusion (cf.~Example~\ref{exmp:A} xii)).

As a result we finally get $(\Xi^\perp)^\perp$ and the dual pairs.
Moreover \eqref{eq:dual1} tells us whether the dual pair is special or 
not.  
We will calculate $\#\{\Theta\subset\Phi\,;\,\langle\Theta\rangle\simeq\Xi\}$ 
in \S\ref{sec:fund}.

Now we prepare the lemma to examine the action of $W_\Sigma$ on the 
imbeddings of a root system $\Xi$ into $\Sigma$.
\begin{lem}
Let $\Xi_1$ and $\Xi_2$ be subsystems of $\Sigma$
with $\Xi_2\subset\Xi_1^\perp$.
Then
\begin{gather}
 \begin{split}
  &\#\bigl(\overline\Hom(\Xi_1,\Sigma)/\!\Out(\Xi_1)\bigr)
   =\#\bigl(\Hom(\Xi_2,\Xi_1^\perp)/\!\Out(\Xi_2)\bigr)
   =1\\
 &\quad\Rightarrow\quad
  \#\bigl(\overline\Hom(\Xi_1+\Xi_2,\Sigma)/\!\Out(\Xi_1+\Xi_2)\bigr)=1,
 \end{split}\label{eq:N1}
 \allowdisplaybreaks\\
 \begin{split}
  &\#\bigl(\overline\Hom(\Xi_1,\Sigma)/\!\Out'(\Xi_1)\bigr)
  =\#\bigl(\overline\Hom(\Xi_2,\Xi^\perp)/\!\Out'(\Xi_2)\bigr)
  =1\\
  &\quad\Rightarrow\quad
  \#\bigl(\overline\Hom(\Xi_1+\Xi_2,\Sigma)/\!\Out'(\Xi_1+\Xi_2)\bigr)=1,
 \end{split}\label{eq:N2}
  \allowdisplaybreaks\\
 \begin{split}
 &\begin{cases}
   \#\bigl(\overline\Hom(\Xi_1,\Sigma)/\!\Out(\Xi_1)\bigr)=
   \#\bigl(\Out(\Xi_1^\perp)\backslash\overline\Hom(\Xi_2,\Xi_1^\perp)/\!
   \Out'(\Xi_2)\bigr)=1,\\
%   \phantom{\#\bigl(\overline\Hom(\Xi_1,\Sigma)/\!\Out(\Xi_1)\bigr)}=1,\\
   \Out'(\Xi_1)\simeq\Out(\Xi_1)
   \text{ and \ }(\Xi_1,\Xi_1^\perp) \text{ is a special dual pair}
  \end{cases}\\
 &\qquad \Rightarrow\ 
    \#\bigl(\overline\Hom(\Xi_1+\Xi_2,\Sigma)/\!\Out'(\Xi_1+\Xi_2)\bigr)=1,
    \label{eq:ND}
 \end{split} \allowdisplaybreaks\\
 \begin{split}
  &\#\bigl(\overline\Hom(\Xi_1,\Sigma)/\!\Out(\Xi_1)\bigr)=
  \#\bigl(\overline\Hom(\Xi_2,\Xi_1^\perp)/\!\Out(\Xi_2)\bigr)=1\\
  &\text{ and \ }\iota(\Xi_2)^\perp\simeq \Xi_1
    \ (\forall\iota\in\Hom(\Xi_2,\Sigma))\ \Rightarrow\ 
  \#\bigl(\overline\Hom(\Xi_2,\Sigma)/\!\Out(\Xi_2)\bigr)=1.
  \label{eq:N3}
 \end{split}
\end{gather}
\end{lem}
\begin{proof}The claims
\eqref{eq:N1} and \eqref{eq:N2} %and \eqref{eq:ND}
are clear because for $\iota\in\Hom(\Xi_1+\Xi_2,\Sigma)$ 
the assumptions assure that 
there exists $w\in W_\Sigma$ such that $\iota(\Xi_1)=w(\Xi_1)$
and hence we may assume $\iota(\Xi_1)=\Xi_1$ in 
$\overline\Hom(\Xi_1+\Xi_2,\Sigma)$.  
Under the assumption in \eqref{eq:ND} 
there exists $w\in W_\Sigma$ such that $w\circ\iota$ stabilizes 
every irreducible component of $\Xi_2$ and therefore it also stabilizes
$\Xi_1$ and we have \eqref{eq:ND}.

The claim \eqref{eq:N3} is also clear because 
for $\iota\in\Hom(\Xi_2,\Sigma)$, $\exists w\in W_\Sigma$
such that $w\circ\iota(\Xi_2)^\perp=\Xi_1$, which implies
$w\circ\iota(\Xi_2)\subset \Xi_1^\perp$.
\end{proof}

\subsection{Type $E_6$}\label{sec:TE6}
The automorphism group of $G(\tilde\Psi)$ is of order 6,
which is generated by a rotation and a reflection.
Since the rotation has order 3, it corresponds to an element
of $W_{E_6}$ and the reflection corresponds to a non-trivial
element of $\Out(E_6)$.
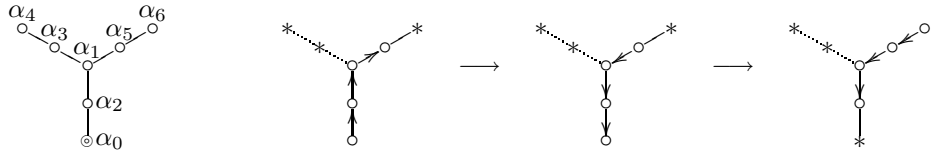
\[
\begin{xy}
 \ar@{-}  *+!D{\alpha_1}  *{\circ}="O";
  (-4.33,2.5)  *+!D{\alpha_3}  *{\circ}="C"
 \ar@{-} "C";(-8.66,5) *+!D{\alpha_4} *{\circ}="D"
 \ar@{-} "O";(4.33,2.5) *+!D{\alpha_5}  *{\circ}="E"
 \ar@{-} "E";(8.66,5) *+!D{\alpha_6}  *{\circ},
 \ar@{-} "O";(0,-5) *+!L{\alpha_2}  *{\circ}="X"
 \ar@{-} "X";(0,-10) *+!L{\alpha_0} *{\ccirc};
\end{xy}
\qquad\qquad
\begin{xy}
 \ar@{.}  *{\circ}="O";
  (-4.33,2.5)   *{\ast}="C"
 \ar@{.} "C";(-8.66,5)  *{\ast}="D"
 \ar@{->} "O";(4.33,2.5)  *{\circ}="E"
 \ar@{-} "E";(8.66,5)  *{\ast},
 \ar@{<-} "O";(0,-5)   *{\circ}="X"
 \ar@{->} "X";(0,-10)  *{\circ};
\end{xy}
\quad\longrightarrow\quad
\begin{xy}
 \ar@{.}  *{\circ}="O";
  (-4.33,2.5)   *{\ast}="C"
 \ar@{.} "C";(-8.66,5)  *{\ast}="D"
 \ar@{<-} "O";(4.33,2.5)  *{\circ}="E"
 \ar@{-} "E";(8.66,5)  *{\ast},
 \ar@{->} "O";(0,-5)   *{\circ}="X"
 \ar@{<-} "X";(0,-10)  *{\circ};
\end{xy}
\quad\longrightarrow\quad
\begin{xy}
 \ar@{.}  *{\circ}="O";
  (-4.33,2.5)   *{\ast}="C"
 \ar@{.} "C";(-8.66,5)  *{\ast}="D"
 \ar@{<-} "O";(4.33,2.5)  *{\circ}="E"
 \ar@{<-} "E";(8.66,5)  *{\circ},
 \ar@{->} "O";(0,-5)   *{\circ}="X"
 \ar@{-} "X";(0,-10)  *{\ast};
\end{xy}
\]
The set $\overline\Hom(A_4,E_6)$ has two elements which are shown in
Example~\ref{exmp:A}~v).  It also shows that
$\Out(E_6)$ non-trivially acts on this set.
If $A_4$ is imbedded to $E_6$ given as in the above 
imbedding $G(A_4)\subset G(\tilde E_6)$ with the starting
vertex $\{\alpha_0\}$,
the non-trivial action by $\Out(A_4)$ changes the starting vertex
as is shown above.
Then by an element of $W_{A_5}$ with
$A_5=\langle\alpha_0,\alpha_2,\alpha_1,\alpha_5,\alpha_6\rangle$
the imbedding is transformed as is shown by the second arrow.
Then the result corresponds to a reflection, which implies that
$\Out(A_4)$ also acts non-trivially on $\overline\Hom(A_4,E_6)$
and hence $\#\bigl(\overline\Hom(A_4,E_6)/\!\Out(A_4)\bigr)=1$.

The same argument works for $\Xi=A_5$, $A_2+A_1$ and $A_3+A_1$.
Similarly $(\alpha_2,\alpha_0,\alpha_5,\alpha_6)$ is transformed
to $(\alpha_6,\alpha_5,\alpha_0,\alpha_2)$ by an element of $W_{A_5}$
and furthermore to $(\alpha_0,\alpha_2,\alpha_4,\alpha_3)$ by a rotation.
Hence a non-trivial element of $\Out(A_2)$ for 
$A_2=\langle\alpha_0,\alpha_2\rangle$ induces the transposition 
of two irreducible components of $A_2^\perp\simeq A_2+A_2$, which implies
$\#\bigl(\overline\Hom(2A_2,E_6)/\!\Out'(2A_2)\bigr)=
%\#\bigl(\overline\Hom(3A_2,E_6)/\!\Out'(3A_2)\bigr)=
1$.

From our construction of the representatives of 
$\overline\Hom(\Xi,E_6)$ it is obvious to have 
$\#\bigl(\overline\Hom(\Xi,E_6)/\!\Out'(\Xi)\bigr)=1$ for
$\Xi=D_5$ (cf.~\eqref{eq:DinE}) and $E_6$ and we can easily calculate
$\#\bigl(\Out(E_6)\backslash\overline\Hom(\Xi,E_6)\bigr)$.
Put $\Sigma=E_6$ and let $(\Xi_1,\Xi_2)$ be any one of the 
pairs $(A_2+A_1,A_1)$, $(2A_2,A_1)$, $(2A_2,A_2)$, 
$(2A_1,A_3)$, $(A_4,A_1)$ and $(A_5,A_1)$.  Then applying \eqref{eq:N2}
to $\Sigma$ and $(\Xi_1,\Xi_2)$, we have 
$\#\bigl(\overline\Hom(\Xi_1+\Xi_2,\Sigma)/
\!\Out'(\Xi_1+\Xi_2)\bigr)=1$.

\subsection{Type $E_7$}\label{sec:TE7}
Note that $G(\tilde E_7)$ has an automorphism of order 2 and it corresponds
to an element of $W_{E_7}$ because $W_{E_7}=\Aut(E_7)$.

Let $\Sigma=E_7$ and let $(\Xi_1,\Xi_2)$ be any one of
$(A_1,D_6)$, $(A_2,A_2)$, $(A_2,A_2+A_1)$, $(A_2,A_3)$, $(A_2,A_3+A_1)$,
$(A_3,A_3)$ and $(A_3,A_3+A_1)$.
We have 
$\#\bigl(\overline\Hom(\Xi_1+\Xi_2, \Sigma)/\!\Out'(\Xi_1+\Xi_2)\bigr)=1$ 
by \eqref{eq:N2}.
Here we note that $A_1^\perp\simeq D_6$, 
$A_2^\perp\simeq A_5$ and $A_3^\perp\simeq A_3+A_1$.
We can apply \eqref{eq:ND} to $(\Xi_1,\Xi_2)=(D_4,kA_1)$ with $1\le k\le 3$
and we have the same conclusion.
Applying \eqref{eq:N1} to $(A_3,3A_1)$, 
we have $\#\bigl(\overline\Hom(A_3+3A_1,E_7)/\!\Out(A_3+3A_1)\bigr)=1$.

The subsystems $\Xi$ of $E_7$ which are isomorphic to $3A_1$
and satisfy $\Xi^\perp\simeq 4A_1$ are mutually equivalent by $\Sigma$.
Hence $\Xi^\perp \simeq 4A_1$ also have this property. 
Namely
\[
 \#\bigl(W_{E_7}\backslash\{\iota\in\Hom(4A_1,E_7)\,;\,
  (\iota(4A_1)^\perp)^\perp = \iota(4A_1)\}/\Aut(4A_1)\bigr)=1.
\]

Put $(A_1)_o=\langle\alpha_0\rangle$.
We have $G\bigl(\widetilde{(A_1)_o^\perp}\bigr)$ 
as is given in the following first diagram.
Put $(2A_1)_o=\langle\alpha_0,\alpha_p\rangle$.
Then the extended Dynkin 
diagrams of the components of $(2A_1)_o^\perp=
\langle\alpha_2,\alpha_3,\alpha_4,\alpha_5,\alpha_7\rangle\simeq D_4+A_1$
are also given by the following second diagram.
These diagrams correspond to the last figure in Remark~\ref{rem:A1}~i).
Here %we put
\begin{gather*}
 \begin{split}
 -\alpha_p
 &:=
  (\alpha_2+\alpha_3+\cdots+\alpha_7)
 +(\alpha_4+\alpha_5+\alpha_6)
 = \epsilon_5+\epsilon_6,\\
 -\alpha_q
 &:=
  \alpha_2+\alpha_3+2\alpha_4+\alpha_5
  =\epsilon_3+\epsilon_4,
 \end{split}
 \allowdisplaybreaks\\[-6pt]
\begin{xy}
 \ar@{.}
   *+!D{\alpha_0} *{\bullet};
  (5,0)  *+!D{\alpha_1} *{\ast}="A",
 \ar@{.} 
 "A";(10,0) *+!D{\alpha_3} *{\circ}="C"
 \ar@{-} "C";(15,0) *+!D{\alpha_4} *{\circ}="D"
 \ar@{-} "D";(20,0) *+!D{\alpha_5} *{\circ}="E"
 \ar@{-} "E";(25,0) *+!D{\alpha_6} *{\circ}="F"
 \ar@{-} "F";(30,0) *+!D{\alpha_7} *{\circ}
 \ar@{-} "D";(15,-5) *+!L{\alpha_2} *{\circ}
 \ar@{-} "F";(25,-5) *+!L{\alpha_p} *{\ccirc}
\end{xy}
\rightarrow
\begin{xy}
 \ar@{.}
  *+!D{\alpha_0} *{\bullet};
  (5,0)  *+!D{\alpha_1} *{\ast}="A",
 \ar@{.} "A";(10,0) *+!D{\alpha_3} *{\circ}="C"
 \ar@{-} "C";(15,0)  *{\circ}="D"
 \ar@{-} "D";(20,0) *+!D{\alpha_5} *{\circ}="E"
 \ar@{.}  "E";(25,0) *+!D{\alpha_6} *{\ast}="F"
 \ar@{.}  "F";(30,0) *+!D{\alpha_7} *{\ccirc}
 \ar@{-} "D";(15,-5) *+!L{\alpha_2} *{\circ}
 \ar@{.}  "F";(25,-5) *+!L{\alpha_p} *{\bullet}
 \ar@{-} "D";(15,5) *+!L{\alpha_q} *{\ccirc}
\end{xy}
\rightarrow
\raisebox{20pt}{
\begin{xy}
 \ar@{.}  *+!D{\alpha_0} *{\bullet};
  (5,0)  *+!D{\alpha_1} *{\ast}="A"
 \ar@{.} "A";(10,0) *+!D{\alpha_3} *{\ccirc}="C"
 \ar@{.} "C";(15,0)                *{\ast}="D"
 \ar@{.} "D";(20,0) *+!D{\alpha_5} *{\ccirc}="E"
 \ar@{.} "E";(25,0) *+!D{\alpha_6} *{\ast}="F"
 \ar@{.} "F";(30,0) *+!D{\alpha_7} *{\ccirc}
 \ar@{.} "D";(15,-5) *+!L{\alpha_2} *{\ccirc}
 \ar@{.} "F";(25,-5) *+!L{\alpha_p} *{\bullet}
 \ar@{.} "D";(15,5) *+!L{\alpha_q} *{\bullet}
 \ar@{.}  (0,-15) *+!D{\alpha_0} *{\bullet};
  (5,-15)  *+!D{\alpha_1} *{\ast}="X"
 \ar@{.}  "X";(10,-15) *+!D{\alpha_3} *{\circ}="Y"
 \ar@{-}  "Y";(15,-15)                *{\circ}="Z"
 \ar@{-} "Z";(20,-15) *+!D{\alpha_5} *{\circ}="P"
 \ar@{.} "P";(25,-15) *+!D{\alpha_6} *{\ast}="Q"
 \ar@{.} "Q";(30,-15) *+!D{\alpha_7} *{\bullet}
 \ar@{-} "Z";(15,-20) *+!L{\alpha_2} *{\circ}
 \ar@{.} "Q";(25,-20) *+!L{\alpha_p} *{\bullet}
 \ar@{-} "Z";(15,-10) *+!L{\alpha_q} *{\ccirc}
\end{xy}
}
\end{gather*}
In the above diagrams the vertices expressed by asterisks are considered
to be removed and the diagrams are (extended) Dynkin diagrams for other roots.

There are two equivalence classes in the imbeddings of $3A_1$ to $E_7$,
whose representatives are
\begin{align*}
 (3A_1)_1 = \langle\alpha_0,\alpha_p,\alpha_q\rangle,
 \qquad
 (3A_1)_2 = \langle\alpha_0,\alpha_p,\alpha_7\rangle,
\end{align*}
which satisfy
\begin{align*}
 (3A_1)_1^\perp = \langle\alpha_2,\alpha_3,\alpha_5,\alpha_7\rangle
 \simeq 4A_1,\qquad
 (3A_1)_2^\perp = \langle\alpha_2,\alpha_3,\alpha_4,\alpha_5\rangle
 \simeq D_4.
\end{align*}
Thus the image of the imbedding of $4A_1$ to $E_7$ is equivalent to
one of the following subsystems $(4A_1)_j$ of $E_7$:
\begin{align*}
 (4A_1)_1 &=\langle\alpha_0,\alpha_p,\alpha_q,\alpha_7\rangle,
&&
 (4A_1)_1^\perp=\langle\alpha_2,\alpha_3,\alpha_5\rangle
  %=\langle\epsilon_1\pm\epsilon_2,\epsilon_4-\epsilon_3\rangle
  ,
\allowdisplaybreaks\\
 (4A_1)_2 &=\langle\alpha_0,\alpha_p,\alpha_q,\alpha_5\rangle,
&&
 (4A_1)_2^\perp=\langle\alpha_2,\alpha_3,\alpha_7\rangle
  %=\langle\epsilon_1\pm\epsilon_2,\epsilon_6-\epsilon_5\rangle
  ,
\allowdisplaybreaks\\
 (4A_1)_3 &=\langle\alpha_0,\alpha_p,\alpha_q,\alpha_2\rangle,
&&
 (4A_1)_3^\perp=\langle\alpha_3,\alpha_5,\alpha_7\rangle
  %=\langle\epsilon_1-\epsilon_2,\epsilon_4-\epsilon_3,
  %\epsilon_6-\epsilon_5\rangle
  ,
\allowdisplaybreaks\\
 (4A_1)_4 &=\langle\alpha_0,\alpha_p,\alpha_q,\alpha_3\rangle,
&&
 (4A_1)_4^\perp=\langle\alpha_2,\alpha_5,\alpha_7\rangle
  %=\langle\epsilon_1+\epsilon_2,\epsilon_4-\epsilon_3,
  % \epsilon_6-\epsilon_5\rangle
  .
\end{align*}
In view of Remark~\ref{rem:thmauto} ii) 
the above procedures for $3A_1\subset E_6$
can be also explained by the following isomorphic ones.
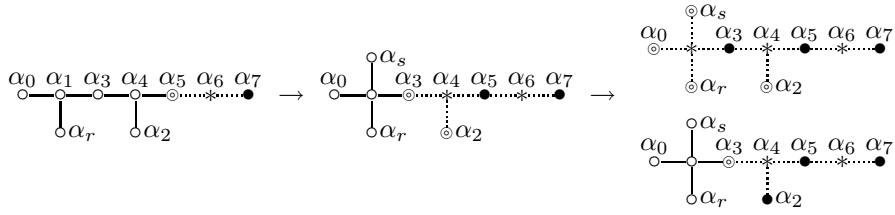
\[
\begin{xy}
 \ar@{-}  *+!D{\alpha_0} *{\circ};
  (5,0)  *+!D{\alpha_1} *{\circ}="A"
 \ar@{-} "A";(10,0) *+!D{\alpha_3} *{\circ}="C"
 \ar@{-} "C";(15,0) *+!D{\alpha_4} *{\circ}="D"
 \ar@{-} "D";(20,0) *+!D{\alpha_5} *{\ccirc}="E"
 \ar@{.} "E";(25,0) *+!D{\alpha_6} *{\ast}="F"
 \ar@{.} "F";(30,0) *+!D{\alpha_7} *{\bullet}
 \ar@{-} "D";(15,-5) *+!L{\alpha_2} *{\circ}
 \ar@{-} "A";(5,-5) *+!L{\alpha_r} *{\circ}
\end{xy}
\rightarrow
\begin{xy}
 \ar@{-}  *+!D{\alpha_0} *{\circ};
  (5,0)   *{\circ}="A"
 \ar@{-} "A";(10,0) *+!D{\alpha_3} *{\ccirc}="C"
 \ar@{.} "C";(15,0) *+!D{\alpha_4} *{\ast}="D"
 \ar@{.} "D";(20,0) *+!D{\alpha_5} *{\bullet}="E"
 \ar@{.} "E";(25,0) *+!D{\alpha_6} *{\ast}="F"
 \ar@{.} "F";(30,0) *+!D{\alpha_7} *{\bullet}
 \ar@{.} "D";(15,-5) *+!L{\alpha_2} *{\ccirc}
 \ar@{-} "A";(5,-5) *+!L{\alpha_r} *{\circ}
 \ar@{-} "A";(5,5) *+!L{\alpha_s} *{\circ}
\end{xy}
\rightarrow
\raisebox{20pt}{
\begin{xy}
 \ar@{.}  *+!D{\alpha_0} *{\ccirc};
  (5,0)     *{\ast}="A"
 \ar@{.} "A";(10,0) *+!D{\alpha_3} *{\bullet}="C"
 \ar@{.} "C";(15,0) *+!D{\alpha_4} *{\ast}="D"
 \ar@{.} "D";(20,0) *+!D{\alpha_5} *{\bullet}="E"
 \ar@{.} "E";(25,0) *+!D{\alpha_6} *{\ast}="F"
 \ar@{.} "F";(30,0) *+!D{\alpha_7} *{\bullet}
 \ar@{.} "A";(5,-5) *+!L{\alpha_r} *{\ccirc}
 \ar@{.} "D";(15,-5) *+!L{\alpha_2} *{\ccirc}
 \ar@{.} "A";(5,5) *+!L{\alpha_s} *{\ccirc}
 \ar@{-}  (0,-15) *+!D{\alpha_0} *{\circ};
  (5,-15)  *{\circ}="X"
 \ar@{-} "X";(10,-15) *+!D{\alpha_3} *{\ccirc}="Y"
 \ar@{.} "Y";(15,-15) *+!D{\alpha_4} *{\ast}="Z"
 \ar@{.} "Z";(20,-15) *+!D{\alpha_5} *{\bullet}="P"
 \ar@{.} "P";(25,-15) *+!D{\alpha_6} *{\ast}="Q"
 \ar@{.} "Q";(30,-15) *+!D{\alpha_7} *{\bullet}
 \ar@{-} "X";(5,-20) *+!L{\alpha_r} *{\circ}
 \ar@{.} "Z";(15,-20) *+!L{\alpha_2} *{\bullet}
 \ar@{-} "X";(5,-10) *+!L{\alpha_s} *{\circ}
\end{xy}
}
\]
Here $\alpha_r=-\alpha_p$ and $\alpha_s=-\alpha_q$.
In fact $\alpha_p$, 
$\alpha_r\in\langle\alpha_2,\alpha_3,\alpha_4,\alpha_5,\alpha_7\rangle^\perp
=(D_4+A_1)^\perp\simeq A_1$.
As $m_7(\alpha_p)<0$ and $m_7(\alpha_r)>0$, we have $\alpha_r=-\alpha_p$.
Similarly we have $\alpha_s=-\alpha_q$ from 
$\alpha_q$, $\alpha_s
\in\langle\alpha_0,\alpha_2,\alpha_3,\alpha_5,\alpha_7,\alpha_p\rangle^\perp
\simeq A_1$.
This is also easily verified by the Dynkin diagrams with the coefficients  
$m_j(-\alpha_0)$ in \S\ref{sec:Dynkin}.

Note that 
\[\langle\alpha_2,\alpha_3,\alpha_5\rangle
 \underset{\langle\alpha_5,\alpha_6,\alpha_7\rangle}{\sim}
 \langle\alpha_2,\alpha_3,\alpha_7\rangle
 \underset{\langle\alpha_1,\alpha_2,\alpha_3,\alpha_4\rangle}{\sim}
 \langle\alpha_1,\alpha_4,\alpha_7\rangle
 \underset{\langle\alpha_1,\alpha_3,\cdots,\alpha_7\rangle}{\sim}
 \langle\alpha_3,\alpha_5,\alpha_7\rangle.
\]

Thus we can conclude
\begin{equation}
  \bigl((4A_1)_j^\perp\bigl)^\perp
   \begin{cases}
   \simeq\langle\alpha_3,\alpha_5,\alpha_7\rangle^\perp
   = \langle\alpha_0,\alpha_2,\alpha_r,\alpha_s\rangle
     \simeq 4A_1 &(j=1,2,3),\\
   = \langle\alpha_0,\alpha_1,\alpha_3,\alpha_r,\alpha_s\rangle\simeq D_4
    &(j=4).
  \end{cases}
\end{equation}
Since $(4A_1)_j^\perp$ for $j=1,2,3$ are equivalent to each other by $E_8$,
so are the subsystems
$(4A_1)_j=((4A_1)_j^\perp)^\perp$ for $j=1,2,3$.
Moreover we have
\begin{equation}\label{eq:4A14}
 \langle\alpha_0,\alpha_p,\alpha_q,\alpha_3\rangle
 \underset{\langle\alpha_2,\alpha_3,\ldots,\alpha_7\rangle}{\sim}
 \langle\alpha_0,\alpha_2,\alpha_q,\alpha_7\rangle
 \underset{\langle\alpha_2,\alpha_3,\alpha_4,\alpha_5\rangle}{\sim}
 \langle\alpha_0,\alpha_3,\alpha_5,\alpha_7\rangle.
\end{equation}

Put $P_\Xi=\{\Theta\subset\Psi\,;\,\langle\Theta\rangle\simeq\Xi\}$.
It is easy to see that if $\Theta\in P_{3A_1}$ satisfies 
$\Theta\cap\{\alpha_1,\alpha_3\}\ne\emptyset$,
$\langle\Theta\rangle\underset{E_7}\sim (3A_1)_1$.
Moreover if $\Theta\in P_{3A_1}$ satisfies 
$\Theta\cap\{\alpha_1,\alpha_3\}=\emptyset$,
then $\Theta=\{\alpha_2,\alpha_5,\alpha_7\}$.
We will have $\# B_{3A_1}=11$ in \S\ref{sec:fund}.

Applying \eqref{eq:N3} to $(\Xi_1,\Xi_2)=(2A_1,5A_1)$ and
$(A_1,6A_1)$ with $\Sigma=E_7$, we have
$\#\bigl(\overline\Hom(\Xi_2,\Sigma)/\!\Out(\Xi_2)\bigr)=1$, respectively.
Similarly applying \eqref{eq:N1} to $(\Xi_1,\Xi_2)=
(5A_1,A_1)$ and $(5A_1,2A_1)$, we have
$\#\bigl(\overline\Hom(\Xi_1+\Xi_2,\Sigma)/\!\Out(\Xi_1+\Xi_2)\bigr)=1$, 
respectively.

\subsection{Type $E_8$}\label{sec:TE8}
Applying \eqref{eq:N3} to $(3A_1,5A_1)$ and
then \eqref{eq:N1} to $(5A_1,A_1)$, $(5A_1,2A_1)$ and $(5A_1,3A_1)$,
respectively, 
we have $\#\bigl(\overline\Hom(kA_1,E_8)/\!\Out(kA_1)\bigr)=1$ 
for $k=5,6,7$ and $8$.
See \S\ref{sec:8A1E8} to get further results on $\Hom(kA_1,E_8)$
with $1\le k\le8$. 

If $(\Xi_1,\Xi_1^\perp)$ is any one of the pairs
$(A_2,E_6)$, $(A_4,A_4)$, $(D_4,D_4)$, $(D_5,A_3)$ and $(D_6,A_1)$,
we have
\[
 \Hom(\Xi_2,\Xi_1^\perp)\ne\emptyset
 \ \Rightarrow\ 
 \#\bigl(\overline\Hom(\Xi_1+\Xi_2,E_8)/\!\Out'(\Xi_1+\Xi_2)\bigr)=1
\]
by applying \eqref{eq:ND}.
Hence if $\Xi$ contains $A_2$, $A_4$, $D_4$, $D_5$ or $D_6$ as an
irreducible component, the value of the column indicated by 
$\#_{\Xi'}$ equals one.
Moreover \eqref{eq:N1} can be applied to
$(\Xi_1,\Xi_2)=(A_3,3A_1)$, $(A_3,4A_1)$ and $(A_5,2A_1)$.
The number $\#\bigl(\overline\Hom(\Xi,E_8)/\!\Out'(\Xi)\bigr)$
for $\Xi=A_3+3A_1$, $A_3+4A_1$ and $A_5+2A_1$ is easily obtained
from $(A_3+A_1)^\perp\simeq A_3+A_1$ and $A_5^\perp\simeq A_2+A_1$.

Put
\begin{align*}
  (A_7)_o&=\langle \alpha_1,\alpha_3,\ldots,\alpha_8\rangle
 \subset (A_8)_o=\langle (A_7)_o,\alpha_0\rangle,\\
  (D_6)_o&=\{\pm\epsilon_i\pm\epsilon_j\,;\,3\le i<j\le 8\}
      =\{\alpha_2,\alpha_3\}^\perp,\\
  P_\Xi&=\{\Theta\subset\{\alpha_1,\ldots,\alpha_8\}\,;\,\langle\Theta
         \rangle\simeq \Xi\}.
\end{align*}
Then we note the following for $\Theta_1$, $\Theta_2\in P_\Xi$. 

$\Theta_1\underset{\Sigma}{\sim} \Theta_2$ if
$\Theta_j\subset(A_8)_o$ for $j=1$, $2$ and $\Theta_1\simeq\Theta_2$.

If $\Theta_1\ni\{\alpha_2,\alpha_3\}$, then
$\Theta_1^\perp = \Theta_1^\perp \cap (D_6)_o$.

Using these facts, we can easily examine $P_\Xi$. 
For example, any $\Theta\in P_{4A_1}$ satisfies 
$\langle\Theta\rangle\underset{E_8}{\sim}(4A_1)_o:
=\langle\alpha_2,\alpha_3,\alpha_6,\alpha_8\rangle$.
Here
$(4A_1)_o^\perp=
\langle\alpha_6,\alpha_8\rangle^\perp\cap D_6\simeq 4A_1$.
%%%
\subsection{Type $F_4$ and $G_2$}
It is easy to examine the cases when $\Sigma=F_4$ and $G_2$
by using Theorem~\ref{thm} together with
Remark~\ref{rem:red}, \eqref{eq:SL} and \eqref{eq:DF4}.
%%%
\subsection{Fundamental subsystems}\label{sec:fund}
We will give the number of the elements
$P_\Xi:=\{\Theta\subset\Psi\,;\,\langle\Theta\rangle\simeq\Xi\}$
for a subsystem $\Xi$ of $\Sigma$ when $\Sigma$ is of 
the exceptional type.
If $\#\Theta=\#\Psi-1$, it is easy to specify $\Xi$ that is isomorphic to
$\langle\Theta\rangle$ and we get the corresponding $\#P_\Xi$.
Other $\langle\Theta\rangle$ are fundamental subsystems of these maximal ones
and hence it is also easy to know whether 
$P_\Xi=\emptyset$ or not.
Note that $\rank\langle\Theta\rangle=\#\Theta$.

The number $\#P_\Xi$ can be inductively calculated as follows.
Let denote the number by $[\Xi,\Sigma]$.
For simplicity $\sum_j m_jA_j$ may be denoted 
$1^{m_1}\cdot 2^{m_2}\cdots$ with omitting the terms satisfying
$m_j=0$.

If $\Sigma$ is of type $E_n$, we divide $P_\Xi$ into the
subsets according to the relation with 
%the irreducible component of $\Xi$ %which contains 
the end root $\alpha_n$.
For example, suppose $\Sigma=E_6$ and $\Theta\subset\Psi$
satisfies $\langle\Theta\rangle\simeq 2A_1$.
Then if $\alpha_6\in\Theta$, the other element of $\Theta$ is in 
$\alpha_6^\perp\simeq A_4$.
If $\alpha_6\not\in\Theta$, $\Theta$ is contained in 
$\Psi\setminus\{\alpha_6\}\simeq D_5$.
Thus we have $[1^2,E_6]=[1,A_4]+[1^2,D_5]$.
Now it is quite easy to have $[1,A_4]=4$ and $[1^2,D_5]=6$.
Note that $[1^2,D_5]=[1,A_2+A_1]+[1^2,A_4]=3+3=6$.
We will show such calculations except for quite easy cases.
\begin{align*}
 [1^2,E_6] &= [1,A_4] +  [1^2,D_5] = 4 + 6 = 10,\\
 [1^2,E_7] &= [1,D_5] +  [1^2,E_6] = 5+10 = 15,\\
 [1^2,E_8] &= [1,E_6] +  [1^2,E_7] = 6+15 = 21,
 \allowdisplaybreaks\\
 [1^3,E_6] &= [1^2,A_4] + [1^3,D_5] = 3 + 2 = 5,\\
 [1^3,E_7] &= [1^2,D_5] + [1^3,E_6] = 6+5 = 11,\\
 [1^3,E_8] &= [1^2,E_6] + [1^3,E_7] = 10+11 = 21, 
 \allowdisplaybreaks\\
 [1^4,E_7] &= [1^3,D_5] + [1^4,E_6] = 2+0=2,\\
 [1^4,E_8] &= [1^3,E_6] +[1^4,E_7] = 5+2 = 7,
 \allowdisplaybreaks\\
 [2\cdot 1,E_6] &= [1,A_2+A_1] + [2,A_4] +[2\cdot 1,D_5]
           = 3 + 3 + 4 = 10,\\
 [2\cdot1,E_7]  &= [1,A_4] + [2,D_5] + [2\cdot1,E_6]
            = 4+4+10 = 18,\\
 [2\cdot1,E_8] &=
  [1,D_5]+[2,E_6]+[2\cdot1,E_7]=5+5+18=28,
 \allowdisplaybreaks\\
 [2\cdot1^2,E_7] &= [1^2,A_4] + [2\cdot 1,D_5] + [2\cdot 1^2,E_6]
            = 3+4+5 = 12,\\
 [2\cdot 1^2,E_8] &= [1^2,D_5]+[2\cdot 1,E_6]+[2\cdot 1^2,E_7]
   =6+10+12=28, 
\allowdisplaybreaks\\
 [2\cdot 1^3,E_8] &=
  [1^3,D_5]+[2\cdot 1^2,E_6]+[2\cdot 1^3,E_7]
   =2+5+0=7, \allowdisplaybreaks\\
 [2^2,E_7] &= [2,A_4] + [2^2,E_6] = 3+1 = 4,\\
 [2^2,E_8] &=
  [2,D_5]+[2^2,E_7] = 4+4 = 8, 
\allowdisplaybreaks\\
  [2^2\cdot1,E_8] &=
  [2\cdot1,D_5]+[2^2,E_6]+[2^2\cdot1,E_7]
   = 4+1+4=9, \allowdisplaybreaks\\
  [2^2\cdot1^2,E_8] &=
  [2\cdot 1^2,D_5]+[2^2\cdot1,E_6]+[2^2\cdot 1,E_7]
   = 1+1+0=2, 
\allowdisplaybreaks\\
 [3\cdot 1,E_7] &= [1,A_2+A_1] + [3, D_5] + [3\cdot1,E_6]
            = 3+4+4 = 11,\\
[3\cdot1,E_8] &=
  [1,A_4]+[3,E_6]+[3\cdot1,E_7] = 4+5+11 = 20,\\
[3\cdot 1^2,E_8] &=
  [1^2,A_4]+[3\cdot1,E_6]+[3\cdot 1^2,E_7]
  = 3+4+3=10, 
\allowdisplaybreaks\\
[3\cdot2,E_8] &=
  [2,A_4] + [3,D_5] + [3\cdot2,E_7] = 3+4+3 = 10,
\allowdisplaybreaks\\
 [3\cdot2\cdot1,E_8] &=
  [2\cdot1,A_4] + [3\cdot1,D_5] + [3\cdot 2,E_6]
  + [3\cdot2\cdot1,E_7]\\
 &= 2+1+0+1 = 4,
\allowdisplaybreaks\\
[3^2,E_8] &= [3,A_4]+[3^2,E_7] = 2+0 = 2, 
\allowdisplaybreaks\\
[4\cdot1,E_7] &= [1,A_1] + [4,D_5] + [4\cdot1,E_6]
            = 1 + 2 + 2 = 5,\\
[4\cdot1,E_8] &=
  [1,A_2+A_1]+[4,E_6]+[4\cdot1,E_7] = 3+4+5 = 12,\\
[4\cdot2,E_8] &=
  [2,A_2+A_1]+[4,D_5]+[4\cdot2,E_7] = 1+2+1 = 4.
\end{align*}
\section{Some remarks}\label{sec:rem}
\subsection{Some results from the tables}
In this section we always assume that $\Xi$ is a subsystem of an irreducible
root system $\Sigma$.

By our classification we have the following remarks.
\begin{rem}\label{rem:T1}
{\rm i)} The numbers of equivalence classes of certain subsystems $\Xi$ 
(cf.\ Remark~\ref{rem:maximal}) and 
their pairs are as follows.
Here we don't count the empty subsystem.

\smallskip
\noindent
\begin{tabular}{|l|c|c|c|c|c|}\hline
 \hfill$\Sigma\quad$ & $E_6$ & $E_7$ & $E_8$ & $F_4$ & $G_2$\\ \hline
%\hline
equivalence (isomorphic) classes 
  & 20 (20) & 46 (40) & 76 (71) & 36 (22) & 6 (4) \\ \hline
\!$S$-closed subsystems\!
  & 20  & 46  & 76  & 23 & 5  \\ \hline
\!$L$-closed ($\perp$-closed) subsystems
 & 16 (7)\phantom{2} & 31 (13) & 40 (18) & 11 (9)\phantom{2} & 3 (3) \\ \hline
$\Xi^\perp = \emptyset\ (\rank\Xi=\rank\Sigma)$
 & 10 (3)\phantom{2} & 19 (7)\phantom{2} & 33 (13) & 20 (16) & 4 (4) \\ \hline
maximal ($S$-closed) subsystems & 3 (3) & 4 (4) & 5 (5) & 3 (3) & 3 (2)
\\ \hline
%fundamental subsystems
%  & 16 & 31 & 40 & 10 & 3 \\ \hline
%isomorphic classes  & 20 & 40 & 71 & 22 & 4\\ \hline
dual pairs (special dual pairs) & 3 (1) & 6 (3) & 11 (11) & 5 (4) &  1 (1) \\ \hline
\end{tabular}
%The number of equivalence classes of non-empty subsystems of $\Sigma$
%is 20, 46, 76, 36 or 6 according to $\Sigma$ is of type $E_6$, $E_7$, $E_8$,
%$F_4$ or $G_2$, respectively.

%Let $\Xi$ be a subsystem of an irreducible root system $\Sigma$.

\smallskip
{\rm ii)} Let $\sigma$ be an outer automorphism of $\Sigma$.
Then $\sigma(\Xi)\underset{\Sigma}\sim\Xi$ if $(\Sigma,\Xi)$ 
does not satisfy the following condition.
\begin{equation}
  \Sigma\simeq D_n \text{ with an even }n,\ 
  \Xi\underset{\Sigma}{\sim}\textstyle\sum_j m_jA_j 
  \text{ and }\sum_j (j+1)m_j=n.
\end{equation}

{\rm iii)}
Suppose $\Xi$ is irreducible.
Then $\Xi^\perp\cap\Sigma$ is also irreducible 
if $(\Sigma,\Xi)$ is not isomorphic to any one in
the following list:

\smallskip
\noindent
\begin{tabular}{|l|l|c||l|l|c|}\hline
 $\Sigma$ & $\Xi$ & $\Xi^\perp$ & $\Sigma$ & $\Xi$ & $\Xi^\perp$
    \\ \hline
 $B_{n\ (n\ge3)}$ & $A_1$ & $B_{n-2}+A_1\text{ or }B_{n-1}$ &
 $C_{n\ (n\ge3)}$ & $A_1$ & $C_{n-2}+A_1\text{ or }C_{n-1}$
   \\ \hline
 $D_{n\ (n\ge4)}$ & $A_1$ & $D_{n-2}+A_1$ & 
  & &
    \\ \hline
 $D_5$ & $A_3$ & $2A_1$ \text{ or } $\emptyset$ &
 $D_6$ & $A_3$ & $A_3$ \text{ or } $2A_1$
   \\ \hline
 $D_{n\ (n\ge 7)}$ & $A_{n-3}$ & $2A_1$ &
 $D_{n\ (n\ge6)}$ & $D_{n-2}$ & $2A_1$
   \\ \hline
 $E_6$ & $A_2$ & $2A_2$ &
 $E_6$ & $A_3$ & $2A_1$
  \\ \hline
 $E_7$ & $A_3$ & $A_3+A_1$ &
 $E_7$ & $D_4$ & $3A_1$
  \\ \hline
 $E_8$ & $A_5$ & $A_2+A_1$ &
 $E_8$ & $D_6$ & $2A_2$
  \\ \hline
\end{tabular}

\smallskip
{\rm iv) }
The $L$-closure $\tilde\Xi$ of $\Xi$ in $\Sigma$ 
(cf.~Definition \ref{defn:S-closed} and Remark \ref{rem:closed}) can be 
easily obtained from the table in \S\ref{sec:table}.
Note that $\tilde\Xi$ is the maximal subsystem satisfying
\begin{equation}
  \Xi\subset \tilde\Xi\subset (\Xi^\perp)^\perp
  \text{ and }\rank\Xi=\rank\tilde\Xi.
\end{equation}
\end{rem}
%%%%%%%%%%%%%%%%% Orthogonal Systems %%%%%%%%%%%%%%%%
\begin{rem}[orthogonal systems]\label{rem:orth}
A subsystem $\Xi$ of $\Sigma$ or the fundamental system of $\Xi$
is called an \textsl{orthogonal system} of $\Sigma$
if $\Xi$ is isomorphic to $mA_1$ with a certain positive integer $m$.
An orthogonal system $\Xi$ is called \textsl{maximal} if $\Xi^\perp=\emptyset$
and called \textsl{strongly orthogonal} if $\Xi$ is $S$-closed.

Suppose $\Sigma$ is irreducible.
Let $\Xi=\langle\alpha_1,\ldots,\alpha_m\rangle$
and $\Xi'=\langle\alpha'_1,\ldots,\alpha'_m\rangle$
be orthogonal systems of $\Sigma$ with rank $m$.

{\rm i)} The rank of a maximal orthogonal system
is given in Corollary~\ref{cor:A}~i) when $\Sigma$ is simply laced.
If $\Sigma$ is not simply laced, the rank equals $\rank \Sigma$.

{\rm ii)} If one of the following conditions is satisfied,
then $\Xi\underset{\Sigma}{\sim}\Xi'$.
\begin{align}
 &\text{$\Xi$ and $\Xi'$ are strongly orthogonal maximal systems},\\
 &\text{$\Sigma$ is of type $A_n$, $E_6$, $E_7$ or $E_8$ and 
   $(\Sigma,\Xi)$ is not isomorphic to}\\
 &\quad\text{$(E_7,3A_1)$, $(E_7,4A_1)$ or $(E_8,4A_1)$.}\notag
\end{align}

{\rm iii)}
Let $\iota$ be a bijective map of $\Xi$ to $\Xi'$
with $(\iota(\alpha_j)|\iota(\alpha_j))=(\alpha_j|\alpha_j)$ 
for $j=1,\ldots,m$.
Suppose $m\ge2$.
Then there exists $w\in W_\Sigma$ with $\iota=w|_{\Xi}$ if one of the
following conditions is satisfied.
\begin{align}
 &\text{$\Sigma$ is of type $A_n$, $B_2$, $E_6$, $F_4$ or $G_2$.}
  \allowdisplaybreaks\\
 &\text{$\Sigma$ is of type $E_7$ with $m\le2$.}
  \allowdisplaybreaks\\
 &\text{$\Sigma$ is of type $E_8$ with $m\le3$.}
\end{align}
See \S\ref{sec:orth} and \S\ref{sec:8A1E8} for more details.
%%%%%%%%%%%%%%  Fundamental Subsystems %%%%%%%%%%%%
\end{rem}
\begin{rem}[fundamental subsystem]\label{rem:fund} 
{\rm i)} %Let $\Xi$ be a subsystem of an irreducible root system $\Sigma$.
%Then
We have
\begin{align}
 &\text{$\Xi$ is $L$-closed}\ \Leftrightarrow\ 
  \text{$\Xi$ is fundamental},\label{eq:LF}\\
 &\text{$\Xi$ is $\perp$-closed}\ \Rightarrow\ 
 \text{$\Xi$ is fundamental}.%,\\
% &\text{$(\Xi_1,\Xi_2)$ is a dual pair in $\Sigma$}\ \Rightarrow\ 
% \text{$\Xi_1$ and $\Xi_2$ are fundamental}.
\end{align} 
The minimal fundamental subsystem containing $\Xi$ is the
$L$-closure of $\Xi$.

{\rm ii)}
For %a subsystem $\Xi$ of an irreducible
%root system $\Sigma$ and 
a subset $\Theta\subset\Psi$
\begin{equation}\label{eq:P}
 \langle\Theta\rangle\simeq \Xi\text{ \ and \ }
 \langle\Theta\rangle\cap\Sigma^L\simeq\Xi\cap\Sigma^L
 \quad \Rightarrow\quad
 \langle\Theta\rangle\underset{\Sigma}{\overset{w}{\sim}}\Xi
\end{equation}
if $(\Sigma,\Xi)$ is not isomorphic to any one of the following list.
\begin{align}
 \begin{split}
  &\Sigma\text{ is of type $B_n\ (n\ge2)$, $C_n\ (n\ge3)$ or $D_n\ (n\ge4)$}
  \\
  &\text{ and $\Xi$ has at least one $A_3$-component or two $A_1$-components.}
 \end{split}\label{eq:P2}
\allowdisplaybreaks\\
 \begin{split}
 &(E_7,4A_1),\ (E_7,A_3+2A_1),\ (E_7,A_5+A_1),\\
 &(E_8,4A_1),\ 
 (E_8,A_3+2A_1),\ (E_8,2A_3),\ (E_8,A_5+A_1),\ (E_8,A_7).\label{eq:P2E}
 \end{split}
\allowdisplaybreaks\\
 &
 (E_7,3A_1),\  (E_7,A_5),\  (E_7,A_3+A_1).\label{eq:P3}
\end{align}

If $(\Sigma,\Xi)$ is isomorphic to one of the pairs in \eqref{eq:P2E}
and $\Xi$ is a fundamental subsystem, then \eqref{eq:P} is valid.

If $(\Sigma,\Xi)$ is isomorphic to one of the pairs in \eqref{eq:P3},
there exist $\Theta_1$, $\Theta_2\subset\Psi$ such that
$\Xi\simeq\langle\Theta_1\rangle\simeq\langle\Theta_2\rangle$,
$\langle\Theta_1\rangle\underset{\Sigma}{\not\sim}\langle\Theta_2\rangle$
and $\Xi\underset{\Sigma}{\sim}\langle\Theta_1\rangle$
or $\langle\Theta_2\rangle$.
%\if0

Hence if $(\Xi,\Sigma)$ is not isomorphic to any one of the pairs in 
\eqref{eq:P2} and \eqref{eq:P2E}, %we have
\begin{equation}\label{eq:PAL}
 \langle\Theta\rangle\simeq \Xi\text{ \ and \ }
 \langle\Theta\rangle\cap\Sigma^L\simeq\Xi\cap\Sigma^L
 \quad \Rightarrow\quad
 \exists\Theta'\subset\Psi\text{ such that }
 \langle\Theta'\rangle\underset{\Sigma}{\sim}\Xi.
\end{equation}
Note that \eqref{eq:PAL} is still valid even if $\Sigma$ is of type
$D_n$ except for the case
\begin{equation}\label{eq:PDA}
\text{$\Sigma$ is of type $D_n$ $(n\ge4)$ and $m_1+2m_3\ge 4$},
\end{equation}
where $m_j$ is the number of $A_j$-components of $\Xi$.
In fact the subsystems
\begin{equation}
\begin{split}
 \langle \epsilon_1\pm \epsilon_2,\, \epsilon_3\pm \epsilon_4\rangle
 &\underset{D_n}{\sim}2D_2\simeq 4A_1,\\
 \langle \epsilon_1\pm \epsilon_2,\, \epsilon_3-\epsilon_4,\,
 \epsilon_4\pm \epsilon_5\rangle&
 \underset{D_n}{\sim}D_2+D_3\simeq 2A_1+A_3,\\
 \langle \epsilon_1- \epsilon_2,\, \epsilon_2\pm \epsilon_3,\,
 \epsilon_4- \epsilon_5,\, \epsilon_5\pm \epsilon_6
 \rangle&\underset{D_n}{\sim}2D_3\simeq 2A_3
\end{split}
\end{equation}
of $D_n$ and the subsystems
\begin{equation}
\begin{split}
 \langle \epsilon_1\pm \epsilon_2\rangle & 
 \underset{B_n}{\sim}D_2\simeq 2A_1,\\
 \langle \epsilon_1-\epsilon_2,\,\epsilon_2\pm \epsilon_3\rangle 
  & \underset{B_n}{\sim}D_3\simeq A_3
\end{split}
\end{equation}
of $B_n$ are not fundamental.

{\rm iii)}  Given a subset $\Phi$ of %an irreducible root system 
$\Sigma$, 
we examine the condition 
\begin{equation}\label{eq:subfund}
 \exists w\in W_\Sigma\text{ such that }w(\Phi)\subset\Psi,
\end{equation}
namely, the condition that 
$\Phi$ can be extended to a fundamental system of $\Sigma$.

A subset $\Phi$ of $\Sigma$ is called a \textsl{$\Pi$-system} by \cite{Dy} if
$\Phi$ satisfies the two conditions
\begin{align}
  &\alpha\in\Phi,\ \beta\in\Phi\Rightarrow\alpha-\beta\notin\Sigma,
  \label{eq:PI1}\\
  &\text{the elements of $\Phi$ are linearly independent.}
  \label{eq:PI2}
\end{align}
It is easy to see that \eqref{eq:subfund} implies that $\Phi$ is a
$\Pi$-system.

Suppose $\Phi$ is a $\Pi$-system of $\Sigma$. 
Put $\Xi=\langle\Phi\rangle$.
Then $\Phi$ is a fundamental system of $\Xi$, 
which is shown by \cite[Theorem~5.1]{Dy}.
Hence we note that the fundamental system of a subsystem of $\Sigma$ 
is a $\Pi$-system if and only if the subsystem is $S$-closed.
Therefore if $G(\Phi)$ is a subdiagram of the Dynkin diagram of $\Sigma$
and the condition
\begin{equation}\label{eq:PIF}
 \begin{cases}
  (\Sigma,\Xi) \text{ is isomorphic to one of \eqref{eq:PDA} and 
  \eqref{eq:P2E},}
  \\
 \text{ or }\Sigma \text{ is of type $B_n$ with \eqref{eq:P2}}
 \end{cases}
\end{equation}
is not satisfied, it follows from our table that \eqref{eq:subfund} is valid.
When $\Xi$ is irreducible, \eqref{eq:PIF} is equivalent to the condition that
$(\Sigma,\Xi)\simeq(E_8,A_7)$ or $(B_n,A_3)$, which coincides with the result
given by \cite{AL}.
%\fi
\end{rem}
\begin{rem}[maximal subsystems]\label{rem:maximal}
{\rm i) }
A proper subsystem $\Xi$ of $\Sigma$ is called \textsl{maximal} if there
is no subsystem $\Xi'$ satisfying $\Xi\subsetneqq\Xi'\subsetneqq\Sigma$.
We have the following list of the maximal subsystems of irreducible root 
systems.

\iftrue
\smallskip
\centerline{\begin{tabular}{|l|l|l|}\hline
$\Sigma$ & $\Xi$: maximal,\ $\rank\Xi=\rank\Sigma$& 
 $\rank\Xi=\rank\Sigma-1$\\ \hline
$A_n$& & $A_{m-1}+A_{n-m}=\Xi_{m\ (2\le 2m\le n+1)}$\! \\ \hline
$B_n$& $B_m+B_{n-m\ (2\le 2m \le n)},\ D_n^L$ & \\ \hline
$C_n$& $C_m+C_{n-m\ (2\le 2m\le n)},\ D_n^S$ & \\ \hline
$D_n$& $D_m+D_{n-m}=\Xi_{m\ (4\le 2m\le n)}$ & 
  $D_{n-1}=\Xi_1,\ A_{n-1}=\Xi_n$ \\ \hline
$E_6$& $A_1+A_5=\Xi_2,\ 3A_2=\Xi_3$ & $D_5=\Xi_1$ \\ \hline
$E_7$& $A_1+D_6=\Xi_1,\ A_7=\Xi_2,\ A_2+A_5=\Xi_3$\! & $E_6=\Xi_7$ \\ \hline
$E_8$& $D_8=\Xi_1,\ A_8=\Xi_2,\ 2A_4=\Xi_5, $ & \\
   & $A_2+E_6=\Xi_7,\ A_1+E_7=\Xi_8$ & \\ \hline
$F_4$& $C_4,\ A_2^L+A_2^S,\ B_4$ & \\ \hline
$G_2$& $A_1^L+A_1^S,\ A_2^L,\ A_2^S$ &\\ \hline
\end{tabular}}
\smallskip
\else
\begin{equation}\label{eq:maximal}
\begin{split}
A_n&\supset A_m+A_{n-m-1}\quad(2\le 2m <n),\\
B_n&\supset B_m+B_{n-m}\quad(2\le 2m \le n),\ D_n^L,\\
C_n&\supset C_m+C_{n-m}\quad(2\le 2m\le n),\ D_n^S,\\
D_n&\supset D_m+D_{n-m}\quad(4\le 2m\le n),\ D_{n-1},\ A_{n-1},\\
E_6&\supset A_1+A_5,\ 3A_2,\ D_5,\\
E_7&\supset A_7,\ A_1+D_6,\ A_2+A_5,\ E_6,\\
E_8&\supset D_8,\ A_8,\ A_1+E_7,\ A_2+E_6,\ 2A_4,\\
F_4&\supset B_4,\ C_4,\ A_2^L+A_2^S,\\
G_2&\supset A_2^L,\ A_2^S,\ A_1^L+A_1^S.
\end{split}\end{equation}
\fi

Note that $D_2\simeq 2A_1$ and $D_3\simeq A_3$ in the above.

{\rm ii) }
A proper subsystem $\Xi$ of $\Sigma$ is called a 
\textsl{maximal $S$-closed subsystem} if $\Xi$ is $S$-closed and 
if there is no $S$-closed subsystem
$\Xi'$ satisfying $\Xi\subsetneqq\Xi'\subsetneqq\Sigma$.
The list of the maximal $S$-closed subsystems of the 
irreducible root system $\Sigma$ is same as in i) %\eqref{eq:maximal}
if $\Sigma$ is simply laced. % = A_n$, $D_n$, $E_6$, $E_7$ or $E_8$.
In the other cases we have

\iftrue
\smallskip
\centerline{
\begin{tabular}{|l|l|l|}\hline
$\Sigma$ & $\Xi$: $S$-closed maximal,\ $\rank\Xi=\rank\Sigma$
& $\rank\Xi=\rank\Sigma-1$ \\ \hline
$B_n$ & $D_m^L+B_{n-m}=\Xi_{m\ (2\le m\le n)}$& 
$B_{n-1}=\Xi_1$ \\ \hline
$C_n$ & $C_m+C_{n-m}=\Xi_{m\ (2\le 2m\le n)}$ & 
$A_{n-1}^S=\Xi_n$ \\ \hline
$F_4$ & $A_1^L+C_3=\Xi_1,\ A_2^L+A_2^S=\Xi_2,\ B_4=\Xi_4$ & \\ \hline
$G_2$ & $A_1^L+A_1^S=\Xi_1,\ A_2^L=\Xi_2$ & \\ \hline
\end{tabular}}
\smallskip
\else
\begin{equation}
\begin{split}
B_n&\supset D_m+B_{n-m}\quad(2\le m\le n),\ D_n^L,\ B_{n-1},\\
C_n&\supset C_m+C_{n-m}\quad(2\le 2m\le n),\ A_{n-1}^S,\\
F_4&\supset B_4,\ A_2^L+A_2^S,\\
G_2&\supset A_2^L,\ A_1^L+A_1^S.
\end{split}
\end{equation}
\fi
They are studied by \cite{BS} and \cite{W} (cf.~\cite{Dy})
and the Dynkin diagram of $\Xi$ is
\begin{equation}\label{eq:maxS}
\begin{cases}
G(\Psi\setminus\{\alpha_j\})&(m_j(\alpha_{max})=1),\\
%\text{ or}\\
G\bigl(\tilde\Psi\setminus\{\alpha_j\})
&(m_j(\alpha_{max})\text{ is a prime number }\ge2)
\end{cases}
\end{equation}
with a suitable $\alpha_j\in\Psi$ satisfying 
$m_j(\alpha_{\max})=1$, $2$, $3$ or $5$.

{\rm iii)}
Let $\Xi_j$ be the maximal subsystem of $\Sigma$ defined by 
\eqref{eq:maxS} with $\alpha_j\in\Psi$.
Then Proposition~\ref{prop:subsystem} (cf.~\cite[Theorem~3.1]{W}) implies
\begin{equation}
 \Xi_j=\bigl\{\alpha=\textstyle\sum_{\alpha_\nu\in\Psi} 
 m_\nu(\alpha)\alpha_\nu\in\Sigma\,;\,
 m_j(\alpha)\equiv0\mod\max\{2,m_j(\alpha_{\max})\}\bigr\}.
\end{equation}
Note that $\iota(A_m)\cap\Xi_j\ne\emptyset$ for any 
$\iota\in\Hom(A_m,\Sigma)$ if $m\ge\max\{2,m_j(\alpha_{max})\}$. 

This claim follows from the following fact with putting 
$n_i=m_j\bigl(\iota(\epsilon_i-\epsilon_{i+1})\bigr)$ for the 
$i$-th root $\epsilon_i-\epsilon_{i+1}$ in the fundamental system of $A_m$.

For a positive integer $m\ge 2$ and a sequence of integers
$n_1,\dots,n_m$ we can choose $1\le k\le k'\le m$ such that
$n_k+n_{k+1}+\dots+ n_{k'}\equiv0\mod m$.

For example, 
when $\Sigma=E_8$, we have $\Xi_5\simeq 2A_4$, $m_5(\alpha_{max})=5$ and
{\def\mL#1#2{\!\displaystyle\mathop{\mbox{$#1$}}_{\mbox{$#2$}}\!}
\begin{equation}
 E_8\setminus \Xi_5\supset \Xi':=
 \bigl\langle 00\mL001000,\, 01\mL211100,\,11\mL101110,\,
 00\mL111111\bigr\rangle\simeq A_4
\end{equation}
under the notation in \S\ref{sec:table}.}
Here the roots $\alpha\in\Sigma$ are indicated by the numbers $m_\nu(\alpha)$ 
arranged according to the Dynkin diagram of $\Sigma$.
We have $W_{\Xi_5}\cap W_{\Xi'}=\{e\}$ because $\Xi'$ is
$L$-closed in $\Sigma$ (cf.~Remark~\ref{rem:closed} iv)).
\end{rem}
\subsection{Further study of the action of the Weyl group}\label{sec:OUT}
As for Q4 in \S\ref{sec:intro}, that is, ``Is $\Out(\Xi)$ realized 
by $W_\Sigma$?" can be answered from the table in \S\ref{sec:table} by
the condition $\#=\#_{\Xi}$ and the answer is ``yes" in most cases in the 
table. We will consider the cases when the answer is ``no", namely, 
we will study the group
$\Out_\Sigma(\Xi)$ in $\Out(\Xi)$ (cf.~Definition~\ref{def:OutX}).
Under the notation in \S\ref{sec:table} we have
\begin{equation}
 \#\bigl(\Out(\Xi)/\!\Out_\Sigma(\Xi)\bigr)=\#/\#_\Xi.
\end{equation}

If $\Sigma$ is of the classical type, it is easy to analyze
$\Out_\Sigma(\Xi)$ because the action of $W_\Sigma$ is easy.
If $\Xi$ is irreducible, $W_\Sigma(\Xi)$ is understood well
by Theorem~\ref{thm} using Dynkin diagrams.
Moreover since $N_{W_{\Sigma}}(\Xi)\subset N_{W_\Sigma}(\Xi+\Xi^\perp)$,
we have
\begin{equation}
 \Out_\Sigma(\Xi)\simeq
 \{g\in\Out_\Sigma(\Xi+\Xi^\perp)\,;\,g(\Xi)=\Xi\}
\end{equation}
by \eqref{eq:OutX} and therefore
the group $\Out_\Sigma(\Xi)$ is described 
by $\Out_\Sigma(\Xi+\Xi^\perp)$.
Hence we may assume $\Xi$ is $\perp$-dense.

\subsubsection{Dual pairs}\label{sec:actW}
If $\Xi$ is not irreducible, $\Out_\Sigma(\Xi)$ may be understood
as a dual pair. 
For example the dual pair $(D_6,2A_1)$ in $E_8$ is special and
the imbedding $D_6+2A_1\subset E_8$ is unique up to the transformations by 
$W_{E_8}$.
Hence there exists $w\in N_{W_{E_8}}(D_6+2A_1)$ which
swaps two $A_1$'s.  Then $w$ always defines a non-trivial element 
of $\Out(D_6)$.  Namely $\Out_\Sigma(\Xi)$ is the diagonal
subgroup of $\Out(D_6+2A_1)$ through the isomorphism 
$\Out(D_6)\simeq\Out(2A_1)$.
The following cases are understood in this way.
\begin{gather*}
 D_m+D_n\subset D_{m+n}\quad(m\ge2,\,n\ge2,\,m\ne4,\,n\ne4),\\
 A_3+2A_1\subset E_6,\\
 2A_3+A_1\subset E_7,\ 
 A_5+A_2\subset E_7,\ 
 D_4+3A_1\subset E_7,\\
 E_6+A_2\subset E_8,\ 
 A_5+A_2+A_1\subset E_8,\ 
 A_4+A_4\subset E_8,\
 D_6+2A_1\subset E_8,\\ 
 D_5+A_3\subset E_8,\ 
 D_4+D_4\subset E_8,\  
 A_2+A_2\subset F_4.
\end{gather*}
For the imbedding $\Xi\subset \Sigma$ in this list, a still more
concrete description of $\Out_\Sigma(\Xi)$ is desirable
if $\Out(\Xi)\not\simeq\mathbb Z/2\mathbb Z\times\mathbb Z/2\mathbb Z$.

For the imbedding $D_m+D_m\simeq D_m^0+D_m^m\subset D_{2m}$
under the notation in \S\ref{sec:table}, the swapping of two 
$D_m$'s under the generators given there is in 
$\Out_{D_{2m}}(D_m^0+D_m^m)$ and therefore $\Out_{D_{2m}}(D_m+D_m)$ is clear.
Similarly for $2A_4\subset E_8$, if we fix $A_3+A_4\subset A_4+A_4\subset E_8$ 
and $A_3+A_4\subset A_8\subset E_8$, we can also specify the swapping of two 
$A_4$'s in $\Out_{E_8}(2A_4)$.
Other cases in the list are described by the study of the imbedding of
$7A_1\subset E_7$ and $8A_1\subset E_8$ through $4A_1\subset D_4$
as is shown later.

\subsubsection{Strongly orthogonal systems of the maximal rank}\label{sec:orth}
Suppose $\Sigma$ is of type $A_n$, $B_n$, $C_n$, $D_n$, $E_6$, $F_4$ or $G_2$
and put 
$m=2[\frac n2]$. 
Under the notation in \S\ref{sec:Dynkin} the strongly 
orthogonal system $\langle\Theta_\Sigma\rangle$ 
of $\Sigma$ with the maximal rank is 
weakly equivalent to
\begin{align}
 \Theta_{A_n} &:= \{\epsilon_1-\epsilon_2,\epsilon_3-\epsilon_4,\ldots,
    \epsilon_{2m-1}-\epsilon_{2m}\},\allowdisplaybreaks\\
 \Theta_{D_n} &:= \{\epsilon_1-\epsilon_2,\epsilon_1+\epsilon_2,
  \epsilon_3-\epsilon_4,\epsilon_3+\epsilon_4,\ldots,
  \epsilon_{2m-1}+\epsilon_{2m}\},\allowdisplaybreaks\\
 \Theta_{B_n} &:= \begin{cases}
             \Theta_{D_n} & (n=2m),\allowdisplaybreaks\\
             \Theta_{D_n}\cup\{\epsilon_n\} &(n=2m+1),
            \end{cases}\allowdisplaybreaks\\
 \Theta_{C_n} &:= \{2\epsilon_1,\ldots,2\epsilon_n\},\allowdisplaybreaks\\
 \Theta_{E_6} &:= \Theta_{F_4} := 
    \{\epsilon_1-\epsilon_2,\epsilon_1+\epsilon_2,
    \epsilon_3-\epsilon_4,\epsilon_3+\epsilon_4\},\allowdisplaybreaks\\
 \Theta_{G_2} &:= \{\epsilon_1-\epsilon_2,\epsilon_1+\epsilon_2-2\epsilon_3\}.
\end{align}
Then $\Out_\Sigma(\Theta_\Sigma):=\Out_\Sigma(\langle\Theta_\Sigma\rangle)$ 
is identified with the subgroup of the permutation group of 
$\Theta_\Sigma$, which is also identified with the permutation group 
$\mathfrak S_{\#\Theta_\Sigma}$ of 
$\{1,\ldots,\#\Theta_\Sigma\}$ according to the expression of 
$\Theta_\Sigma$ by the above ordered set.
Then
\begin{align*}
 \Out_\Sigma(\Theta_\Sigma) &\simeq
 \mathfrak S_{\#\Theta_\Sigma} \qquad(\Sigma \text{ is of type }
 A_n,\,C_n,\,E_6,\,F_4),\\
 \Out_{G_2}(\Theta_{G_2}) &=\{1\},\\
 \Out_{D_n}(\Theta_{D_n}) 
    &=\begin{cases}
       \bigl\langle(1 \ 2)(3 \ 4),(1 \ 3)(2 \ 4),\\
       \quad(1 \ 3\,\cdots\, 2m-1)(2 \ 4\, \cdots\, 2m)\bigr\rangle
       \simeq W_{D_m}&(n=2m),\\
       \bigl\langle(1 \ 2),(1 \ 3)(2 \ 4),\\
       \quad
       (1 \ 3\,\cdots\, 2m-1)(2 \ 4\,\cdots\, 2m)\bigr\rangle
         \simeq W_{B_m}&(n=2m+1).
      \end{cases}
\end{align*}
Here the generators of $\Out_{D_n}(\Theta_{D_n})$ are expressed by 
products of circular permutations.
Note that the group $\Out_{D_n}(\Theta_{D_n})$ is isomorphic to 
$W_{D_m}$ or $W_{B_m}$ if $n=2m$ or $2m+1$, respectively.

\subsubsection{$8A_1\subset E_8$ and $7A_1\subset E_7$}\label{sec:8A1E8}
Since $\Out(8A_1)$ is isomorphic to the symmetric group 
$\mathfrak S_8$, $\Out_{E_8}(8A_1)$ is identified with a subgroup of 
$\mathfrak S_8$.  Since $\#/\#_\Xi=30$, 
$\#\Out_{E_8}(8A_1)= 8!/30=1344=2^6\cdot 3\cdot 7$.
To be more precise, we fix $8A_1\subset E_8$:
\begin{gather}
 \begin{split}
 (8A_1)_o&:=\{\pm\alpha_2,\pm\alpha_3,\pm\alpha_5,\pm\alpha_q\pm\alpha_7,
  \pm\alpha_p,\pm\alpha_0^7,\pm\alpha_0^8\}\subset E_8,\\
 \alpha_2 &=\epsilon_1+\epsilon_2,\ \alpha_3=\epsilon_2-\epsilon_1,\ 
 \alpha_5=\epsilon_4-\epsilon_3,\ \alpha_7=\epsilon_6-\epsilon_5,\\
 \alpha_p &= -\epsilon_5-\epsilon_6,\ \alpha_q = -\epsilon_3-\epsilon_4,\ 
 \alpha_0^8=-\epsilon_7-\epsilon_8,\ \alpha_0^7=\epsilon_7-\epsilon_8.
 \end{split}\label{eq:8A1E8}\allowdisplaybreaks\\
\begin{xy}
 \ar@{.}  *+!D{\alpha_0^7} *{\bullet};
  (5,0)  *+!D{\alpha_1} *{\ast}="A"
 \ar@{.} "A";(10,0) *+!D{\alpha_3} *{\circ}="C"
 \ar@{-} "C";(15,0) *+!D{\alpha_4} *{\circ}="D"
 \ar@{-} "D";(20,0) *+!D{\alpha_5} *{\circ}="E"
 \ar@{-} "E";(25,0) *+!D{\alpha_6} *{\circ}="F"
 \ar@{-} "F";(30,0) *+!D{\alpha_7} *{\circ}="G"
 \ar@{.} "G";(35,0) *+!D{\alpha_8} *{\ast} ="H"
 \ar@{.} "H";(40,0) *+!D{\alpha_0^8} *{\bullet}
 \ar@{-} "D";(15,-5) *+!L{\alpha_2} *{\circ}
 \ar@{-} "F";(25,-5) *+!L{\alpha_p} *{\ccirc}
\end{xy}
\quad
\begin{xy}
 \ar@{.}  *+!D{\alpha_1} *{\circ};
  (5,0)   *+!D{\alpha_3} *{\circ}="A"
 \ar@{-} "A";(10,0) *+!D{\alpha_4} *{\circ}="C"
 \ar@{-} "C";(15,0) *+!D{\alpha_5} *{\circ}="D"
 \ar@{-} "D";(20,0) *+!D{\alpha_6} *{\circ}="E"
 \ar@{-} "E";(25,0) *+!D{\alpha_7} *{\circ}="F"
 \ar@{-} "F";(30,0) *+!D{\alpha_8} *{\circ}="G"
 \ar@{-} "G";(35,0) *+!D{\alpha_0^8} *{\circ}
 \ar@{-} "C";(10,-5) *+!L{\alpha_2} *{\circ}
 \ar@{-} "G";(30,-5) *+!L{\alpha_t} *{\ccirc}
\end{xy}\label{eq:DiagE8}
\end{gather}
Here we used the notation in \S\ref{sec:TE7} and \S\ref{sec:table}.
In particular $E_7\subset E_8$.
Note that
$\alpha_0^8$ and $\alpha_0^7$ are negatives of maximal roots of 
$E_8$ and $E_7$, respectively.
We identify $\mathfrak S_8$ with the permutation group of
the set $\{1,2,3,4,5,6,7,8\}$ of numbers and this ordered set is also
identified with the ordered set given as generators of $(8A_1)_o$
in \eqref{eq:8A1E8}.

Since $(\alpha_0^8)^\perp=E_7$, the left figure above corresponds
to the first diagram in \S\ref{sec:TE7}.
Since 
$(D_8)_o:=\langle\alpha_2,\alpha_3,\ldots,\alpha_8,\alpha_0^8\rangle$
is of type $D_8$, its extended diagram is given in the right figure 
of \eqref{eq:DiagE8} with the negative $\alpha_t$ of its maximal root.
Here we note that $\alpha_0^7<0$ and $\alpha_t>0$.
Since $(D_6)_o^\perp\cap(\alpha_0^8)^\perp\simeq A_1$ by denoting
$(D_6)_o:=\langle\alpha_2,\ldots,\alpha_7\rangle$, we have
$\alpha_t=-\alpha_0^7$. 
Then $\Out_{(D_8)_o}\bigl((8A_1)_o\bigr)$ is generated by
\begin{align}
  g_1 &=(1\ 3\ 5\ 7)(2\ 4\ 6\ 8),\label{eq:E8g1}\\
  g_2 &=(1\ 3)(2\ 4),\\
  g_3 &=(1\ 2)(3\ 4).
\end{align}
Here the generators $g_j$ are expressed by products of cyclic
permutations in $\mathfrak S_8$.
Note that $\#\Out_{(D_8)_o}\bigl((8A_1)_o\bigr)=2^3\cdot 4!=2^6\cdot 3$.

Now we will consider the other diagram in \S\ref{sec:TE7}.
\[
\begin{xy}
 \ar@{-}  *+!D{\alpha_0^7} *{\circ};
  (5,0)   *{\circ}="A"
 \ar@{-} "A";(10,0) *+!D{\alpha_3} *{\ccirc}="C"
 \ar@{.} "C";(15,0) *+!D{\alpha_4} *{\ast}="D"
 \ar@{.} "D";(20,0) *+!D{\alpha_5} *{\bullet}="E"
 \ar@{.} "E";(25,0) *+!D{\alpha_6} *{\ast}="F"
 \ar@{.} "F";(30,0) *+!D{\alpha_7} *{\bullet}="G"
 \ar@{.} "G";(35,0) *+!D{\alpha_8} *{\ast}="H"
 \ar@{.} "H";(40,0) *+!D{\alpha_0^8} *{\bullet}
 \ar@{.} "D";(15,-5) *+!L{\alpha_2} *{\ccirc}
 \ar@{-} "A";(5,-5) *+!L{\alpha_r} *{\circ}
 \ar@{-} "A";(5,5) *+!L{\alpha_s} *{\circ}
 \ar@{}  (70,0) *{\alpha_r=-\alpha_p,\ \alpha_s=-\alpha_q}
\end{xy}
\]
This shows that the element
\begin{equation}
  g_4 = (2\ 4)(6\ 7)
\end{equation}
which corresponds to an element of the 
$W_{\langle\alpha_0^7,\alpha_1,\alpha_3,\alpha_r\rangle}\simeq W_{D_4}$ 
is in $\Out_{E_8}\bigl((8A_1)_o\bigr)$ and not in 
$\Out_{(D_8)_o}\bigl((8A_1)_o\bigr)$.
Then we can conclude that $\Out_{E_8}\bigl((8A_1)_o\bigr)$ 
is generated by $g_j$ $(j=1,2,3,4)$, which is clear by considering 
the order of the groups.

Put $(7A_1)_o = (8A_1)_o\setminus\{\pm\alpha_o^8\}$.
Since $E_7=(\alpha_0^8)^\perp$, it is easy to see that 
$\Out_{E_7}((7A_1)_o)$ is generated by $g_2$, $g_3$, $g_4$ and
\begin{equation}
 g_1' = (1\ 3\ 5)(2\ 4\ 6).\label{eq:E7g1}
\end{equation}
Here we naturally identify $\Out(7A_1)$ with $\mathfrak S_7$
and we have
\begin{align}
 \Out_{E_8}\bigl((8A_1)_o\bigr) &= \langle g_1,g_2,g_3,g_4\rangle,
 & \#\Out_{E_8}\bigl((8A_1)_o\bigr)&=2^6\cdot 3\cdot 7=1344,\\
 \Out_{E_7}\bigl((7A_1)_o\bigr) &= \langle g_1',g_2,g_3,g_4\rangle,
 & \#\Out_{E_7}\bigl((7A_1)_o\bigr)&=2^3\cdot 3\cdot 7=168.
\end{align}
Put $(6A_1)_o=\{\pm\alpha_2,\pm\alpha_3,\pm\alpha_5,
\pm\alpha_q,\pm\alpha_7,\pm\alpha_p\}\subset\{\alpha_0^7,\alpha_0^8\}^\perp
\simeq D_6$ and $(5A_1)_o=\{\pm\alpha_2,\pm\alpha_3,\pm\alpha_5,
\pm\alpha_q,\pm\alpha_7\}\subset\{\alpha_p,\alpha_o^7,\alpha_o^8\}^\perp
\simeq D_4+A_1$. 
Note that $(D_6,2A_1)$ and $(D_4+A_1,3A_1)$ are special dual pairs
in $E_8$ and therefore $\Out_{E_8}(D_6)=\Out(D_6)$ and 
$\Out_{E_8}(D_4+A_1)=\Out(D_4+A_1)$.  
Then we have easily
\begin{align}
 \Out_{E_8}\bigl((7A_1)_o\bigr) &\overset{\sim}{\leftarrow}\Out_{E_7}\bigl((7A_1)_o\bigr),\\
 \Out_{E_8}\bigl((6A_1)_o\bigr) &= \bigl\langle g_1',g_2,(1 \ 2)\bigr\rangle\simeq W_{B_3},
  & \#\Out_{E_8}\bigl((6A_1)_o\bigr) &=48,\\
 \Out_{E_7}\bigl((6A_1)_o\bigr) &= \bigl\langle g_1',g_2,g_3\bigr\rangle\simeq W_{D_3},
  & \#\Out_{E_8}\bigl((6A_1)_o\bigr) &=24,\\
 \Out_{E_8}\bigl((5A_1)_o\bigr) &= \bigl\langle(1 \ 2),(2 \ 3),(3 \ 4)\bigr\rangle\simeq W_{A_3},
  & \#\Out_{E_8}\bigl((5A_1)_o\bigr) &=24,\\
 \Out_{E_7}\bigl((5A_1)_o\bigr) &= \bigl\langle g_2,(1 \ 2)\bigr\rangle\simeq W_{B_2},
  & \#\Out_{E_8}\bigl((5A_1)_o\bigr) &=8.
\end{align}
Put $(4A_1)_o=\{\pm\alpha_2,\pm\alpha_3,\pm\alpha_5,\pm\alpha_q\}$
and $(4A_1)_1 =\{\pm\alpha_2,\pm\alpha_3,\pm\alpha_5,\pm\alpha_7\}$.
Then $(4A_1)_o^\perp\cap E_8\simeq D_4$, 
$(4A_1)_o^{\perp\perp}\simeq D_4$
and $(4A_1)_1^\perp\cap E_8\simeq 4A_1$ and 
$(4A_1)_1^{\perp\perp}=(4A_1)_1$.
%%%
\subsubsection{$2D_4\subset E_8$, $D_4+4A_1\subset E_8$ and $D_4+3A_1\subset E_7$}%
\label{sec:2D4E8}
Retain the notation in the previous section (cf.~\eqref{eq:DiagE8}) 
and put
\begin{align}
  (2D_4)_o &= \langle\alpha_2,\alpha_3,\alpha_4,\alpha_5,\alpha_q\rangle
             +\langle\alpha_7,\alpha_t,\alpha_8,\alpha_0^8,\alpha_p\rangle
             &&\subset E_8,\\
  (D_4+4A_1)_o &= \langle\alpha_2,\alpha_3,\alpha_4,\alpha_5,\alpha_q\rangle
             +\langle\alpha_7,\alpha_t,\alpha_0^8,\alpha_p\rangle
             &&\subset E_8,\\
  (D_4+3A_1)_o &= \langle\alpha_2,\alpha_3,\alpha_4,\alpha_5,\alpha_q\rangle
             +\langle\alpha_7,\alpha_t,\alpha_p\rangle&&\subset E_7.
\end{align}
Then we have the natural identification
\begin{align*}
 & \Out_{E_8}\bigl((2D_4)_o\bigr)\supset
  \Out_{E_7}\bigl((D_4+3A_1)_o\bigr)\\
 &\quad\simeq
 \bigl\{g\in\Out_{E_8}\bigl((8A_1)_o\bigr)\,;\,
   g(\alpha_q)=\alpha_q
   \text{ and }
   g(\alpha_0^8)=\alpha_0^8\bigr\}
\end{align*}
together with \eqref{eq:OAut} 
and therefore $\Out_{E_7}\bigl((D_4+3A_1)_o\bigr)$ is generated by
\begin{equation*}
  (1 \ 2)(5 \ 6) :\ 
    \alpha_2\leftrightarrow\alpha_3,\ 
    \alpha_7\leftrightarrow\alpha_0^7\text{ \ and \ }
  (1 \ 3)(6 \ 7) :\ 
    \alpha_2\leftrightarrow\alpha_5,\ 
    \alpha_7\leftrightarrow\alpha_t.
\end{equation*}
Here the first element corresponds to an element in $W_{(D_8)_o}$
and the second element equals $g_2g_4$.
Moreover $\Out_{E_8}\bigl((2D_4)_o\bigr)$ contains
$(1 \ 5)(2 \ 6)(3 \ 7)(4 \ 8)$ and $\Out_{E_8}\bigl((D_4+4A_1)_o\bigr)$
contains $\Out_{(2D_4)_o}\bigl((D_4+4A_1)_o\bigr)$.  Hence
\begin{align*}
 \Out_{E_7}\bigl((D_4+3A_1)_o\bigr) &=
 \bigl\langle(1 \ 2)(5 \ 6),\ (1 \ 3)(6 \ 7)\bigr\rangle,\\
 \#\Out_{E_7}\bigl((D_4+3A_1)_o\bigr) &= 6,\\
 \Out_{E_8}\bigl((2D_4)_o\bigr) &=
 \bigl\langle(1 \ 2)(5 \ 6),\ (1 \ 3)(6 \ 7),\ 
  (1 \ 5)(2 \ 6)(3 \ 7)(4 \ 8)\bigr\rangle,\\
 \#\Out_{E_8}\bigl((2D_4)_o\bigr) &= 12,\\
 \Out_{E_8}\bigl((D_4+4A_1)_o\bigr) &=
  \bigl\langle(5 \ 6)(7 \ 8),\ (5 \ 8)(6 \ 7),\ 
  \Out_{E_8}\bigl((2D_4)_o\bigr)\bigr\rangle,\\
 \#\Out_{E_8}\bigl((D_4+4A_1)_o\bigr) &= 48.
\end{align*}
%%%%%%%%%%%%%
\subsubsection{$4A_2\subset E_8$, $3A_2\subset E_7$ and 
$3A_2\subset E_6$}\label{sec:4A2E8}
First note that as groups, $\Out(4A_2)$ and $\Out(3A_2)$ are isomorphic
to $W_{B_4}$ and $W_{B_3}$ and their orders of the groups are
$2^4\cdot 4!$ and $2^3\cdot 3!$, respectively.
Fix a representative $4A_2\subset E_8$:
\begin{gather*}
 \begin{aligned}
 \alpha_2 &= \epsilon_1+\epsilon_2,&
 \alpha_0^6 &= -\tfrac12(\epsilon_1+\epsilon_2+\epsilon_3
     +\epsilon_4+\epsilon_5-\epsilon_6-\epsilon_8+\epsilon_8),\\
 \alpha_3 &= \epsilon_2-\epsilon_1,&
 \alpha_1 &= \tfrac12(\epsilon_1+\epsilon_8)
  -\tfrac12(\epsilon_2+\epsilon_3+\epsilon_4+\epsilon_5
  +\epsilon_6+\epsilon_7),\\
 \alpha_5 &= \epsilon_4-\epsilon_3,&
 \alpha_6 &= \epsilon_5-\epsilon_4,\\
 \alpha_8 &= \epsilon_7-\epsilon_6,&
 \alpha_0^8 &= -\epsilon_7-\epsilon_8,
 \end{aligned}\\[-32pt]
\begin{aligned}
 (4A_2)_0 &=\bigl\langle\{\alpha_2,\alpha_0^6,\alpha_3,\alpha_1,\alpha_5,\alpha_6,\alpha_8,\alpha_0^8\}\bigr\rangle,\\
 (3A_2)_0 &=\bigl\langle\{\alpha_2,\alpha_0^6,\alpha_3,\alpha_1,\alpha_5,\alpha_6\}\bigr\rangle.
\end{aligned}\qquad%\\
\raisebox{22pt}{\begin{xy}
 \ar@{.}  *+!D{\alpha_4}  *{\ast}="O";
  (-4.33,2.5)  *+!D{\alpha_3}  *{\bullet}="C"
 \ar@{-} "C";(-8.66,5) *+!D{\alpha_1} *{\bullet}="D"
 \ar@{.} "O";(4.33,2.5) *+!D{\alpha_5}  *{\bullet}="E"
 \ar@{-} "E";(8.66,5) *+!D{\alpha_6}  *{\bullet}="Y"
 \ar@{.} "O";(0,-5) *+!L{\alpha_2}  *{\bullet}="X"
 \ar@{-} "X";(0,-10) *+!L{\alpha_0^6} *{\bullet}
 \ar@{.} "Y";(13.66,5) *+!D{\alpha_7} *{\ast}="Z"
 \ar@{.} "Z";(18.66,5) *+!D{\alpha_8} *{\bullet}="W"
 \ar@{-} "W";(23.66,5) *+!D{\alpha_0^8} *{\bullet}
\end{xy}}
\end{gather*}
Then the permutation group of the 8 generators of $(4A_2)_o$ is identified
with $\mathfrak S_8$ as in the case of $(8A_1)_o\subset E_8$.
Then $\Out\bigl((4A_2)_o\bigr)=\langle g_1,g_2,g_3\rangle$ and $\Out(3(A_2)_o\bigr)
=\langle g_1',g_2,g_3\rangle$ (cf.~\eqref{eq:E8g1}--\eqref{eq:E7g1}).
Note that
\begin{align*}
  \#\Out_{E_6}\bigl((3A_2)_o\bigr)&=\#\Out\bigl((3A_2)_o\bigr)/8=6,\\
  \#\Out_{E_7}\bigl((3A_2)_o\bigr)&=\#\Out\bigl((3A_2)_o\bigr)/4=12,\\
  \#\Out_{E_8}\bigl((4A_2)_o\bigr)&=\#\Out\bigl((4A_2)_o\bigr)/8=48,\\
  \Out_{E_6}\bigl((3A_2)_o\bigr)&\subset\Out_{E_8}\bigl((4A_2)_o\bigr).
\end{align*}
Since $(2A_2,2A_2)$ is a special dual pair in $E_8$, we have
\begin{equation}\label{eq:2A2inE8}
 \bigl(\exists 2A_2\subset (4A_2)_o\text{ such that }
 w|_{2A_2}=id\bigr)\ \Rightarrow w=id
\end{equation}
for $w\in\Out_{E_8}\bigl((4A_2)_o\bigr)$.
We will choose elements in $\Out_{E_8}\bigl((2A_2)_o\bigr)$.
The rotation of the extended Dynkin diagram of $E_6$ comes 
from $W_{E_6}$ and therefore the element
\begin{equation}
 h_1 = (1\ 3\ 5)(2\ 4\ 6)
\end{equation}
is contained in $\Out_{E_6}\bigl((3A_2)_o\bigr)$.
The argument \S\ref{sec:TE6} shows that in view of the 
transformation of an element of 
$W_{\langle\alpha_0^8,\alpha_2,\alpha_4,\alpha_5,\alpha_6\rangle}$,
$(1 \ 6)(2 \ 5)$ or $(1 \ 6)(2 \ 5)(3 \ 4)$ should be in 
$W_{E_6}\bigl((3A_2)_o\bigr)$.  Owing to \eqref{eq:2A2inE8}, 
we can conclude that
\begin{equation}
 h_2 = (1 \ 6)(2 \ 5)(3 \ 4)
\end{equation}
is contained in the group and 
$\Out_{E_6}\bigl((3A_2)_o\bigr)=\langle h_1,h_2\rangle$.

Since $\overline\Hom(E_6,E_7)=1$, 
$\Out_{E_7}\bigl((3A_2)_o\bigr)$ contains
\begin{equation}
 h_3 = (3 \ 5)(4 \ 6).
\end{equation}
Considering in 
$(A_8)_o=\langle\alpha_1,\alpha_3,\alpha_4,\alpha_5,\alpha_6,\alpha_7,\alpha_8,\alpha_0^8\rangle\simeq A_8$,
there is an element $w\in W_{(A_8)_o}$
such that
\[
 w(\alpha_1)=\alpha_1,\ 
 w(\alpha_3)=\alpha_3\ 
 w(\alpha_4)=\alpha_6,\ 
 w(\alpha_5)=\alpha_7,\ 
 w(\alpha_6)=\alpha_4,\ 
 w(\alpha_7)=\alpha_5.
\]
Then it also follows from \eqref{eq:2A2inE8} that
\begin{equation}
  h_4 = (1 \ 2)(5 \ 7)(6 \ 8)
\end{equation}
is in $\Out_{E_8}\bigl((4A_2)_o\bigr)$.
Calculating the order of the group, we have
\begin{align}
 \Out_{E_6}\bigl((3A_2)_o\bigr) &= \langle h_1,h_2\rangle, &
 &\#\Out_{E_6}\bigl((3A_2)_o\bigr)=6,\\
 \Out_{E_7}\bigl((3A_2)_o\bigr) &= \langle h_1,h_2,h_3\rangle, &
 &\#\Out_{E_7}\bigl((3A_2)_o\bigr)=12,\\
 \Out_{E_8}\bigl((4A_2)_o\bigr) &= \langle h_1,h_2,h_4\rangle, &
 &\#\Out_{E_8}\bigl((4A_2)_o\bigr)=48.
\end{align}

%%%%%%%%%%%%%%%%%%%%%%%%%%  Classical  %%%%%%%%%%%%%%%%%%%%%%%%%%%%%%
\section{List of irreducible root systems}\label{sec:Dynkin}
%%%%%  A_n %%%%%
%%%
$A_n$\quad
\begin{xy}
 \ar@{-}  *++!D{\alpha_1} *++!U{{}_1} *\cir<3pt>{}="A";
  (10,0)  *++!D{\alpha_2} *++!U{{}_1} *\cir<3pt>{}="B"
 \ar@{-} "B";(20,0) *++!D{\alpha_3} *++!U{{}_1} *\cir<3pt>{}="C"
 \ar@{-} "C";(25,0) \ar@{.} (25,0);(30,0)^*!U{\cdots}
 \ar@{-} (30,0);(35,0) *+<10pt>!D{\alpha_{n-2}} *++!U{{}_1} *\cir<3pt>{}="H"
 \ar@{-} "H";(45,0) *+<10pt>!D{\alpha_{n-1}} *++!U{{}_1} *\cir<3pt>{}="I"
 \ar@{-} "I";(55,0) *++!D{\alpha_{n}} *++!U{{}_1} *\cir<3pt>{}="J"
 \ar@{-} "A";(27.5,-11) *++!D{{}_1} *!L++!U{\alpha_0} 
          *\cir<1.5pt>{} *\cir<3pt>{}="H"
 \ar@{-} "H";"J"
\end{xy}
\qquad
$A_1$\quad
\begin{xy}
  \ar@{=}
% \ar@<.7mm>@{=} 
 *++!D{\alpha_0} *++!U{{}_1} *\cir<1.5pt>{} *\cir<3pt>{}="A";
  (10,0) *++!D{\alpha_1} *++!U{{}_1} *\cir<3pt>{}="B"
% \ar@<.7mm>@{=} "B";"A"
\end{xy}
\begin{align*}
 \Sigma &= \bigl\{\pm(\epsilon_i-\epsilon_j)\,;\,1\le i<j\le n+1\bigr\},
 \quad\#\Sigma = n(n+1),\\
 \alpha_j &= \epsilon_j - \epsilon_{j+1\ (j=1,\ldots,n)},\ 
 \alpha_0 = \epsilon_{n+1} - \epsilon_1,\quad
 \#W = (n+1)!.
\end{align*}
%%%%%

\smallskip
%%%%%  B_n %%%%%
$B_n$\quad
\begin{xy}
 \ar@{-}  *++!D{\alpha_1} *++!U{{}_1} *\cir<3pt>{} ;
  (10,0)  *++!D{\alpha_2} *+!L+!U{{}_2} *\cir<3pt>{}="B"
 \ar@{-} "B";(20,0) *++!D{\alpha_3} *++!U{{}_2} *\cir<3pt>{}="C"
 \ar@{-} "C";(25,0) \ar@{.} (25,0);(30,0)^*!U{\cdots}
 \ar@{-} (30,0);(35,0) *+<10pt>!D{\alpha_{n-2}} *++!U{{}_2} *\cir<3pt>{}="H"
 \ar@{-} "H";(45,0) *+<10pt>!D{\alpha_{n-1}} *++!U{{}_2} *\cir<3pt>{}="I"
 \ar@{=>} "I";(55,0) *++!D{\alpha_{n}} *++!U{{}_2} *\cir<3pt>{}
 \ar@{-} "B";(10,-10) *++!R{{}_1} *++!L{\alpha_0} *\cir<1.5pt>{} *\cir<3pt>{}
\end{xy}
\qquad$B_2=C_2$
\begin{xy}
 \ar@{=>}  *++!D{\alpha_0} *++!U{{}_1} *\cir<1.5pt>{} *\cir<3pt>{} ;
  (10,0)  *++!D{\alpha_2} *++!U{{}_2} *\cir<3pt>{}="A"
 \ar@{<=} "A";(20,0) *++!D{\alpha_1} *++!U{{}_1} *\cir<3pt>{}="B"
\end{xy}
\begin{align*}
 \Sigma &= \bigl\{\pm(\epsilon_i-\epsilon_j),\ 
 \pm(\epsilon_i+\epsilon_j),\ \pm\epsilon_k \,;\,1\le i<j\le n,\ 1\le k\le n\bigr\},\quad\#\Sigma=2n^2,\\
 \alpha_j &= \epsilon_j - \epsilon_{j+1\ (j=1,\ldots,n-1)},\ 
 \alpha_n = \epsilon_n,\ 
 \alpha_0 = -\epsilon_1 - \epsilon_2,\ 
 \alpha_0'=-\epsilon_1,\quad
 \#W = 2^n\cdot n!.
\end{align*}
%%%%%

\smallskip
%%%%% C_n %%%%%
$C_n$\quad
\begin{xy}
 \ar@{=>}  *++!U{{}_1} *++!D{\alpha_0} *\cir<1.5pt>{}*\cir<3pt>{} ;
  (10,0)  *++!D{\alpha_1} *++!U{{}_2} *\cir<3pt>{}="A"
 \ar@{-} "A";(20,0) *++!D{\alpha_2} *++!U{{}_2} *\cir<3pt>{}="B"
 \ar@{-} "B";(25,0) \ar@{.} (25,0);(30,0)^*!U{\cdots}
 \ar@{-} (30,0);(35,0) *+<10pt>!D{\alpha_{n-2}} *++!U{{}_2} *\cir<3pt>{}="H"
 \ar@{-} "H";(45,0) *+<10pt>!D{\alpha_{n-1}} *++!U{{}_2} *\cir<3pt>{}="I"
 \ar@{<=} "I";(55,0) *++!D{\alpha_{n}} *++!U{{}_1} *\cir<3pt>{}
\end{xy}
\begin{align*}
 \Sigma &= \bigl\{\pm(\epsilon_i-\epsilon_j),\ 
 \pm(\epsilon_i+\epsilon_j),\ \pm2\epsilon_k
  \,;\,1\le i<j\le n,\ 1\le k\le n\bigr\},\quad\#\Sigma=2n^2,\\
 \alpha_j &= \epsilon_j - \epsilon_{j+1\ (j=1,\ldots,n-1)},\ 
 \alpha_n = 2\epsilon_n,\ 
 \alpha_0 = -2\epsilon_1,\ 
 \alpha_0'=-\epsilon_1-\epsilon_2,
 \  \#W = 2^n\cdot n!.
\end{align*}
%%%%%

\smallskip
%%%%% D_n %%%%%
$D_n$\quad
\begin{xy}
 \ar@{-}  *++!D{\alpha_1} *++!U{{}_1} *\cir<3pt>{} ;
  (10,0)  *++!D{\alpha_2} *+!L+!U{{}_2} *\cir<3pt>{}="B"
 \ar@{-} "B";(20,0) *++!D{\alpha_3} *++!U{{}_2} *\cir<3pt>{}="C"
 \ar@{-} "C";(25,0) \ar@{.} (25,0);(30,0)^*!U{\cdots}
 \ar@{-} (30,0);(35,0) *+<10pt>!D{\alpha_{n-3}} *++!U{{}_2} *\cir<3pt>{}="H"
 \ar@{-} "H";(45,0) *+<10pt>!D{\alpha_{n-2}} *+!L+!U{{}_2} *\cir<3pt>{}="I"
 \ar@{-} "I";(55,0) *+<10pt>!D{\alpha_{n-1}} *++!U{{}_1} *\cir<3pt>{}
 \ar@{-} "B";(10,-10) *++!R{{}_1} *++!L{\alpha_0} *\cir<1.5pt>{} *\cir<3pt>{}
 \ar@{-} "I";(45,-10) *++!L{\alpha_n} *++!R{{}_1} *\cir<3pt>{}
\end{xy}
\qquad
$D_4$
\begin{xy}
 \ar@{-}  *++!D{\alpha_1} *++!U{{}_1} *\cir<3pt>{} ;
  (10,0)  *+!L+!D{\alpha_2} *+!L+!U{{}_2} *\cir<3pt>{}="B"
 \ar@{-} "B";(20,0) *++!D{\alpha_3} *++!U{{}_1} *\cir<3pt>{}="C"
 \ar@{-} "B";(10,-10) *++!L{\alpha_4} *++!R{{}_1} *\cir<3pt>{}
 \ar@{-} "B";(10,10) *++!L{\alpha_0} *++!R{{}_1} *\cir<1.5pt>{} *\cir<3pt>{}
\end{xy}
\begin{align*}
 \Sigma &= \bigl\{\pm(\epsilon_i-\epsilon_j),\ 
 \pm(\epsilon_i+\epsilon_j)\,;\,1\le i<j\le n\bigr\},
 \quad\#\Sigma=2n(n-1),\\
 \alpha_j &= \epsilon_j - \epsilon_{j+1\quad(j=1,\ldots,n-1)},\qquad
 \alpha_n = \epsilon_{n-1}+\epsilon_n,\\
 \alpha_0 &= -\epsilon_1-\epsilon_2,\quad \#W=2^{n-1}\cdot n!.
\end{align*}
%%%%%
%%%%%%%%%%%%%%%%%%%%%%%%%  Exceptional  %%%%%%%%%%%%%%%%%%%%%%%%%%%%%

\medskip
%%%%% E_6 %%%%%
$E_6$\quad
\begin{xy}
 \ar@{-}  *++!D{\alpha_1} *++!U{{}_1} *\cir<3pt>{} ;
  (10,0)  *++!D{\alpha_3} *++!U{{}_2} *\cir<3pt>{}="C"
 \ar@{-} "C";(20,0) *++!D{\alpha_4} *+!L+!U{{}_3} *\cir<3pt>{}="D"
 \ar@{-} "D";(30,0) *++!D{\alpha_5} *++!U{{}_2} *\cir<3pt>{}="E"
 \ar@{-} "E";(40,0) *++!D{\alpha_6} *++!U{{}_1} *\cir<3pt>{},
 \ar@{-} "D";(20,-10) *++!L{\alpha_2} *++!R{{}_2} *\cir<3pt>{}="X"
 \ar@{-} "X";(20,-20) *++!L{\alpha_0} *++!R{{}_1} *\cir<1.5pt>{} *\cir<3pt>{};
\end{xy}
\begin{align*}
 \Sigma &= \bigl\{\pm(\epsilon_i-\epsilon_j),\ 
 \pm(\epsilon_i+\epsilon_j),\ 
  \pm\tfrac12(\epsilon_8-\epsilon_7-\epsilon_6
   +\textstyle\sum_{k=1}^5(-1)^{\nu(k)}\epsilon_k)\,;\\
  &\qquad
   \,1\le i<j\le 5,\ \textstyle\sum_{k=1}^5\nu(k)\text{ is even}\bigr\},
  \quad\#\Sigma=72,\\
 \alpha_1 &= \tfrac12(\epsilon_1 + \epsilon_8)
     -\tfrac12(\epsilon_2+\epsilon_3+\epsilon_4+\epsilon_5
     +\epsilon_6 + \epsilon_7),\\
 \alpha_2 &= \epsilon_1+\epsilon_2,\quad
 \alpha_j  = \epsilon_{j-1} - \epsilon_{j-2}\quad(3\le j\le 6),\\
 \alpha_0 &= -\tfrac12(\epsilon_1+\epsilon_2+\epsilon_3+\epsilon_4
             +\epsilon_5-\epsilon_6-\epsilon_7+\epsilon_8),\quad
             \#W=2^7\cdot 3^4\cdot 5.
\end{align*}
%%%%%

\smallskip
%%%%% E_7 %%%%%
$E_7$\quad
\begin{xy}
 \ar@{-}  *++!U{{}_1} *++!D{\alpha_0} *\cir<1.5pt>{} *\cir<3pt>{};
  (10,0)  *++!D{\alpha_1} *++!U{{}_2} *\cir<3pt>{}="A"
 \ar@{-} "A";(20,0) *++!D{\alpha_3} *++!U{{}_3} *\cir<3pt>{}="C"
 \ar@{-} "C";(30,0) *++!D{\alpha_4} *+!L+!U{{}_4} *\cir<3pt>{}="D"
 \ar@{-} "D";(40,0) *++!D{\alpha_5} *++!U{{}_3} *\cir<3pt>{}="E"
 \ar@{-} "E";(50,0) *++!D{\alpha_6} *++!U{{}_2} *\cir<3pt>{}="F"
 \ar@{-} "F";(60,0) *++!D{\alpha_7} *++!U{{}_1} *\cir<3pt>{},
 \ar@{-} "D";(30,-10) *++!L{\alpha_2} *++!R{{}_2} *\cir<3pt>{}
\end{xy}
\begin{align*}
 \Sigma &= \bigl\{\pm(\epsilon_i-\epsilon_j),\ 
 \pm(\epsilon_i+\epsilon_j),\ \pm(\epsilon_7-\epsilon_8),\ 
  \pm\tfrac12(\epsilon_7-\epsilon_8
   +\textstyle\sum_{k=1}^6(-1)^{\nu(k)}\epsilon_k)\,;\\
  &\qquad
   \,1\le i<j\le 6,\ \textstyle\sum_{k=1}^6\nu(k)\text{ is odd}\bigr\},
  \quad\#\Sigma=126,\\
 \alpha_1 &= \tfrac12(\epsilon_1 + \epsilon_8)
     -\tfrac12(\epsilon_2+\epsilon_3+\epsilon_4+\epsilon_5
     +\epsilon_6 + \epsilon_7),\\
 \alpha_2 &= \epsilon_1+\epsilon_2,\quad
 \alpha_j  = \epsilon_{j-1} - \epsilon_{j-2}\quad(3\le j\le 7),\\
 \alpha_0 &= \epsilon_7-\epsilon_8,\quad
 \#W = 2^{10}\cdot 3^4\cdot 5\cdot 7.
\end{align*}
%%%%%

\smallskip
%%%%% E_8 %%%%%
$E_8$\quad
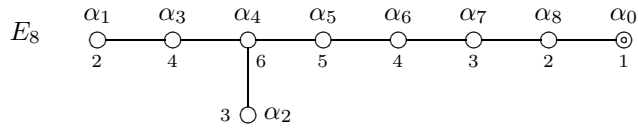
\begin{xy}
 \ar@{-}  *++!D{\alpha_1} *++!U{{}_2} *\cir<3pt>{} ;
  (10,0)  *++!D{\alpha_3} *++!U{{}_4} *\cir<3pt>{}="C"
 \ar@{-} "C";(20,0) *++!D{\alpha_4} *+!L+!U{{}_6} *\cir<3pt>{}="D"
 \ar@{-} "D";(30,0) *++!D{\alpha_5} *++!U{{}_5} *\cir<3pt>{}="E"
 \ar@{-} "E";(40,0) *++!D{\alpha_6} *++!U{{}_4} *\cir<3pt>{}="F"
 \ar@{-} "F";(50,0) *++!D{\alpha_7} *++!U{{}_3} *\cir<3pt>{}="F"
 \ar@{-} "F";(60,0) *++!D{\alpha_8} *++!U{{}_2} *\cir<3pt>{}="G"
 \ar@{-} "G";(70,0) *++!D{\alpha_0} *++!U{{}_1} *\cir<1.5pt>{}
  *\cir<3pt>{},
 \ar@{-} "D";(20,-10) *++!L{\alpha_2} *++!R{{}_3} *\cir<3pt>{}
\end{xy}
\begin{align*}
 \Sigma &= \bigl\{\pm(\epsilon_i-\epsilon_j),\ 
 \pm(\epsilon_i+\epsilon_j),\ 
  \tfrac12\textstyle\sum_{k=1}^8(-1)^{\nu(k)}\epsilon_k\,;\,
%  \qquad
  1\le i<j\le 8,\\ &\qquad \textstyle\sum_{k=1}^8\nu(k)\text{ is even}\bigr\},
  \quad\#\Sigma=240,\\
 \alpha_1 &= \tfrac12(\epsilon_1 + \epsilon_8)
     -\tfrac12(\epsilon_2+\epsilon_3+\epsilon_4+\epsilon_5
     +\epsilon_6 + \epsilon_7),\\
 \alpha_2 &= \epsilon_1+\epsilon_2,\quad
 \alpha_j  = \epsilon_{j-1} - \epsilon_{j-2}\quad(3\le j\le 8),\\
 \alpha_0 &= -\epsilon_7-\epsilon_8,\quad
 \#W = 2^{14}\cdot 3^5\cdot 5^2\cdot 7.
\end{align*}
%%%%%

\smallskip
%%%%% F_4 %%%%%
$F_4$\quad
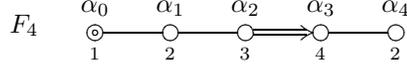
\begin{xy}
 \ar@{-}  *++!D{\alpha_0} *++!U{{}_1} *\cir<1.5pt>{} *\cir<3pt>{};
  (10,0)  *++!D{\alpha_1} *++!U{{}_2} *\cir<3pt>{}="A"
 \ar@{-} "A";(20,0) *++!D{\alpha_2} *++!U{{}_3} *\cir<3pt>{}="B"
 \ar@{=>} "B";(30,0) *++!D{\alpha_3} *++!U{{}_4} *\cir<3pt>{}="C"
 \ar@{-} "C";(40,0) *++!D{\alpha_4} *++!U{{}_2} *\cir<3pt>{}
\end{xy}
\begin{align*}
 \Sigma &= \bigl\{\pm(\epsilon_i-\epsilon_j),\ 
  \pm(\epsilon_i+\epsilon_j),\ 
  \pm\epsilon_k,\ 
  \pm\tfrac12(\epsilon_1\pm\epsilon_2\pm\epsilon_3\pm\epsilon_4)\,;\\
   &\qquad
   1\le i<j\le 4,\ 1\le k\le 4\bigr\},\quad\#\Sigma=48,\\
 \alpha_1 &= \epsilon_2-\epsilon_3,\ 
 \alpha_2 = \epsilon_3-\epsilon_4,\ 
 \alpha_3 = \epsilon_4,\ 
 \alpha_4 = \tfrac12(\epsilon_1-\epsilon_2-\epsilon_3-\epsilon_4),\\
 \alpha_0 &= -\epsilon_1-\epsilon_2,\ 
 \alpha_0' = -\epsilon_1,
 \quad \#W = 2^7\cdot 3^2.
\end{align*}%%%%%

\smallskip
%%%%% G_2 %%%%%
$G_2$\quad
\begin{xy}
 \ar@{-}  *++!D{\alpha_0} *++!U{{}_1} *\cir<1.5pt>{} *\cir<3pt>{};
  (10,0)  *++!D{\alpha_1} *++!U{{}_2} *\cir<3pt>{}="A"
 \ar@3{->} "A";(20,0) *++!U{{}_3} *++!D{\alpha_2} *\cir<3pt>{}
\end{xy}
\begin{align*}
 \Sigma &= \bigl\{\pm(\epsilon_i-\epsilon_j),\ 
  \mp(2\epsilon_1-\epsilon_2-\epsilon_3),\ 
  \mp(2\epsilon_2-\epsilon_1-\epsilon_3),\ 
  \pm(2\epsilon_3-\epsilon_1-\epsilon_2)\,;\\
 &\qquad 
   1\le i<j\le 3\bigr\},\quad\#\Sigma=12,\\
 \alpha_1 &= -2\epsilon_1+\epsilon_2+\epsilon_3,\ 
 \alpha_2 = \epsilon_1-\epsilon_2,\ 
 \alpha_0 = \epsilon_1+\epsilon_2-2\epsilon_3,\ 
 \alpha_0' = \epsilon_2-\epsilon_3,\ \ 
 \#W=12.
\end{align*}
\begin{rem} {\rm i) }
There are natural identifications of root systems
\begin{gather}
 D_1=\emptyset,\quad D_2\simeq A_1+A_1,\quad D_3\simeq A_3,\quad
 A_1\simeq B_1\simeq C_1,\quad B_2\simeq C_2\\
\intertext{and Weyl groups}
\begin{gathered}
 \begin{matrix}
 \mathfrak S_n\ltimes(\mathbb Z/2\mathbb Z)^n & \overset{\sim}{\to} &
 W_{B_n}=W_{C_n}\\
 \vin & & \vin\\
 \bigl(\sigma,(c_1,\dots,c_n)\bigr)& \mapsto &
 w_{\sigma,c}:
  \mathbb R^n\ni\epsilon_j\mapsto (-1)^{c_j}\epsilon_{\sigma(j)}
  \quad(j=1,\dots,n),
 \end{matrix}\\
 W_{D_n} = \{w_{\sigma,c}\in W_{B_n}\,;\,(-1)^{\sum c_j}=1\},\\
 W_{A_{n-1}} = \{w_{\sigma,c}\in W_{B_n}\,;\,c_1=\dots=c_n=0\}.
\end{gathered}
\end{gather}
{\rm ii) }
For the fundamental system $\Psi$ of an irreducible root 
system $\Sigma$ we have
\begin{align}
 \# W_\Psi &= \#\Sigma^L\cdot
 \#W_{\Psi\cap\alpha_0^\perp}\label{eq:NW}\\
 &=(\#\Psi)!\cdot
  \#\{\alpha_j\in\tilde\Psi\,;\,m_j(\alpha_{max})=1\}
  \cdot\prod_{\alpha_j\in\tilde\Psi}m_j(\alpha_{max}).
\end{align}
{\rm iii) }
There exist roots $\alpha_i\in\Psi$ satisfying $m_i(\alpha_{max})=1$,
which are determined by \eqref{eq:mult1}. 
Moreover we have
\begin{equation}\label{eq:sum}
 m_j(\alpha_{max}) =\sum_{\alpha_\nu\in\tilde\Psi\setminus\{\alpha_j\}}
 m_\nu(\alpha_{max})
 \Bigl(-\frac{(\alpha_{\nu}|\alpha_j)}{(\alpha_j|\alpha_j)}\Bigr)
 \qquad(\alpha_j\in\tilde\Psi)
\end{equation}
because of \eqref{eq:zerosum}.
In particular
\begin{equation}
 2m_j(\alpha_{\max}) = 
 \sum_{\substack{\alpha_\nu\in\tilde\Psi\setminus\{\alpha_j\}
 \\(\alpha_\nu|\alpha_j)\ne0}}
                      m_\nu(\alpha_{max})
 \qquad(\alpha_j\in\tilde\Psi\cap\Sigma^L).
\end{equation}
Note that these conditions determine the extended Dynkin diagrams
with the numbers $\{m_j(\alpha_{\max})\,;\,\alpha_j\in\tilde\Psi\}$
as in the following proposition, from which the classification
of the root systems follows.
\end{rem}
%
%For the completeness of this note we will give a proof of the 
%classification of root systems or extended Dynkin diagrams.  
%The proof here based on \eqref{eq:sum} is easy and seems to be known.

\begin{defn}\label{defn:affineDy}
A diagram $\mathcal G$ consisting of finite vertices and lines 
and/or arrows is an \textsl{affine Dynkin diagram} 
if $\mathcal G$ satisfies the following conditions.
Each line or arrow links a vertex to another vertex and each arrow has 
a stem formed by multiple lines.  
Moreover each vertex has an attached positive real number with the following 
property.

Fix any vertex $P$ in $\mathcal G$ and let $m$ be the number attached to $P$.
Let $L_1,\dots,L_p$ be the lines and arrows linking $P$ to other vertices.
We denote the vertices and their attached numbers by 
$Q_1,\dots,Q_p$ and $m_1,\dots,m_p$, respectively.
%.
Then
\begin{gather}\label{eq:ExtNumber}
  2 m = k_1m_1+\cdots+k_pm_p
\intertext{by putting}
   k_j = 
  \begin{cases}
    1 \text{ if $L_j$ is an arrow starting from $P$ or a line},\\
    \text{the number ($\ge 2$) of lines in $L_j$ 
    if $L_j$ is an arrow pointing toward $P$}
  \end{cases}\notag
\end{gather}
and the minimal number attached in any connected component of $\mathcal G$
equals $1$.
\end{defn}
Then we have the following proposition, which is probably known.
Its proof is elementary and easy and %therefore 
we give it for the completeness.
%%%
\begin{prop}\label{prop:clasify}
The connected affine Dynkin diagram $\mathcal G$ is $\tilde R$ with an 
irreducible root system of type $R$ or one of the following diagrams
(cf.~Remark~\ref{rem:cls}~iii)).
\noindent
\iftrue
$\tilde B_n'\!\!$~
\begin{xy}
 \ar@{<=} *+!U{-\epsilon_1} *+!D{{}_1} *{\ccirc} ; (5,0) *+!D{{}_1} *{\circ}="A"
 \ar@{-} "A";(10,0) *+!D{{}_1} *{\circ}="B"
 \ar@{-} "B";(12.5,0) \ar@{.} (12.5,0);(15,0)^*!U{{}_{\cdots}}
 \ar@{-} (15,0);(17.5,0) *+!D{{}_1} *{\circ}="H"
 \ar@{-} "H";(22.5,0) *+!D{{}_1} *{\circ}="I"
 \ar@{=>} "I";(27.5,0) *+!D{{}_1} *{\circ}
\end{xy}
\quad\ 
$\tilde C_n'\ (n\ge3)$~
\begin{xy}
 \ar@{-} *+!D{{}_1}  *{\circ} ; (5,0) *+!D{{}_2} *{\circ}="B"
 \ar@{-} "B";(10,0) *+!D{{}_2} *{\circ}="C"
 \ar@{-} "C";(12.5,0) \ar@{.} (12.5,0);(15,0)^*!U{{}_{\cdots}}
 \ar@{-} (15,0);(17.5,0) *+!D{{}_2} *{\circ}="H"
 \ar@{-}  "H";(22.5,0) *+!D{{}_2} *{\circ}="I"
 \ar@{<=} "I";(27.5,0) *+!D{{}_1} *{\circ}
 \ar@{-}  "B";(5,-5) *+!R{{}_1} *+!L{-\epsilon_1-\epsilon_2} *{\ccirc}
\end{xy}
\quad\ 
$\tilde C_2'\simeq\tilde B_2'$~
\begin{xy}
 \ar@{<=} *+!D{{}_1} *{\ccirc} ; (5,0)  *+!D{{}_1} *{\circ}="A"
 \ar@{=>} "A";(10,0) *+!D{{}_1} *{\circ}="B"
\end{xy}\\
$\tilde F_4'$~
\begin{xy}
 \ar@{-} *+!D{{}_1} *{\circ}; (5,0)  *+!D{{}_2} *{\circ}="A"
 \ar@{=>} "A";(10,0) *+!D{{}_3} *{\circ}="B"
 \ar@{-} "B";(15,0) *+!D{{}_2} *{\circ}="C"
 \ar@{-} "C";(20,0) *+!D{{}_1} *+!U{-\epsilon_1} *{\ccirc}
\end{xy}
%%%%%
\ \ 
%%%%% G_2 %%%%%
$\tilde G_2'$~
\begin{xy}
 \ar@3{->} *+!D{{}_1} *{\circ};
  (5,0)  *+!D{{}_2} *{\circ}="A"
 \ar@{-} "A";(10,0) *+!D{{}_1} 
  *+!U{\epsilon_2-\epsilon_3} *{\ccirc}
\end{xy}\ \ \ 
$\widetilde{BC}_n\!\!\!$~
\begin{xy}
 \ar@{=>} *+!D{{}_1} *+!U{-2\epsilon_1} *{\ccirc} ; (5,0) *+!D{{}_2} *{\circ}="B"
 \ar@{-} "B";(10,0) *+!D{{}_2} *{\circ}="C"
 \ar@{-} "C";(12.5,0) \ar@{.} (12.5,0);(15,0)^*!U{{}_{\cdots}}
 \ar@{-} (15,0);(17.5,0) *+!D{{}_2} *{\circ}="H"
 \ar@{-}  "H";(22.5,0) *+!D{{}_2} *{\circ}="I"
 \ar@{=>} "I";(27.5,0) *+!D{{}_2} *{\circ}
\end{xy}
\ \ \ 
$\widetilde{BC}_1$~
\begin{xy}
 \ar@<.7mm>@{=} *+!D{{}_1} *{\ccirc}="A";
  (5,0) *+!D{\;{}_2} *{\!\text{\Large$\rangle$}\!\circ}="B"
 \ar@<.7mm>@{=} "B";"A"
\end{xy}
\else
\tilde B_n'$
\begin{xy}
 \ar@{<=}  *++!D{\alpha_0} *++!U{-\epsilon_1} *\cir<1.5pt>{}*\cir<3pt>{} ;
  (10,0)  *++!D{\alpha_1} *++!U{{}_1} *\cir<3pt>{}="A"
 \ar@{-} "A";(20,0) *++!D{\alpha_2} *++!U{{}_1} *\cir<3pt>{}="B"
 \ar@{-} "B";(25,0) \ar@{.} (25,0);(30,0)^*!U{\cdots}
 \ar@{-} (30,0);(35,0) *+<12pt>!D{\alpha_{n-2}} *++!U{{}_1} *\cir<3pt>{}="H"
 \ar@{-} "H";(45,0) *+<12pt>!D{\alpha_{n-1}} *++!U{{}_1} *\cir<3pt>{}="I"
 \ar@{=>} "I";(55,0) *++!D{\alpha_{n}} *++!U{{}_1} *\cir<3pt>{}
\end{xy}\\
$A_{2n-1}^{(2)}=\tilde C_n'$\!\!\!\!\!\!
\begin{xy}
 \ar@{-}  *++!D{\alpha_1} *++!U{{}_1} *\cir<3pt>{} ;
  (10,0)  *++!D{\alpha_2} *+!L+!U{{}_2} *\cir<3pt>{}="B"
 \ar@{-} "B";(20,0) *++!D{\alpha_3} *++!U{{}_2} *\cir<3pt>{}="C"
 \ar@{-} "C";(25,0) \ar@{.} (25,0);(30,0)^*!U{\cdots}
 \ar@{-} (30,0);(35,0) *+<12pt>!D{\alpha_{n-2}} *++!U{{}_2} *\cir<3pt>{}="H"
 \ar@{-} "H";(45,0) *+<12pt>!D{\alpha_{n-1}} *++!U{{}_2} *\cir<3pt>{}="I"
 \ar@{<=} "I";(55,0) *++!D{\alpha_{n}} *++!U{{}_1} *\cir<3pt>{}
 \ar@{-} "B";(10,-10) *++!L{\alpha_0} *++!R{-\epsilon_1-\epsilon_2}*\cir<1.5pt>{} *\cir<3pt>{}
\end{xy}
\ \ $\tilde C_2'\simeq\tilde B_2'$\!\!
\begin{xy}
 \ar@{<=}  *++!D{\alpha_0} *\cir<1.5pt>{} *\cir<3pt>{} ;
  (10,0)  *++!D{\alpha_2} *++!U{{}_1} *\cir<3pt>{}="A"
 \ar@{=>} "A";(20,0) *++!D{\alpha_1} *++!U{{}_1} *\cir<3pt>{}="B"
\end{xy}\\
$A_{2n}^{(2)}=\widetilde {BC}_n$\!\!
\begin{xy}
 \ar@{=>}  *++!D{\alpha_0} *++!U{-2\epsilon_1} *\cir<1.5pt>{}*\cir<3pt>{} ;
  (10,0)  *++!D{\alpha_1} *++!U{{}_2} *\cir<3pt>{}="A"
 \ar@{-} "A";(20,0) *++!D{\alpha_2} *++!U{{}_2} *\cir<3pt>{}="B"
 \ar@{-} "B";(25,0) \ar@{.} (25,0);(30,0)^*!U{\cdots}
 \ar@{-} (30,0);(35,0) *+<12pt>!D{\alpha_{n-2}} *++!U{{}_2} *\cir<3pt>{}="H"
 \ar@{-} "H";(45,0) *+<12pt>!D{\alpha_{n-1}} *++!U{{}_2} *\cir<3pt>{}="I"
 \ar@{=>} "I";(55,0) *++!D{\alpha_{n}} *++!U{{}_2} *\cir<3pt>{}
\end{xy}\ \ 
$A_2^{(2)}=\widetilde{BC_1}$\!\!
\begin{xy}
 \ar@<.7mm>@{=}
 *++!D{\alpha_0} *\cir<1.5pt>{} *\cir<3pt>{}="A";
  (10,0) *++!D{\ \ \alpha_1} *++!U{\ \ {}_2} 
  *{\!\text{\large$>$}\!\raisebox{3pt}{\cir<3pt>{}}}="B"
 \ar@<.7mm>@{=} "B";"A"
\end{xy}\\
%%%%% F_4 %%%%%
$E_6^{(2)}=\tilde F_4'$\quad
\begin{xy}
 \ar@{-}  *++!D{\alpha_1} *++!U{{}_1} *\cir<3pt>{};
  (10,0)  *++!D{\alpha_2} *++!U{{}_2} *\cir<3pt>{}="A"
 \ar@{=>} "A";(20,0) *++!D{\alpha_3} *++!U{{}_3} *\cir<3pt>{}="B"
 \ar@{-} "B";(30,0) *++!D{\alpha_4} *++!U{{}_2} *\cir<3pt>{}="C"
 \ar@{-} "C";(40,0) *++!D{\alpha_0} *++!U{-\epsilon_1} 
 *\cir<1.5pt>{} *\cir<3pt>{}
\end{xy}\qquad
$D_4^{(3)}=\tilde G_2'$\!\!
\begin{xy}
 \ar@3{->}  *++!D{\alpha_1}  *++!U{{}_1} *\cir<3pt>{};
  (10,0)  *++!D{\alpha_2} *++!U{{}_2} *\cir<3pt>{}="A"
 \ar@{-} "A";(20,0) *++!D{\alpha_0} *++!U{\epsilon_2-\epsilon_3} *\cir<1.5pt>{} *\cir<3pt>{}
\end{xy}
\fi
\\
Here $\tilde R$ denotes the extended Dynkin diagram of type $R$
with the numbers $m_j(\alpha_{max})$ attached to $\alpha_j\in\tilde\Psi$, 
which has been given in this section.
\end{prop}
\begin{proof}
Fix any vertex $P$ in $\mathcal G$ and retain the notation in 
Definition~\ref{defn:affineDy}.

If $p\ge2$ and $Q_1=Q_2$, then $2m\ge 2m_1$ and 
we similarly have $2m_1\ge 2m$ and
hence $k_1=k_2=1$ and $p=2$.
This only happens when $\mathcal G=\tilde A_1$.

Now we may assume $m_1\ge\dots\ge m_p$ and 
$Q_i\ne Q_j$ if $i\ne j$.

\underline{Claim}:
Let $P_1,\dots,P_\ell$ be vertices in $\mathcal G$ such that for any 
$j=2,\dots\ell-1$,
$P_j$ is linked only to both $P_{j-1}$ and $P_{j+1}$
and no arrow points to $P_j$.
Then the corresponding attached numbers $m_1,\dots,m_\ell$ form an 
arithmetical progression series.

Since $2m_j\ge m$, the relation \eqref{eq:ExtNumber} assures $\sum k_j\le 4$.

Note that if $2m_j=m$, $Q_j$ is an end vertex to which
no arrow points.

\underline{Case $\sum k_j=4$}:
Then $m_1=\cdots=m_p=\frac m2$, 
$Q_j$ are end vertices and there is no arrow starting from $P$.
We may assume $k_1\ge\dots\ge k_p$.
Hence $(k_1,\dots,k_p)=(1,1,1,1)$, $(2,1,1)$, $(2,2)$, $(3,1)$
or $(4)$ and $\mathcal G=\tilde D_4$, $\tilde B_2'$, 
$\tilde B_2$, $\tilde G_2'$ or $\widetilde{BC}_1$, respectively.

\underline{Case $\sum k_j=3$ and $p=1$}:
We have $k_1=3$ and $m=\frac32m_1$.
Then there exists a vertex $Q'\in G\setminus\{P\}$ 
with the number $m'\ge\frac{m_1}2$ such that $Q'$ is linked to $Q$. 
Since $2m_1\ge m+m'$, $m'=\frac 12m_1$ and $Q'$ is
an end vertex.
Hence $\mathcal G=\tilde G_2$.

\underline{Case $\sum k_j=3$, $p\ge2$ and $m_1\ge m$}:
Note that $k_1=1$ and $p=2$ or $3$.
When $p=2$, $k_2=2$ and $(m_1,m_2)=(m,\frac12m)$.
When $p=3$, $k_2=k_3=1$ and $(m_1,m_2,m_3)=(m,\frac12m,\frac12m)$.
Hence if $2\le j\le p$, $Q_j$ is an end vertex as follows.
\[
p=2\quad
\begin{xy}
 \ar@{=>}  *+!D{Q_2} *+!U{\frac m2} *{\circ} ;
  (7,0)  *+!D{P_{\phantom{1}}} *+!U{m} *{\circ}="A"
 \ar@{.} "A";(14,0) *+!D{Q_1} *+!U{m} *{\circ}="B"
 \ar@{.} "B";(20,0) *{}
\end{xy}\qquad\qquad\qquad
p=3\quad
\begin{xy}
 \ar@{-}  *+!D{Q_2} *+!U{\frac m2} *{\circ} ;
  (7,0)  *+!D{P_{\phantom{1}}} *+!L!U{m} *{\circ}="B"
 \ar@{.} "B";(14,0) *+!D{Q_1} *+!U{m} *{\circ}="C"
 \ar@{.} "C";(20,0) *{};
 \ar@{-} "B";(7,-7) *+!R{\frac m2} *++!L{Q_3} *{\circ}
\end{xy}
\]
Denoting $P_1=P$ and $P_2=Q_1$, we choose the maximal sequence of vertices
$P_1,\dots,P_\ell$ given in the above claim.
The numbers attached to these vertices are $m$.
If an arrow links $P_{\ell-1}$ to $P_\ell$, 
it has double lines and points toward the end vertex $P_\ell$
and hence $\mathcal G=\widetilde{BC}_n$ or $\tilde B_n$ 
according to $p=2$ or $3$, respectively.

Now we may assume that $P_j$ is linked  to $P_{j+1}$ by a line
if $1\le j< \ell$. 
Then $P_\ell$ is not an end vertex and therefore $P_\ell$ is a branching 
vertex or there exists an arrow pointing toward $P_\ell$.
Applying the argument we have done to the vertex $P_\ell$ in place of $P$,
we conclude the following.
If $P_\ell$ is a branching vertex, $P_\ell$ is linked to two end vertices 
together with $P_{\ell-1}$ by lines and we have $\tilde C_n'$ or
$\tilde D_n$ according to $p=2$ or $3$, respectively.
If $P_\ell$ is not a branching vertex, an arrow starting from an end 
vertex points toward $P_\ell$ and we have accordingly 
$\mathcal G=\tilde C_n$ or $\tilde C_n'$. 

\underline{Case $\sum k_j=3$, $p\ge2$ and $m_1<m$}:
Fix $j$ with $1\le j\le p$ and let $P_1,\dots,P_{\ell_j}$ be 
the maximal sequence given in the claim such that $P_1=P$ and $P_2=Q_j$.
The corresponding attached numbers $m'_1=m,m'_2=m_j,m'_3,\dots,m'_{\ell_j}$ 
form a strictly decreasing arithmetical progression series and
the argument in the preceding case assures that $P_{\ell_j}$ is not
a branching vertex. Moreover $P_{\ell_j}$ is not linked to any arrow 
but it is an end vertex. Therefore
$m_j=\frac{\ell_j-1}{\ell_j}m$.  
If $p=3$, \eqref{eq:ExtNumber} implies
\begin{equation}
  \tfrac1{\ell_1}+\tfrac1{\ell_2}+\tfrac1{\ell_3}=1,
  \quad \ell_1\ge \ell_2\ge \ell_3
\end{equation}
and hence $(\ell_1,\ell_2,\ell_3)=(3,3,3)$, $(4,4,2)$ or $(6,3,2)$ and 
$\mathcal G=\tilde E_6$, $\tilde E_7$ or $\tilde E_8$,
respectively.
Similarly if $p=2$, we have $k_1=2$, $k_2=1$ and
$(\ell_1,\ell_2)=(3,3)$ or $(4,2)$ and  
$\mathcal G=\tilde F_4'$ or $\tilde F_4$, respectively.

\underline{Other cases}:
Now we may assume that $\mathcal G$ has no branching vertex and moreover that
if $\mathcal G$ contains an arrow, the arrow has double lines and points 
toward an end vertex.
Hence it follows from the claim that $\mathcal G$ equals
$\tilde B_n'$ if $\mathcal G$ contains an arrow and $\tilde A_n$ 
($n\ge 2$) if otherwise.
\end{proof}
\begin{rem}\label{rem:cls}
Retain the notation in Proposition~\ref{prop:clasify}.

{\rm i)} 
$\tilde A_1$ 
% \begin{xy}
%  \ar@{=} *\cir<1.5pt>{} *\cir<3pt>{};(8,0) *++!D{{}_1} *\cir<3pt>{}
% \end{xy}
is sometimes denoted by
 \begin{xy}
  \ar@{<=>} *++!D{{}_1} *\cir<1.5pt>{} *\cir<3pt>{};(8,0) *++!D{{}_1} *\cir<3pt>{}
 \end{xy}
or
 \begin{xy}
  \ar@<.7mm>@{=} *++!D{{}_1} *\cir<1.5pt>{} *\cir<3pt>{}="A";
    (8,0) *++!D_{{}_1}*\cir<3pt>{}="B"
  \ar@<-.7mm>@{=} "A";"B" \end{xy}
or
 \begin{xy}
  \ar@{-}^*{{}_\infty} *++!D{{}_1} *\cir<1.5pt>{} *\cir<3pt>{};(8,0) 
  *++!D{{}_1} *\cir<3pt>{}
 \end{xy}.

{\rm ii)} The proposition is known as the classification of the generalized
Cartan matrices of affine type (cf.~\cite[Ch.~4]{Ka}), where $\tilde R$,
$\tilde B_n'$, $\tilde C_n'$, $\widetilde{BC}_n$, $\tilde F_4'$ and 
$\tilde G_2'$ are 
denoted by $R^{(1)}$, $D_{n+1}^{(2)}$, $A_{2n-1}^{(2)}$, $A_{2n}^{(2)}$, 
$E_6^{(2)}$ and $D_4^{(3)}$, respectively.

{\rm iii)} Suppose $R$ is not simply laced and let $\alpha'_{\max}$ be the
maximal root in $\Sigma\setminus\Sigma^L$ and 
put $\alpha'_0=-\alpha_{max}'$.
Then the root $\alpha_0'$ corresponds to the vertex in $\tilde R'$ 
indicated by $\ccirc$ and we get the Dynkin diagram of type $R$ 
by deleting the vertex.
%Under the realization of the root systems and the fundamental systems 
%in this section 
%$\alpha_0'=-\epsilon_1$, $-\epsilon_1-\epsilon_2$, $-\epsilon_1$ or
%$\epsilon_2-\epsilon_3$ when $R=B_n$, $C_n$, $F_4$ or $G_2$, respectively.

If we change the arrows in $\tilde R'$ by those with the opposite directions,
$\tilde R'$ is changed into the extended Dynkin diagram of the dual root 
system of $R$. 

The root system
\begin{equation}
\{\pm(\epsilon_i-\epsilon_j),\ \pm(\epsilon_i+\epsilon_j),\ 
\pm\epsilon_k,\ \pm2\epsilon_k\,;\,1\le i<j\le n,\ 1\le k\le n\}
\end{equation}
in $\mathbb R^n$ is the non-reduced root system of type $BC_n$.
If we denote the maximal root of the non-reduced root system of 
type $BC_n$ by $\alpha_{max}$, the root $\alpha_0:=-\alpha_{max}$ 
corresponds to the vertex in $\widetilde{BC}_n$ indicated by $\ccirc$.

Note that for an irreducible root system $\Sigma$ with a fundamental system
$\Psi$, the set 
$\{\beta\in\Sigma\,;\,(\alpha|\beta)\le0\ (\forall\alpha\in\Psi)\}$
is the complete representatives of the orbits under the action of its Weyl 
group or equivalently the decomposition of $\Sigma$ according to the length of
the roots.
The attached numbers in the above diagrams are the coefficients 
$m_j(-\alpha_0)$ in the expression 
$-\alpha_0=\sum_{\alpha_j\in\Psi}m_j(-\alpha_0)\alpha_j$
for the element $\alpha_0$ in the set.
Here we don't assume $\Sigma$ is reduced.

{\rm iv)}
We easily get a realization of $\mathcal G$ in an Euclidean space as
follows.  For example, if $\mathcal G=\tilde E_8$, we first realize
$\langle\alpha_2,\dots,\alpha_8,\alpha_0\rangle\simeq D_8$ 
as in the standard way given in this section and
then $\alpha_1$ is determined by \eqref{eq:zerosum} and moreover
$E_7=\{\alpha_0\}^\perp$ and $E_7=\{\alpha_0,\alpha_8\}^\perp$.
If $\mathcal G=\tilde F_4$, we first realize 
$\langle\alpha_0,\alpha_1,\alpha_2,\alpha_3\rangle\simeq B_4$ and then 
$\alpha_4$.

{\rm v)}
The connected diagram $\mathcal G$ corresponds to an indecomposable finite 
subset $\Phi$ of $\mathbb R^n\setminus\{0\}$ with 
$\#\Phi=n+1$ such that
\begin{equation}
 \alpha\in\Phi,\ \beta\in\Phi\text{ and }\alpha\ne\beta
 \Rightarrow -2\frac{(\alpha|\beta)}{(\alpha|\alpha)}\in\{0,1,2,\ldots\}.
\end{equation}
Here the subset $\Phi$ of $\mathbb R^n$
is indecomposable if there exists no non-trivial 
decomposition $\Phi=\Phi_1\cup\Phi_2$ with $\Phi_1\perp\Phi_2$.

{\rm vi)}
The diagram $\mathcal G$ with arrows is constructed as a quotient under
the action of an automorphism of an extended Dynkin diagram of a simply laced
root system.
The arrow between two vertices in the quotient represents the difference 
of the numbers of the corresponding original vertices as follows.
%which include all connected affine Dynkin diagrams.
\begin{gather*}
 \tilde A_{2n-1}\Rightarrow\tilde B_{n\ (n\ge2)}'
 :
 \begin{xy}
 \ar@{-} *+!D{{}_1} *{\circ}="A"; (5,1) *+!D{{}_1} *{\circ}="B"
 \ar@{-} "B";(7.5,1) \ar@{.} (7.5,1);(10,1)^*!U{{}_{\cdots}}
 \ar@{-} (10,1);(12.5,1) *+!D{{}_1} *{\circ}="C"
 \ar@{-} "C"; (17.5,0) *+!D{{}_1} *{\circ}="D";
 \ar@{-} "A"; (5,-1) *{\circ}="E"
 \ar@{-} "E";(7.5,-1) \ar@{.} (7.5,-1);(10,-1)
 \ar@{-} (10,-1);(12.5,-1) *{\circ}="F"
 \ar@{-} "F"; "D" *{\circ};
\end{xy}
\to
\begin{xy}
 \ar@{<=} *+!D{{}_1} *{\circ}="A"; (5,0) *+!D{{}_1} *{\circ}="B"
 \ar@{-} "B";(7.5,0) \ar@{.} (7.5,0);(10,0)^*!U{{}_{\cdots}}
 \ar@{-} (10,0);(12.5,0) *+!D{{}_1} *{\circ}="C"
 \ar@{=>} "C"; (17.5,0) *+!D{{}_1} *{\circ}
\end{xy}
\allowdisplaybreaks\\
%%%%%%%%%%%%%%%%%%%%
 \tilde D_{2n}\Rightarrow\tilde B_n
 \Rightarrow\widetilde{BC}_{n-1\ (n\ge 3)}
  :
\begin{xy}
 \ar@{-} (0,-2) *{\circ}="A"; (5,1) *+!D{{}_2} *{\circ}="B"
 \ar@{-} (0,3) *+!D{{}_1} *{\circ} ; "B"
 \ar@{-} "B";(7.5,1) \ar@{.} (7.5,1);(10,1)^*!U{{}_{\cdots}}
 \ar@{-} (10,1);(12.5,1) *+!D{{}_2} *{\circ}="C"
 \ar@{-} "C"; (17.5,0) *+!D{{}_2} *{\circ}="D";
 \ar@{-} (0,-4) *{\circ} ; (5,-1) *{\circ}="E"
 \ar@{-} (0,1) *{\circ} ; "E"
 \ar@{-} "E";(7.5,-1) \ar@{.} (7.5,-1);(10,-1)
 \ar@{-} (10,-1);(12.5,-1) *{\circ}="F"
 \ar@{-} "F"; "D" *{\circ};
\end{xy}
\to
\begin{xy}
 \ar@{-} (0,-1) *{\circ}="A"; (5,0) *+!D{{}_2} *{\circ}="B"
 \ar@{-} (0,1) *+!D{{}_1} *{\circ} ; "B"
 \ar@{-} "B";(7.5,0) \ar@{.} (7.5,0);(10,0)^*!U{{}_{\cdots}}
 \ar@{-} (10,0);(12.5,0) *+!D{{}_2} *{\circ}="C"
 \ar@{=>} "C"; (17.5,0) *+!D{{}_2} *{\circ}="D"
\end{xy}
\to
\begin{xy}
 \ar@{=>} *+!D{{}_1} *{\circ}="A"; (5,0) *+!D{{}_2} *{\circ}="B"
 \ar@{-} "B";(7.5,0) \ar@{.} (7.5,0);(10,0)^*!U{{}_{\cdots}}
 \ar@{-} (10,0);(12.5,0) *+!D{{}_2} *{\circ}="C"
 \ar@{=>} "C"; (17.5,0) *+!D{{}_2} *{\circ}="D"
\end{xy}
\allowdisplaybreaks\\
%%%%%%%%%%%%%%%%%
 \tilde D_{n+1}\Rightarrow\tilde C_n'\Rightarrow
 \tilde C_{n-1\ (n\ge 3)}
  :
\begin{xy}
 \ar@{-} (0,-2) *{\circ}="A"; (5,0) *+!D{{}_2} *{\circ}="B"
 \ar@{-} (0,2) *+!D{{}_1} *{\circ} ; "B"
 \ar@{-} "B";(7.5,0) \ar@{.} (7.5,0);(10,0)^*!U{{}_{\cdots}}
 \ar@{-} (10,0);(12.5,0) *+!D{{}_2} *{\circ}="D"
 \ar@{-} "D";(17.5,1) *+!D{{}_1} *{\circ}
 \ar@{-} "D";(17.5,-1) *{\circ}
\end{xy}
\to
\begin{xy}
 \ar@{-} (0,-1) *{\circ}="A"; (5,0) *+!D{{}_2} *{\circ}="B"
 \ar@{-} (0,1) *+!D{{}_1} *{\circ} ; "B"
 \ar@{-} "B";(7.5,0) \ar@{.} (7.5,0);(10,0)^*!U{{}_{\cdots}}
 \ar@{-} (10,0);(12.5,0) *+!D{{}_2} *{\circ}="D"
 \ar@{<=} "D"; (17.5,0) *+!D{{}_1} *{\circ}
\end{xy}
\to
\begin{xy}
 \ar@{=>} *+!D{{}_1} *{\circ}="A"; (5,0) *+!D{{}_2} *{\circ}="B"
 \ar@{-} "B";(7.5,0) \ar@{.} (7.5,0);(10,0)^*!U{{}_{\cdots}}
 \ar@{-} (10,0);(12.5,0) *+!D{{}_2} *{\circ}="D"
 \ar@{<=} "D"; (17.5,0) *+!D{{}_1} *{\circ}
\end{xy}
\allowdisplaybreaks\\
%%%%%%%%%%%%%%%%%%%%
\begin{xy}\ar@3{->} *+{\tilde D_4};   (10,0) *+{\tilde G_2'}\end{xy}\!\!
 :
\begin{xy}
 \ar@{-} *{\circ}="A" ; 
         (5,0) *{\circ}="B"
 \ar@{-} "B";(10,0) *+!D{{}_1} *{\circ}="C"
 \ar@{-} (0,-2) *{\circ} ; "B" *{\circ}
 \ar@{-} (0,2) *+!D{{}_1} *{\circ} ; "B" *+!D{{}_2} *{\circ}
\end{xy}
\to
\begin{xy}
 \ar@3{->}  *+!D{{}_1} *{\circ} ;
  (6,0)   *+!D{{}_2} *{\circ}="C"
 \ar@{-} "C";(11,0)  *+!D{{}_1} *{\circ}
\end{xy}
\qquad\qquad\qquad\qquad
\begin{xy}
  \ar@<.7mm>@{=} *+{\tilde D_4}="A";
    (11,0) *+{\!\!\!\text{\Large$\rangle$} \widetilde{BC}_1}="B"
  \ar@<-.7mm>@{=} "A";"B" \end{xy}\!\!
 :
\begin{xy}
 \ar@{-} (0,3)  *+!D{{}_1} *{\circ} ; (5,0) *+!D{{}_2} *{\circ}="B"
 \ar@{-} (0,1)  *{\circ};"B" *{\circ}
 \ar@{-} (0,-1) *{\circ};"B" *{\circ}
 \ar@{-} (0,-3) *{\circ};"B" *{\circ}
\end{xy}
\to
\begin{xy}
 \ar@<.7mm>@{=} *+!D{{}_1} *{\ccirc}="A";
  (5,0) *+!D{\;{}_2} *{\!\text{\Large$\rangle$}\!\circ}="B"
 \ar@<.7mm>@{=} "B";"A"
\end{xy}
\allowdisplaybreaks\\
 \tilde E_6\Rightarrow\tilde F_4'
 :
\begin{xy}
 \ar@{-} (0,1) *+!D{{}_1} *{\circ}; (5,1) *+!D{{}_2} *{\circ}="B"
 \ar@{-} "B";(10,0) *+!D{{}_3} *{\circ}="C"
 \ar@{-} (0,-1) *{\circ}; (5,-1) *{\circ}="D"
 \ar@{-} "D";"C"
 \ar@{-} "C";(15,0) *+!D{{}_2} *{\circ}="E"
 \ar@{-} "E";(20,0) *+!D{{}_1} *{\circ}
\end{xy}
 \to
\begin{xy}
 \ar@{-} *+!D{{}_1} *{\circ}; (5,0) *+!D{{}_2} *{\circ}="B"
 \ar@{=>} "B";(10,0) *+!D{{}_3} *{\circ}="C"
 \ar@{-} "C";(15,0) *+!D{{}_2} *{\circ}="E"
 \ar@{-} "E";(20,0) *+!D{{}_1} *{\circ}
\end{xy}
%\allowdisplaybreaks\\
\ \ \quad
\begin{xy}\ar@3{->} *+{\tilde E_6};   (10,0) *+{\tilde G_2}\end{xy}
\!\!
 :
\begin{xy}
 \ar@{-} *{\circ}; (5,0) *{\circ}="B"
 \ar@{-} "B";(10,0) *{\circ}="C"
 \ar@{-} (0,-2) *{\circ} ;(5,-2) *{\circ}="X",
 \ar@{-} (0,2) *+!D{{}_1} *{\circ} ;(5,2) *+!D{{}_2} *{\circ}="Y"
 \ar@{-} "X";"C" *{\circ};
 \ar@{-} "Y";"C" *+!D{{}_3} *{\circ}
\end{xy}\ \to
\begin{xy}
 \ar@{-} *+!D{{}_1} *{\circ}="A" ; 
         (5,0) *+!D{{}_2} *{\circ}="C"
 \ar@3{->} "C";(10,0) *+!D{{}_3} *{\circ}
\end{xy}
\allowdisplaybreaks\\
 \tilde E_7\Rightarrow\tilde F_4
  :
\begin{xy}
 \ar@{-} (0,1) *+!D{{}_1} *{\circ}; (5,1) *+!D{{}_2} *{\circ}="B"
 \ar@{-} "B";(10,1) *+!D{{}_3} *{\circ}="C"
 \ar@{-} "C";(15,0) *+!D{{}_4} *{\circ}="E"
 \ar@{-} (0,-1) *{\circ}; (5,-1) *{\circ}="D"
 \ar@{-} "D";(10,-1) *{\circ}="F"
 \ar@{-} "F";"E"
 \ar@{-} "E";(20,0) *+!D{{}_2} *{\circ}
\end{xy}
\to
\begin{xy}
 \ar@{-} *+!D{{}_1} *{\circ}; (5,0) *+!D{{}_2} *{\circ}="B"
 \ar@{-} "B";(10,0) *+!D{{}_3} *{\circ}="C"
 \ar@{=>} "C";(15,0) *+!D{{}_4} *{\circ}="E"
 \ar@{-} "E";(20,0) *+!D{{}_2} *{\circ}
\end{xy}
\qquad\quad
\tilde E_8:
\begin{xy}
 \ar@{-} *+!D{{}_2} *{\circ}; (5,0) *+!D{{}_4} *{\circ}="B"
 \ar@{-} "B";(10,0) *+!D{{}_6} *{\circ}="C"
 \ar@{-} "C";(15,0) *+!D{{}_5} *{\circ}="E"
 \ar@{-} "E";(20,0) *+!D{{}_4} *{\circ}="F"
 \ar@{-} "F";(25,0) *+!D{{}_3} *{\circ}="G"
 \ar@{-} "G";(30,0) *+!D{{}_2} *{\circ}="H"
 \ar@{-} "H";(35,0) *+!D{{}_1} *{\circ}
 \ar@{-} "C";(10,-5) *+!L{{}_3} *{\circ}
\end{xy}
\begin{xy}
\end{xy}
\end{gather*}

{\rm vii)}
Allowing a line linking a vertex to the same vertex in $\mathcal G$, 
we have the following extra ones.
\begin{equation}
\begin{gathered}
\begin{xy}
 \ar@{-}  *+!D{{}_1} *{\circ} ; (5,0)  *+!D{{}_2} *{\circ}="B"
 \ar@{-} "B";(10,0) *+!D{{}_2} *{\circ}="C"
 \ar@{-} "C";(12.5,0) \ar@{.} (12.5,0);(15,0)^*!U{{}_{\cdots}}
 \ar@{-} (15,0);(17.5,0) *+!D{{}_2} *{\circ}="H"
 \ar@{-} "H";(22.5,0) *+!D{{}_2} *{\circ}="I"
 \ar@{-} @(ur,dr) "I";"I"
 \ar@{-} "B";(5,-5) *+!L{{}_1} *{\circ}
\end{xy}
\quad
\begin{xy}
 \ar@{<=}  *+!D{{}_1} *{\circ} ; (5,0)  *+!D{{}_1} *{\circ}="B"
 \ar@{-} "B";(10,0) *+!D{{}_1} *{\circ}="C"
 \ar@{-} "C";(12.5,0) \ar@{.} (12.5,0);(15,0)^*!U{{}_{\cdots}}
 \ar@{-} (15,0);(17.5,0) *+!D{{}_1} *{\circ}="H"
 \ar@{-} "H";(22.5,0) *+!D{{}_1} *{\circ}="I"
 \ar@{-} @(ur,dr) "I";"I"
\end{xy}
\quad
\begin{xy}
 \ar@{=>}  *+!D{{}_1} *{\circ} ; (5,0)  *+!D{{}_2} *{\circ}="B"
 \ar@{-} "B";(10,0) *+!D{{}_2} *{\circ}="C"
 \ar@{-} "C";(12.5,0) \ar@{.} (12.5,0);(15,0)^*!U{{}_{\cdots}}
 \ar@{-} (15,0);(17.5,0) *+!D{{}_2} *{\circ}="H"
 \ar@{-} "H";(22.5,0) *+!D{{}_2} *{\circ}="I"
 \ar@{-} @(ur,dr) "I";"I"
\end{xy}\\[-12pt]
\begin{xy}
 \ar@{-} (5,0) *+!D{{}_1} *{\circ}="B";(10,0) *+!D{{}_1} *{\circ}="C"
 \ar@{-} "C";(12.5,0) \ar@{.} (12.5,0);(15,0)^*!U{{}_{\cdots}}
 \ar@{-} (15,0);(17.5,0) *+!D{{}_1} *{\circ}="H"
 \ar@{-} "H";(22.5,0) *+!D{{}_1} *{\circ}="I"
 \ar@{-} @(ur,dr) "I";"I"
 \ar@{-} @(ul,dl) "B";"B"
\end{xy}
\qquad
\begin{xy}
 \ar@{-} *+!D{{}_1} *{\circ}="A"
 \ar@{-} @(ur,dr) "A";"A"
 \ar@{-} @(ul,dl) "A";"A"
\end{xy}
\end{gathered}
\end{equation}

{\rm viii)}
If the assumption of the finiteness of the vertices is dropped
in the proposition, we moreover have the following $\mathcal G$, which
easily follows from the proof of the proposition.
\begin{equation}\begin{split}
 A_{+\infty}\ 
  &\begin{xy}
   \ar@{-} *+!D{{}_1} *{\circ} ; (5,0) *+!D{{}_2} *{\circ}="B"
   \ar@{-} "B"; (10,0) *+!D{{}_3} *{\circ}="C"
   \ar@{-} "C";(13,0) \ar@{.} (13,0);(16,0)
  \end{xy}\qquad
 A_\infty\ \ 
  \begin{xy}
   \ar@{.} ;(3,0) \ar@{-} (3,0);(6,0) *+!D{{}_1} *{\circ}="A"
   \ar@{-} "A";(11,0) *+!D{{}_1} *{\circ}="B"
   \ar@{-} "B";(14,0) \ar@{.} (14,0);(17,0)
  \end{xy}\qquad
 D_\infty\ 
  \begin{xy}
   \ar@{-} *+!D{{}_1} *{\circ} ; (5,0) *+!D{{}_2} *{\circ}="B"
   \ar@{-} "B"; (10,0) *+!D{{}_2} *{\circ}="C"
   \ar@{-} "C";(13,0) \ar@{.} (13,0);(16,0)
   \ar@{-}  "B";(5,-5) *+!L{{}_1} *{\circ}
  \end{xy}\\[-6pt]
 B_\infty\ \ 
  &\begin{xy}
   \ar@{<=} *+!D{{}_1} *{\circ} ; (5,0) *+!D{{}_1} *{\circ}="B"
   \ar@{-} "B"; (10,0) *+!D{{}_1} *{\circ}="C"
   \ar@{-} "C";(13,0) \ar@{.} (13,0);(16,0)
  \end{xy}\qquad
 C_\infty\ 
  \begin{xy}
   \ar@{=>} *+!D{{}_1} *{\circ} ; (5,0) *+!D{{}_2} *{\circ}="B"
   \ar@{-} "B"; (10,0) *+!D{{}_2} *{\circ}="C"
   \ar@{-} "C";(13,0) \ar@{.} (13,0);(16,0)
  \end{xy}
  \qquad
 \begin{xy}
 \ar@{-} *+!D{{}_1} *{\circ}="A" ; (5,0) *+!D{{}_1} *{\circ}="B"
 \ar@{-} "B" ;(10,0) *+!D{{}_1} *{\circ}="C"
 \ar@{-} "C";(12.5,0) \ar@{.} (12.5,0);(15,0)^*!U{{}_{\cdots}}
 \ar@{-} @(ul,dl) "A";"A"
 \end{xy}
\end{split}\end{equation}
\end{rem}
%%%%%%%%%%%%%%%%%% TABLE %%%%%%%%%
\section{Tables}\label{sec:table}
In this section we assume that $\Sigma$ is an irreducible and reduced 
root system. We will classify the elements of $\Hom(\Xi,\Sigma)$ under 
a suitable isomorphisms for every root system $\Xi$.
%%%
\begin{defn}
For $\iota$, $\iota'\in\Hom(\Xi,\Sigma)$ we define that
$\iota$ is \textsl{weakly equivalent} to $\iota'$
if and only if there exists $g\in\Aut(\Sigma)=\Hom(\Sigma,\Sigma)$ with
$
 \iota'(\Xi) = g\circ\iota(\Xi),
$
that it, $\iota'(\Xi)\underset{\Sigma}{\overset{w}{\sim}}\Xi$.
Then $\Hom(\Xi,\Sigma)$ is decomposed into the equivalence classes.
\end{defn}
%%%
In many cases $\Hom(\Xi,\Sigma)$ itself is the equivalence class
but if it is not so, we will identify every equivalence class
$\Hom(\Xi,\Sigma)_o$ contained in $\Hom(\Xi,\Sigma)$ by a suitable 
geometric condition.

In the tables in this section we will list up all $\Xi$ with 
$\Hom(\Xi,\Sigma)\ne\emptyset$ and classify them with some data
under the following notation.
\begin{align*}
\Aut(\Xi) &:= \Hom(\Xi,\Xi),\qquad
\Aut(\Sigma) := \Hom(\Sigma,\Sigma),\allowdisplaybreaks\\
\Aut'(\Xi)&:=\prod_{j=1}^m\Aut(\Xi_j)\subset\Aut(\Xi)\text{ for the 
irreducible decomposition }\\[-5pt]
 &\qquad\qquad
\Xi = \Xi_1+\cdots+\Xi_m,
  \allowdisplaybreaks\\
\#\,&:\,\#\bigl(W_{\Sigma}\backslash\!\Hom(\Xi,\Sigma)_o\bigr),
\\
\#_{\Xi}\,&:\,\#\bigl(W_{\Sigma}\backslash\!\Hom(\Xi,\Sigma)_o/\!\Aut(\Xi)\bigr),\\
\#_{\Xi'}\,&:\,\#\bigl(W_{\Sigma}\backslash\!\Hom(\Xi,\Sigma)_o/\!\Aut'(\Xi)\bigr),\\
\#_{\Sigma}\,&:\,\#\bigl(\Aut(\Sigma)\backslash\!\Hom(\Xi,\Sigma)_o\bigr),
 \allowdisplaybreaks\\
\Xi^{\perp\perp}\,&:\,
 \begin{cases}
   \circ&\bigl((\Xi^\perp)^\perp=\Xi\text{ and \eqref{eq:special} is valid}\bigr),\\
   \times&\bigl((\Xi^\perp)^\perp=\Xi\text{ but \eqref{eq:special} is not valid}\bigr),\\
   (\Xi^\perp)^\perp&\bigl((\Xi^\perp)^\perp\ne\Xi\bigr),
 \end{cases}
 \allowdisplaybreaks\\
P\,&:\,
 \begin{cases}
   \#\{\Theta\subset\Psi\,;\,\langle\Theta\rangle\overset{w}{\underset{\Sigma}{\sim}}\Xi\}
    \qquad\bigl(\text{if }\Xi\text{ is fundamental}\bigr),\\
%    \circ\ \ \qquad\bigl(\text{if $\Xi$ is not fundamental but $L$-closed}\bigr),\\
    \leftarrow\qquad\bigl(\text{$L$-closure of $\Xi$ is given by $(\Xi^\perp)^\perp$}\bigr),\\
    \to\qquad\bigl(\text{$L$-closure of $\Xi$ is given in the right column}\bigr),
 \end{cases}
 \allowdisplaybreaks\\
L\,&:\,
  L\text{-closure}
     \qquad\bigl(\text{if }\rank(\Xi^\perp)^\perp>\rank\Xi\text{ and }\Xi
    \text{ is not $L$-closed}\bigr),
 \allowdisplaybreaks\\
S\,&:\, S\text{-closure}
     \qquad\bigl(\text{if }\Xi\text{ is not $S$-closed (cf.\ Definition~\ref{defn:S-closed})}\bigr),
 \allowdisplaybreaks\\
 \langle j_1,\ldots,j_m\rangle&\,:\,
   \langle\alpha_{j_1},\ldots,\alpha_{j_m}\rangle
  \qquad\bigl(\text{under the notation in \S\ref{sec:Dynkin}}\bigr),\\
 \langle\setminus j\rangle&\,:\,
   \langle\Psi\setminus\{\alpha_j\}\rangle,\allowdisplaybreaks
\intertext{For subsystem $\Xi\subset\Sigma$ and a subgroup 
$G$ of $\Aut(\Sigma)$ we put}
%\begin{align*}
  \Hom(\Xi,\Sigma)_o &:= \{\iota\in\Hom(\Xi,\Sigma)\,;\,
   \iota(\Xi)\underset{\Sigma}{\overset{w}\sim}\Xi\},\allowdisplaybreaks\\
%\end{align*}
  O_\Xi^w &:=\{
    \Theta\subset \Sigma\,;\,W_\Theta\Theta=\Theta\text{ and }
    \Theta\underset{\Sigma}{\overset{w}\sim}
    \Xi\}\qquad\text{(cf.~Definition~\ref{def:weq})},\\
  O_\Xi  &:=\{
      \Theta\subset \Sigma\,;\,W_\Theta\Theta=\Theta\text{ and }
    \Theta\underset{\Sigma}{\sim}
    \Xi\},\allowdisplaybreaks\\
  N_G(\Xi) &:= \{g\in G\,;\, g(\Xi)=\Xi\}. %,\\
%  Z_G(\Xi) &:= \{g\in G\,;\, g|_\Xi = id\}.
\end{align*}
Then
\begin{align}
 O_\Xi^w &\simeq
 \Hom(\Xi,\Sigma)_o/\!\Aut(\Xi) \simeq
 \Aut(\Sigma)/N_{\Aut(\Sigma)}(\Xi),\notag\\
 O_\Xi &\simeq W_\Sigma/N_{W_\Sigma}(\Xi),\notag\allowdisplaybreaks\\
 \# O_\Xi^w/\# O_\Xi
  &= \#\bigl(W_\Xi\backslash\! \Hom(\Xi,\Sigma)_o/\!\Aut(\Xi)\bigr)
  = (\#_\Xi),\allowdisplaybreaks\\
 (\#)/(\#_\Xi)
  &=
  \#\bigl(W_\Sigma\backslash\!\Hom(\Xi,\Sigma)_o\bigr)\bigm/
  \#\bigl(W_\Sigma\backslash\!\Hom(\Xi,\Sigma)_o/\!\Aut(\Xi)\bigr)
  \label{eq:NOut} \allowdisplaybreaks\\
  &=\#\bigl(\Out(\Xi)/\!\Out_\Sigma(\Xi)\bigr)\notag \allowdisplaybreaks\\
  &=\#\Bigl(\Out(\Xi)\bigm/\bigl(N_{W_\Sigma}(\Xi)
   /(W_\Xi\times W_{\Xi^\perp})\bigr)
    \Bigr),\notag
  \allowdisplaybreaks\\
 \# O_\Xi 
  &= (\#)\cdot\# W_{\Sigma}\bigm/
     \bigl((\#_\Sigma)\cdot\#\Out(\Xi)\cdot \# W_\Xi\cdot \# W_{\Xi^\perp}\bigr).\label{eq:Nisom}
\end{align}
Here $(\#)$, $(\#_\Xi)$, $(\#_{\Xi'})$ and $(\#_{\Sigma})$
are numbers given in the columns indicated by $\#$, $\#_\Xi$, $\#_{\Xi'}$ 
and $\#_\Sigma$  in the table below, respectively.
Since $1\le(\#_\Xi)\le(\#_{\Xi'})$, $(\#_\Xi)$ may not
be written if $(\#_{\Xi'})=1$.
If $\Out(\Sigma)$ is trivial, $(\#_{\Sigma})=(\#)$ and therefore
$(\#_{\Sigma})$ may not be written.

Note that \eqref{eq:NW} corresponds to the special case of \eqref{eq:Nisom}
with $\Xi=\{\pm\alpha_0\}$.
\begin{rem}\label{rem:ans}
We obtain the answers to the questions in the introduction from the 
table in this section as follows.

Answers to Q1 and Q2 are given by the table.

The number in Q3 is given by \eqref{eq:Nisom} with the table.

The answer to Q4 is yes if and only if $(\#_{\Xi})=(\#)$
(cf.~Remark~\ref{rem:orth} iii), Remark~\ref{rem:last} iii) and \S\ref{sec:OUT}).

The answer to Q5 is yes if and only if $(\#_{\Xi'})=(\#)$.

The answer to Q6 is yes if and only if $(\#_{\Xi})=1$ 
(cf.~Remark~\ref{rem:T1} ii)).

The answer to the first question of Q7 is given by Remark~\ref{rem:fund}
and the number in Q7 is obtained from the column $P$ in the table.
\end{rem}

%\newpage
\subsection{Classical type} 
($\Sigma$ : $A_n$, $B_n$, $C_n$, $D_n$)

\centerline{\textbf{$\Xi$ : Irreducible}}

%\medskip
\noindent
\begin{longtable}{|l|l|c|c|c|l|l|c|} \hline
$\Sigma$ & $\Xi$ & $\#$ & $\#_{\Xi}$ & $\#_{\Sigma}$
  & $\Xi^\perp$ & $\Xi^{\perp\perp}$ & P \\ \hline
  $A_n$ & $A_1$ & 1 & 1 & 1 & $A_{n-2}$ & $\times$ & $n$ \\ \hline
  $A_{n\ (1<m\le n-2)}$ & $A_m$ & 2 & 1 & 1 & $A_{n-m-1}$ &
  $\times$ & $n-m+1$ \\ \hline
  $A_{n\ (n\ge 3)}$ & $A_{n-1}$ & 2 & 1 & 1 & $\emptyset$ &
  $\Sigma$ &$2$ \\ \hline
  $A_{n\ (n\ge2)}$ & $A_n$ & 2 & 1 & 1 & $\emptyset$ &
  $\Sigma$ &$1$\\ \hline
\hline
 $\Sigma\ (n\ge 5)$ & $\Xi$ & $\#$ & $\#_{\Xi}$ & $\#_{\Sigma}$
  & $\Xi^\perp$ & $\Xi^{\perp\perp}$ & P\\ \hline
$D_n$ & $A_1$ &1 &1 &1 &$D_{n-2}+A_1$ & $\times$ &$n$ \\ \hline
$D_n$ & $A_2$ &1 &1 &1 &$D_{n-3}$ & $A_3$&$n-1$ \\ \hline
$D_n$ & $A_3$ &1 &1 &1 &$D_{n-4}$ & $D_4$& $n-2$ \\ \cline{3-8}
    &   $(D_3)$ &1 &1 &1 &$D_{n-3\ (n\ne 7)}$ & $\circ$  &1 \\ 
  \cline{6-7}
& & & & & $D_{4\ (n=7)}$ & $\times$ &  \\ \hline
$D_{n\ (4\le k\le n-3)}$ & $A_k$ &1 &1 &1 &$D_{n-k-1}$ & $D_{k+1}$&
 $n-k+1$
  \\ \hline
$D_n$ & $A_{n-2}$ &1 &1 &1 &$\emptyset$ & $\Sigma$&$3$
  \\ \hline
$D_{n\ (n:\text{odd})}$ & $A_{n-1}$ &2 &1 &1 &$\emptyset$ & $\Sigma$&$2$ 
  \\ \cline{1-1} \cline{3-5}
$D_{n\ (n:\text{even})}$ &  &2 &2 &1 & & & \\ \hline
$D_n$ & $D_4$ &3 &1 &3 &$D_{n-4\ (n\ge 6)}$ & $2$&$1$ \\ \cline{6-7}
 &&&&&$\emptyset_{\ (n=5)}$ &$\Sigma$& \\ \hline
$D_{n\ (4<k\le n-2)}$ & $D_k$ &1 &1 &1 &$D_{n-k}$ & $\circ_{\ \;(k\ne n-4)}$ & $1$\\ \cline{7-7}
& & & & & & $\times_{\ (k=n-4)}$& \\ \hline
$D_{n\ (n\ge6)}$ & $D_{n-1}$ &1 &1 &1 &$\emptyset$ & $\Sigma$ &$1$ \\ \hline
$D_n$ & $D_n$ &2 &1 &1 &$\emptyset$ & $\Sigma$ &$1$ \\ \hline
%\newpage
\hline
%\begin{longtable}{|l|c|c|c|c|c|c|c|} \hline
$\Sigma\ (n\ge2)$ & $\Xi$ & $\#$ & $\#_{\Xi}$ & $\#_{\Sigma}$
  & $\Xi^\perp$ & $\Xi^{\perp\perp}$ & P \\ \hline
$B_n$ & $A_1^L$& 1 & 1 & 1 &  $B_{n-2}+A_1^L$ & $\circ$ & $n-1$ 
\\ \cline{2-8}
      & $A_1^S$ &1 & 1 & 1 &  $B_{n-1}$ & $\circ$ & $1$ 
\\ \hline
$B_{n\ (n\ge3)}$ & $A_2$ & 1 & 1 & 1 &  $B_{n-3}$ & $B_3$ & $n-2$ 
\\ \hline
$B_{n\ (n\ge4)}$ & $A_3$ & 1 & 1 & 1 &  $B_{n-4}$ & $B_4$ & $n-3$ 
\\ \cline{1-1} \cline{3-8}
$B_{n\ (n\ge 3)}$ &$(D_3)$& 1 & 1 & 1 & $B_{n-3}$ &$B_3$ &$\leftarrow$ \\ \hline
$B_{n\ (4\le m<n)}$ & $A_m$ & 1 & 1 & 1 & $B_{n-m-1}$ &$B_{m+1}$ &$n-m$ \\ \hline
$B_{n\ (n\ge4)}$ & $D_4$ & 3 & 1 & 3 & $B_{n-4}$ &$B_4$ &$\leftarrow$ \\ \hline
$B_{n\ (4<m\le n)}$ & $D_m$ & 1 & 1 & 1 & $B_{n-m}$ &$B_m$ &$\leftarrow$ \\ \hline
$B_{n\ (2\le m\le n)}$ & $B_m$ & 1 & 1 & 1 & $B_{n-m}$ &$\circ$ & $1$ \\ 
\hline
%\newpage
\hline
%\end{tabular}
%
%
%\smallskip
%\noindent
%\begin{tabular}{|l|c|c|c|c|c|c|c|l|} \hline
$\Sigma\ (n\ge2)$ & $\Xi$ & $\#$ & $\#_{\Xi}$ & $\#_{\Xi'}$
  & $\Xi^\perp$ & $\Xi^{\perp\perp}$ & P \\ \hline
$C_n$ & $A_1^S$& 1 & 1 & 1 &  $C_{n-2}+A_1^S$ & $\circ$ & $n-1$ 
 \\ \cline{2-8}
      & $A_1^L$ &1 & 1 & 1 &  $C_{n-1}$ & $\circ$ & $1$ 
 \\ \hline
$C_{n\ (n\ge3)}$ & $A_2$ & 1 & 1 & 1 &  $C_{n-3}$ & $C_3$ & $n-2$ 
 \\ \hline
$C_{n\ (n\ge4)}$ & $A_3$ & 1 & 1 & 1 &  $C_{n-4}$ & $C_4$ & $n-3$ 
 \\ \cline{1-1} \cline{3-8}
$C_{n\ (n\ge 3)}$ &$(D_3)$& 1 & 1 & 1 & $C_{n-3}$ &$C_3$ & $\leftarrow\ S:C_3$ \\ \hline
$C_{n\ (4\le m<n)}$ & $A_m$ & 1 & 1 & 1 & $C_{n-m-1}$ &$C_{m+1}$ &$n-m$ \\ \hline
$C_{n\ (n\ge4)}$ & $D_4$ & 3 & 1 & 3 & $C_{n-4}$ &$C_4$ & $\leftarrow\ S:C_4$ \\ \hline
$C_{n\ (4<m\le n)}$ & $D_m$ & 1 & 1 & 1 & $C_{n-m}$ &$C_m$ & $\leftarrow\ S:C_m$ \\ \hline
$C_{n\ (2\le m\le n)}$ & $C_m$ & 1 & 1 & 1 & $C_{n-m}$ &$\circ$ & $1$ \\ \hline
\end{longtable}
The symbol $(D_3)$ is the above table corresponds to $D_3$ in \eqref{eq:BD}.
The subsystems $A_1^L$ and $A_1^S$ of $B_n$ in the above table correspond to
$A_1$ and $B_1$ in \eqref{eq:BD}, respectively.

Applying Remark~\ref{rem:red}~iii) to the table for $\Sigma=B_n$, we
have the table for $\Sigma=C_n$.

Suppose $n>4$.  Then $\#\overline\Hom(A_{n-1},D_n)=2$
and the non-trivial element $g\in\Out(D_n)$ maps its element 
to the other.
Let $A_{n-1}\subset D_n$ by the notation in \S\ref{sec:table}.
Then $h\in\Aut(D_n)$ defined by
$h(\epsilon_j)=-\epsilon_j$ $(j=1,\ldots,n)$ induces the non-trivial
element of $\Out(A_{n-1})$. 
Here $h$ is not an element of $W_{D_n}$ if and only if $n$ is odd.
Hence $\#\bigl(\overline\Hom(A_{n-1},D_n)/\!\Out(A_{n-1})\bigr)=1$
if and only if $n$ is odd.

\bigskip
\centerline{\textbf{$\Sigma=D_4$}}
\nopagebreak

\centerline{
\begin{tabular}{|c|c|c|c|c|c|c|c|c|l|} \hline
$\Sigma$ & $\Xi$ & $\#$ & $\#_{\Xi}$ & $\#_{\Xi'}$ & $\#_{\Sigma}$
  & $\Xi^\perp$ & $\Xi^{\perp\perp}$ & P  \\ \hline
$D_4$ & $A_1$  &1 & 1 & 1 & 1 & $3A_1$ & $\times$ & $4$ \\ \hline
$D_4$ & $A_2$  &1 &1 & 1 & 1 & $\emptyset$ & $\Sigma$ &$3$ \\ \hline
$D_4$ & $A_3$  &3 &3 & 3 & 1 & $\emptyset$ & $\Sigma$ &$3$ \\ \hline
$D_4$ & $D_4$  &6 &1 & 1 & 1 & $\emptyset$ & $\Sigma$ &$1$ \\ \hline
$D_4$ & $2A_1$ &3 &3 & 3 & 1 & $2A_1$ & $\circ$ & $3$ \\ \hline
$D_4$ & $3A_1$ &6 &1 & 6 & 1 & $A_1$ & $\times$ & $1$ \\ \hline
$D_4$ & $4A_1$ &6 &1 & 6 & 1 & $\emptyset$ & $\Sigma$ &$\leftarrow$ \\ \hline
\end{tabular}
}
\bigskip

\centerline{\textbf{$\Sigma$: not of type $D_4$}}
\nopagebreak
We still assume that $\Sigma$ is irreducible and of classical type.
We will examine $\Hom(\Xi,\Sigma)$ when $\Xi$ may not be irreducible.
It is not difficult because the root system and its Weyl group are
easy to describe.
The subsystems of $\Sigma$ can be imbedded in the root space $B_N$
with a sufficiently large $N$.
We should distinguish two subsystems which are isomorphic as root
systems but they are not equivalent by $B_N$.

Under the notation in \S\ref{sec:Dynkin} they are the followings:
\begin{equation}\label{eq:BD}
\begin{split}
   A_1&=%\underset{B_N}{\sim}
    \{\pm(\epsilon_1-\epsilon_2)\},\\
   B_1&=%\underset{B_N}{\sim}
    \{\pm \epsilon_1\}\simeq A_1,\\
   D_2&=%&\underset{B_N}{\sim}
     \langle \epsilon_1-\epsilon_2,\epsilon_1+\epsilon_2\rangle
   \simeq 2A_1,\\
   A_3&=%\underset{B_N}{\sim}
     \langle \epsilon_1-\epsilon_2,\epsilon_2-\epsilon_3,\epsilon_3-\epsilon_4\rangle,\\
   D_3&=%\underset{B_N}{\sim}
     \langle \epsilon_1-\epsilon_2,\epsilon_2-\epsilon_3,\epsilon_2+\epsilon_3\rangle\simeq A_3.
\end{split}
\end{equation}
Let $\{\epsilon_1,\dots,\epsilon_N\}$ be an orthonormal basis
of $\mathbb R^N$ with a sufficiently large positive integer $N$.
Let $\sigma$ be an element of $O(N)$ defined by
$\sigma(\epsilon_j)=\epsilon_{j+1}$ for $1\le j< N$ and $\sigma(\epsilon_N)=\epsilon_1$.
Let $A_n$, $B_n$, $C_n$ and $D_n$ denote the corresponding root spaces 
given in \S \ref{sec:Dynkin} and we identify them with finite subsets of 
$\mathbb R^N$ and put $Q^i_n:=\sigma^i(Q_n)$ for $Q=A$, $B$, $C$ and $D$.
For example
\[
  A_4^3 = \langle \epsilon_4-\epsilon_5,\epsilon_5-\epsilon_6,\epsilon_7-\epsilon_8, \epsilon_8-\epsilon_9\rangle
  \subset\mathbb R^N
\]
For $\mathbf m=(m_1,m_2,\ldots)$, $\mathbf k=(k_1,k_2,\ldots)$, 
$\mathbf n=(n_1,n_2,\ldots)\in\mathbb N^{\mathbb N}$ with
\begin{equation}\label{eq:ABD}
 k_1=0\text{ and }\sum_{j=1}^\infty |m_j+k_j+n_j|<\infty
\end{equation}
define
\begin{align*}
 \Xi_{\mathbf m}
  &:= \bigcup_{j\ge 1}\bigcup_{\nu=0}^{k_j-1}
    A_j^{(j+1)\nu+\sum_{i=1}^{j-1}(i+1)m_i}
  \simeq\sum_{j\ge 1}m_jA_j,\allowdisplaybreaks\\
 M(\mathbf m)
  &:= \sum_{j\ge1}(j+1)m_j, \allowdisplaybreaks\\
 \Xi_{\mathbf m,\mathbf k}
  &=\Xi_{\mathbf m}\cup
    \bigcup_{j\ge 2}\bigcup_{\nu=0}^{k_i-1}
    D_j^{M(\mathbf m)+j\nu+\sum_{i=1}^{j-1}ik_i}\ \text{with }
    D_2=\langle \epsilon_1-\epsilon_2,\epsilon_1+\epsilon_2\rangle,
  \allowdisplaybreaks\\
 M(\mathbf m,\mathbf k)
  &:= M(\mathbf m)+\sum_{j\ge 2}jk_j,
 \allowdisplaybreaks\\
 p_D(\mathbf m,\mathbf k)
   &:=\bigl((m_1+2k_2,m_2,m_3+k_3,m_4,\ldots),(0,0,0,k_4,k_5,\ldots)\bigr),
  \allowdisplaybreaks\\
 \Xi_{\mathbf m,\mathbf k,\mathbf n}
  &:= \Xi_{\mathbf m,\mathbf k}\cup
    \bigcup_{j\ge 1}\bigcup_{\nu=0}^{k_i-1}
    B_j^{M(\mathbf m,\mathbf k)+j\nu+\sum_{i=1}^{j-1}ik_i}\ \text{with }
    B_1=\langle\epsilon_1\rangle\\
  &\underset{B_N}{\sim}\sum_{j\ge 1}m_jA_1+\sum_{j\ge 2}k_jD_j
   +\sum_{j\ge1}n_jB_j,
  \allowdisplaybreaks\\
  M(\mathbf m,\mathbf k,\mathbf n)
    &:=M(\mathbf m,\mathbf k)+\sum_{j\ge 1}jn_j=
       \sum_{j\ge1}(j+1)m_j+\sum_{j\ge1}j(k_j+n_j),
  \allowdisplaybreaks\\
  p_B(\mathbf m,\mathbf k,\mathbf n)
    &:=\bigl((m_1+n_1+2k_2,m_2,m_3+k_3,m_4,\ldots),
        (0,0,0,k_4,k_5,\ldots),\\
    &\qquad (0,n_2,n_3,\ldots)\bigr).
\end{align*}
Suppose $n\ge M(\mathbf m,\mathbf k,\mathbf n)$.
Then $\Xi_{\mathbf m,\mathbf k,\mathbf n}$ is naturally a subsystem
of $B_n$ and
\begin{equation}
 \Xi_{\mathbf m,\mathbf k,\mathbf n}^\perp\cap B_n
 \simeq k_1A_1+B_{n-M(\mathbf m,\mathbf k,\mathbf n)}
\end{equation}
 and if there exists $w\in\Aut(B_n)=W_{B_n}$ such that
\[
  \Xi_{\mathbf m,\mathbf k,\mathbf n}=
  w(\Xi_{\mathbf m',\mathbf k',\mathbf n'}),
\]
then $(\mathbf m,\mathbf k,\mathbf n)=(\mathbf m',\mathbf k',\mathbf n')$.

Fix elements $\bar{\mathbf m}=(\bar m_1,\bar m_2,\ldots)$,
$\bar{\mathbf k}=(\bar k_1,\bar k_2,\ldots)$
and $\bar{\mathbf n}=(\bar n_1,\bar n_2,\ldots)$
in $\mathbb N^{\mathbb N}$
satisfying
\begin{equation}
  \bar k_1=\bar k_2=\bar k_3=\bar n_1=0.
\end{equation}
%%%%%%%%%%%%%%%%%%%%%%%%%%   B_n  %%%%%%%%%%%%%%%%%%%%%
\begin{prop}[type $B_n$ $(n\ge 2)$]\label{prop:Bn}
Let
\begin{equation}
 \Xi_{\bar{\mathbf m},\bar{\mathbf k},\bar{\mathbf n}}
 = \sum_{j\ge1}\bar m_jA_j + \sum_{j\ge 4}\bar k_jD_j
   + \sum_{j\ge2}\bar n_jB_j.
\end{equation}
Then
\begin{align}
 \Hom(\Xi_{\bar{\mathbf m},\bar{\mathbf k},\bar{\mathbf n}},B_n)
 &=\coprod_{\substack{p_B(\mathbf m,\mathbf k,\mathbf n)
  = (\bar{\mathbf m},\bar{\mathbf k},\bar{\mathbf n})\\
   M(\bar{\mathbf m},\bar{\mathbf k},\bar{\mathbf n})\le n}}
   \Hom(\Xi_{\bar{\mathbf m},\bar{\mathbf k},\bar{\mathbf n}},B_n)_
   {(\mathbf m,\mathbf k,\mathbf n)},\\
 \Hom(\Xi_{\bar{\mathbf m},\bar{\mathbf k},\bar{\mathbf n}},B_n)_
   {(\mathbf m,\mathbf k,\mathbf n)}&:=
 \{\iota\in \Hom(\Xi_{\bar{\mathbf m},\bar{\mathbf k},\bar{\mathbf n}},B_n)
  \,;\, \text{there exists }\notag\\
  &\qquad w\in W_{B_n}\text{such that }
    w(\Xi_{\mathbf m,\mathbf k,\mathbf n})
   =\iota(\Xi_{\bar{\mathbf m},\bar{\mathbf k},\bar{\mathbf n}}\}
  \}\notag
\end{align}
and
\begin{multline}
   \Hom(\Xi_{\bar{\mathbf m},\bar{\mathbf k},\bar{\mathbf n}},B_n)_
   {(\mathbf m,\mathbf k,\mathbf n)} \ne\emptyset\\
  \quad \Leftrightarrow \ 
   p_B(\mathbf m,\mathbf k,\mathbf n)
  = (\bar{\mathbf m},\bar{\mathbf k},\bar{\mathbf n})
  \text{ and }M(\mathbf m,\mathbf k,\mathbf n)\le n.
\end{multline}
Assume $\Hom(\Xi_{\bar{\mathbf m},\bar{\mathbf k},\bar{\mathbf n}},B_n)_
   {(\mathbf m,\mathbf k,\mathbf n)} \ne\emptyset$. 
Then
\begin{gather}
   \#\bigl(W_{B_n}\!\backslash\!
    \Hom(\Xi_{\bar{\mathbf m},\bar{\mathbf k},\bar{\mathbf n}},B_n)_
   {(\mathbf m,\mathbf k,\mathbf n)}\bigr)
   =\frac{3^{k_4}\cdot(m_1+n_1+2k_2)!\cdot(m_3+k_3)!}{2^{k_2}\cdot m_1!\cdot n_1!\cdot k_2!\cdot m_3!\cdot k_3!},\\ 
  \#\bigl(W_{B_n}\!\backslash\!
    \Hom(\Xi_{\bar{\mathbf m},\bar{\mathbf k},\bar{\mathbf n}},B_n)_
   {(\mathbf m,\mathbf k,\mathbf n)}/\!
    \Aut(\Xi_{\bar{\mathbf m},\bar{\mathbf k},\bar{\mathbf n}})\bigr)
   = 1,\\
  \Xi_{\mathbf m,\mathbf k,\mathbf n}^\perp\cap B_n
   \simeq m_1A_1+B_{n-M(\mathbf m,\mathbf k,\mathbf n)},\\
  \bar \Xi_{\mathbf m,\mathbf k,\mathbf n} 
   = \Xi_{\mathbf m,\mathbf k,\mathbf n}\ \Leftrightarrow\ 
   m_2=m_3=\cdots=k_2=k_3=\cdots=0,\ 
   \textstyle\sum_{j\ge1}n_j\le 1,\\
  \Xi_{\mathbf m,\mathbf k,\mathbf n}\text{ is fundamental}
   \ \Leftrightarrow\ 
   k_2=k_3=\cdots=0\text{ and }\textstyle\sum_{j\ge1}n_j\le 1.
\end{gather}
The $S$-closure of\/ 
$\Xi_{\mathbf m,\mathbf k,\mathbf n}$ equals
$\Xi_{\mathbf m,\mathbf k,(\delta_{\nu,\sum_j jn_j})_\nu}$.
Here $\sum n_jB_j$ changes into $B_{\sum_j jn_j}$.
The $L$-closure of\/ 
$\Xi_{\mathbf m,\mathbf k,\mathbf n}$ equals
the fundamental subsystem
$\Xi_{\mathbf m,0,(\delta_{\nu,\sum_j j(k_j+n_j)})_\nu}$.
Here $\sum k_jD_j+\sum n_jB_j$ changes into $B_{\sum_j j(k_j+n_j)}$.
\end{prop}
%%%%%%%%%%%%%%%%%%%%%   C_n  %%%%%%%%%%%%%%%%%%%
Considering the dual root systems, we have the proposition
for $C_n$:
\begin{prop}[type $C_n\ (n\ge 3)$]\label{prop:Cn}
Let
\begin{equation}
 \Xi_{\bar{\mathbf m},\bar{\mathbf k},\bar{\mathbf n}}
 = \sum_{j\ge1}\bar m_jA_j + \sum_{j\ge 4}\bar k_jD_j
   + \sum_{j\ge2}\bar n_jC_j.
\end{equation}
Then the statements in Proposition~\ref{prop:Bn} 
with replacing $B_n$ and $B_{n-M(\mathbf m,\mathbf k,\mathbf n)}$
by $C_n$ and $C_{n-M(\mathbf m,\mathbf k,\mathbf n)}$, respectively,
are valid except for the last statement on $S$-closure.

The $S$-closure of this 
$\Xi_{\mathbf m,\mathbf k,\mathbf n}$
is $\Xi_{\mathbf m,0,(n_1,k_2+n_2,k_3+n_3,\ldots)}$, which is obtained by replacing $\sum k_jD_j$ by $\sum k_jC_j$.
%The $L$-closure of\/ 
%$\Xi_{\mathbf m,\mathbf k,\mathbf n}$ equals
%$\Xi_{\mathbf m,0,(\delta_{\nu,\sum_j j(k_j+n_j)})_\nu}$.
\end{prop}

We have the following propositions when $\Sigma$ 
is of type $D_n$ or of type $A_n$.
%%%%%%%%%%%%%%%%%%%%%%%  D_n  %%%%%%%%%%%%%%%%%%%%%%
\begin{prop}[type $D_n$ $(n\ge 5)$]
Let
\begin{equation}
  \Xi_{\bar{\mathbf m},\bar{\mathbf k}} =
  \sum_{j\ge 1}\bar m_jA_j + \sum_{j\ge 4}\bar k_jD_j.
\end{equation}
Then
\begin{gather}
 \Hom(\Xi_{\bar{\mathbf m},\bar{\mathbf k}},D_n)
  =\coprod_{\substack{p_D(\mathbf m,\mathbf k)
  = (\bar{\mathbf m},\bar{\mathbf k})\\
   M(\mathbf m,\mathbf k)\le n}}
   \Hom(\Xi_{\bar{\mathbf m},\bar{\mathbf k}},D_n)_
   {(\mathbf m,\mathbf k)},\\
 \Hom(\Xi_{\bar{\mathbf m},\bar{\mathbf k}},D_n)_
   {(\mathbf m,\mathbf k)}:=
 \{\iota\in \Hom(\Xi_{\bar{\mathbf m},\bar{\mathbf k}},D_n)
  \,;\, \text{there exists }w\in W_{B_n}\notag\\
  \qquad\quad \text{such that }
    w(\Xi_{\mathbf m,\mathbf k})
   =\iota(\Xi_{\bar{\mathbf m},\bar{\mathbf k}})
  \},\notag\\
  \Hom(\Xi_{\bar{\mathbf m},\bar{\mathbf k}},D_n)_
   {(\mathbf m,\mathbf k)}\ne\emptyset
  \ \Leftrightarrow \ p_D(\mathbf m,\mathbf k)
  = (\bar{\mathbf m},\bar{\mathbf k})\text{ and }
   M(\mathbf m,\mathbf k)\le n.\label{eq:Dexists}
\end{gather}
When $\Hom(\Xi_{\bar{\mathbf m},\bar{\mathbf k}},D_n)_
   {(\mathbf m,\mathbf k)}\ne\emptyset$,
\begin{gather}
   \#\bigl(W_{D_n}\!\backslash\!
     \Hom(\Xi_{\bar{\mathbf m},\bar{\mathbf k}},D_n)_
   {(\mathbf m,\mathbf k)}\bigr)
   =\varepsilon_1
    \frac{3^{k_4}\cdot(m_1+2k_2)!\cdot(m_3+k_3)!}{2^{k_2}\cdot m_1!\cdot k_2!\cdot m_3!\cdot k_3!},\\ 
  \#\bigl(W_{D_n}\!\backslash\!
   \Hom(\Xi_{\bar{\mathbf m},\bar{\mathbf k}},D_n)_
   {(\mathbf m,\mathbf k)}/\!\Aut
      (\Xi_{\bar{\mathbf m},\bar{\mathbf k}})\bigr)
   = \varepsilon_2,\\
  \#\bigl(\Aut(D_n)\backslash\!
   \Hom(\Xi_{\bar{\mathbf m},\bar{\mathbf k}},D_n)_
   {(\mathbf m,\mathbf k)}/\!\Aut
      (\Xi_{\bar{\mathbf m},\bar{\mathbf k}})\bigr)
   = 1,\\
 \Xi_{\mathbf m,\mathbf k}^\perp \simeq mA_1+D_{n-M(\mathbf m,\mathbf k)},\\
 \bar\Xi_{\mathbf m,\mathbf k}=\Xi_{\mathbf m,\mathbf k}
  \ \Leftrightarrow\ k_2=k_3=\cdots=0\text{ and }
  M(\mathbf m,\mathbf k)\ne n-1,\\
 \Xi_{\mathbf m,\mathbf k} \text{ is fundamental}
 \ \Leftrightarrow \ \textstyle\sum_{j\ge 2}k_j\le1 \ \Leftrightarrow \ 
  \Xi_{\mathbf m,\mathbf k} \text{ is $L$-closed}.
\end{gather}
Here
\begin{align*}
\varepsilon_1&=
  \begin{cases}
   2 &\text{ if \ }M(\mathbf m,\mathbf k)=n,\\
   1 &\text{ if \ }M(\mathbf m,\mathbf k)<n,
  \end{cases}\\
\varepsilon_2&=
  \begin{cases}
   2 &\text{ if \ }M(\mathbf m,\mathbf k)=n
      \text{ and }m_{2\nu}=k_{\nu+1}=0\quad (\nu=1,2,\dots),\\ 
   1 &\text{ otherwise}.
  \end{cases}
\end{align*}
The $L$-closure of $\Xi_{\mathbf m,\mathbf k}$ is obtained by 
replacing $\sum_{j\ge2}k_jD_j$ by $D_{\sum_{j\ge2}jk_j}$.
\end{prop}
\begin{prop}[type $A_n$]
Let $\Xi_{\mathbf m}=\sum_{j\ge1} m_jA_j$. Then
\begin{equation}
 \Hom(\Xi_{\mathbf m},A_n)\ne \emptyset
 \ \Leftrightarrow \ M(\mathbf m)\le n+1
\end{equation}
and if $M(\mathbf m)\le n+1$, we have
\begin{gather}
  \#\overline\Hom(\Xi_{\mathbf m},A_n) =2^{\sum_{j\ge2} m_j},
  \allowdisplaybreaks\\
  \#\bigl(\overline\Hom(\Xi_{\mathbf m},A_n)/\!\Out(\Xi_{\mathbf m})\bigr) =1,
  \allowdisplaybreaks\\
  \#\bigl(\Out(A_n)\backslash
    \overline\Hom(\Xi_{\mathbf m},A_n)\bigr) =
    \begin{cases}
        1 &(\sum_{j\ge2} m_j = 0),\\
        2^{(\sum_{j\ge2} m_j)-1}&(\sum_{j\ge2}m_j>0),\\
    \end{cases}\allowdisplaybreaks\\
  \Xi_{\mathbf m}^\perp\cap A_n\simeq A_{n-M(\mathbf m)},
  \allowdisplaybreaks\\
  \bar \Xi_{\mathbf m} = \Xi_{\mathbf m}
  \ \Leftrightarrow \ 
  \textstyle\sum_{j\ge 1} m_j\le 1\text{ and }M(\mathbf m)\ne n.
\end{gather}
Any subsystem of $A_n$ is fundamental and hence $L$-closed.
\end{prop}

%\bigskip
%\hfill
%\newpage

\subsection{Exceptional type}
($\Sigma$ : $E_6$, $E_7$, $E_8$, $F_4$, $G_2$)

%\begin{table}[h]
%\caption{$E_6$}

\begin{longtable}{|c|c|c|c|c|l|l|c|l|} \hline
$\Sigma$ & $\Xi$ & $\#$ & $\!\#_{\Xi'}\!$ & $\!\#_{\Sigma}\!$
  & $\Xi^\perp$ & $\Xi^{\perp\perp}$ & P & \\ \hline
$E_6$ & $A_1$      & 1 & 1 & 1 & $A_5$ & $\times$ & 6 &\\ \hline
$E_6$ & $A_2$      & 1 & 1 & 1 & $2A_2$ & $\times$ & 5 &\\ \hline
$E_6$ & $A_3$      & 1 & 1 & 1 & $2A_1$      & $\circ$ &  $5$  &\\ \hline
$E_6$ & $A_4$      & 2 & 1 & 1 & $A_1$       &$A_5$&  $4$    &\\ \hline
$E_6$ & $A_5$      & 2 & 1 & 1 & $A_1$       & $\times$ &  $1$       &
 $\langle\setminus2\rangle$\\ \hline
$E_6$ & $D_4$  & 1 & 1 & 1 & $\emptyset$ &$\Sigma$& 1 &\\ \hline
$E_6$ & $D_5$  & 2 & 1 & 1 & $\emptyset$ &$\Sigma$& 2 &
 $\langle\setminus1\rangle,\ 
  \langle\setminus6\rangle$\\ \hline
$E_6$ & $E_6$  & 2 & 1 & 1 & $\emptyset$ &$\Sigma$& 1 &\\ \hline
$E_6$ & $2A_1$     & 1 & 1 & 1 & $A_3$ & $\circ$ & 10 &\\ \hline
$E_6$ & $3A_1$     & 1 & 1 & 1 & $A_1$ & $A_5$ & 5 &\\ \hline
$E_6$ & $4A_1$     & 1 & 1 & 1 & $\emptyset$ & $\Sigma$ & $\to$ & $L:D_4$\\ \hline
$E_6$ & $A_2+A_1$  & 2 & 1 & 1 & $A_2$ & $2A_2$ & 10 &\\ \hline
$E_6$ & $A_2+2A_1$ & 2 & 1 & 1 & $\emptyset$ & $\Sigma$ & 5 & 
 \text{\tiny$\subset 3A_2$}\\ \hline
$E_6$ & $2A_2$     & 4 & 1 & 2 & $A_2$ & $\times$ & 1 &\\ \hline
$E_6$ & $2A_2+A_1$ & 4 & 1 & 2 & $\emptyset$ & $\Sigma$ & 1 &
 $\langle\setminus4\rangle$\ \text{\tiny$\subset 3A_2$}\\ \hline
$E_6$ & $3A_2$     & 8 & 1 & 4 & $\emptyset$ & $\Sigma$ & $\leftarrow$ &
 \S\ref{sec:4A2E8}\\ \hline
$E_6$ & $A_3+A_1$  & 2 & 1 & 1 & $A_1$       & $A_5$ &  4  &\\ \hline
$E_6$ & $A_3+2A_1$ & 2 & 1 & 1 & $\emptyset$ &$\Sigma$ & $\to$ &
  \S\ref{sec:actW},\ $L:D_5$ \\ \hline
$E_6$ & $A_4+A_1$  & 2 & 1 & 1 & $\emptyset$ &$\Sigma$& 2 &
 $\langle\setminus3\rangle,\
  \langle\setminus5\rangle$\\ \hline
$E_6$ & $A_5+A_1$  & 2 & 1 & 1 & $\emptyset$ &$\Sigma$&$\leftarrow$&\\ 
\hline
%\newpage
\hline
$\Sigma$ & $\Xi$ & $\#$ & $\!\#_{\Xi}\!$ & $\!\#_{\Xi'}\!$
  & $\Xi^\perp$ & $\Xi^{\perp\perp}$ & P & \\ \hline
$E_7$& $A_1$  & 1 & 1 & 1 & $D_6$ & $\times$ & 7 &\\ \hline
$E_7$& $A_2$  & 1 & 1 & 1 & $A_5$ & $\circ$ & 6 &\\ \hline
$E_7$& $A_3$ & 1 & 1 & 1 & $A_3+A_1$ & $\circ$ & 6 &\\ \hline
$E_7$& $A_4$ & 1 & 1 & 1 & $A_2$ & $A_5$ & 5 &\\ \hline
$E_7$& $\phantom{\ \ ]''}A_5\ \ ]''$ & 1 & 1 & 1 & $A_2$ & $\circ$ & 1 &
  $\langle2,4,5,6,7\rangle$\\ \cline{3-9}
     &$\phantom{\ \ ]'A_5}\ \ ]'$& 1 & 1 & 1 & $A_1$ & $D_6$ & 2 &
  $\langle3,4,5,6,7\rangle$\\ \hline
$E_7$& $A_6$ & 1 & 1 & 1 & $\emptyset$ & $\Sigma$ & 1 &
  $\langle\setminus2\rangle$\\ \hline
$E_7$& $A_7$ & 1 & 1 & 1 & $\emptyset$ & $\Sigma$ &$\leftarrow$&\\ \hline
$E_7$& $D_4$ & 1 & 1 & 1 & $3A_1$ & $\circ$ & 1 &\\ \hline
$E_7$& $D_5$ & 1 & 1 & 1 & $A_1$ & $D_6$ & 2 &\\ \hline
$E_7$& $D_6$ & 2 & 1 & 1 & $A_1$ & $\times$ & 1 &
  $\langle\setminus1\rangle$\\ \hline
$E_7$& $E_6$ & 1 & 1 & 1 & $\emptyset$ & $\Sigma$ & 1 &
  $\langle\setminus7\rangle$\\ \hline
$E_7$& $E_7$ & 1 & 1 & 1 & $\emptyset$ & $\Sigma$ & 1 &\\ \hline
$E_7$& $2A_1$ & 1 & 1 & 1 & $D_4+A_1$ & $\times$ & 15 &\\ \hline
$E_7$&  $\phantom{\ \ ]''}3A_1\ \ ]''$ & 1 & 1 & 1 & $D_4$ & $\circ$ & 1 &
  $\langle2,5,7\rangle$\\ \cline{3-9}
     & $\phantom{\ \ ]'3A_1}\ \ ]'$  & 1 & 1 & 1 & $4A_1$ & $\times$ & 10 &
  $\langle3,5,7\rangle$\\ \hline
$E_7$& $\phantom{\ \ ]''}4A_1\ \ ]''$& 4 & 1 & 4 & $3A_1$ & $\times$      & 2 &
  $\langle2,3,5,7\rangle$ \\ \cline{3-9}
     & $\phantom{\ \ ]'4A_1}\ \ ]'$ & 1 & 1 & 1  & $3A_1$ & $D_4$&$\leftarrow$&\\ \hline
$E_7$& $5A_1$ & 15 & 1 & 15 & $2A_1$ & $D_4+A_1$ &$\leftarrow$&\S\ref{sec:8A1E8}
\\ \hline
$E_7$& $6A_1$ & 30 & 1 & 30 & $A_1$ & $D_6$ &$\leftarrow$&\S\ref{sec:8A1E8}\\ \hline
$E_7$& $7A_1$ & 30 & 1 & 30 & $\emptyset$ & $\Sigma$ & $\leftarrow$&
 \S\ref{sec:8A1E8} \\ \hline
$E_7$& $A_2+A_1$  & 1 & 1 & 1 & $A_3$ & $A_3+A_1$ & 18 &\\ \hline
$E_7$& $A_2+2A_1$  & 1 & 1 & 1 & $A_1$ & $D_6$ & 12 &\\ \hline
$E_7$& $A_2+3A_1$  & 1 & 1 & 1 & $\emptyset$ & $\Sigma$ &1&\\ \hline
$E_7$& $2A_2$   & 2 & 1 & 1 & $A_2$ & $A_5$ & 4  & \\ \hline
$E_7$& $2A_2+A_1$ & 2 & 1 & 1 & $\emptyset$ & $\Sigma$ & 3  &
 \text{\tiny$\subset 3A_2$}\\ \hline
$E_7$& $3A_2$ & 4 & 1 & 1 & $\emptyset$ & $\Sigma$ &$\to$& \S\ref{sec:4A2E8},\ $L:E_6$\\ \hline
$E_7$& $\phantom{\ ]''}A_3+A_1\ ]''$ & 1 & 1 & 1 & $A_3$ & $\circ$ & 2 &
  $\langle2,5,6,7\rangle$\\ \cline{3-9}
     & $\phantom{\ ]'A_3+A_1}\ ]'$ & 1 & 1 & 1 & $2A_1$ & $D_4+A_1$ & 9 &
  $\langle3,5,6,7\rangle$\\ \hline
$E_7$& $\phantom{\ ]''}A_3+2A_1\ ]''$ & 2 & 1 & 2 & $A_1$ & $D_6$ & 3 & 
  {\tiny$\exists(A_3+A_1)^\perp=A_3$}\\ \cline{3-5}\cline{8-9}
     & $\phantom{\ ]'A_3+2A_1}\ ]'$ & 1 & 1 & 1 &       &       &  &
  {\tiny$\forall(A_3+A_1)^\perp=2A_1$}\\
     & & & & & & &$\to$&\text{\tiny$\subset D_3+D_2$}\,$L:D_5$
\\ \hline
$E_7$& $A_3+3A_1$ & 3 & 1 & 3 & $\emptyset$ & $\Sigma$ & &
  \text{\tiny$\subset 2A_3+A_1$} \\
     & & & & & & &$\to$& $L:D_5+A_1$\\ \hline
$E_7$& $A_3+A_2$ & 2 & 1 & 1 & $A_1$ & $D_6$ & 3
 &\text{\tiny$\subset 2A_3+A_1$} \\ \hline
$E_7$& $A_3+A_2+A_1$ & 2 & 1 & 1 & $\emptyset$ & $\Sigma$ & 1 &
  $\langle\setminus4\rangle$\ \text{\tiny$\subset 2A_3+A_1$}\\ \hline
$E_7$& $2A_3$   & 2 & 1 & 1 & $A_1$ & $D_6$ &$\leftarrow$ &
 \text{\tiny$\subset 2A_3+A_1$}\\ \hline
$E_7$& $2A_3+A_1$ & 2 & 1 & 1 & $\emptyset$ & $\Sigma$ &$\leftarrow$&
  \S\ref{sec:actW}\\ \hline
$E_7$& $A_4+A_1$ & 1 & 1 & 1 & $\emptyset$ & $\Sigma$ & 5 &\\ \hline
$E_7$& $A_4+A_2$ & 1 & 1 & 1 & $\emptyset$ & $\Sigma$ & 1 &
  $\langle\setminus5\rangle$\\ \hline
$E_7$& $\phantom{\ ]''}A_5+A_1\ ]''$ & 1 & 1 & 1 & $\emptyset$ & $\Sigma$ & 1 &
 {\tiny$A_5^\perp=A_2$}, $\langle\setminus3\rangle$ \\ \cline{3-9}
     & $\phantom{\ ]'A_5+A_1}\ ]'$ & 1 & 1 & 1 & $\emptyset$ & $\Sigma$ & $\to$
     &{\tiny$A_5^\perp=A_1$}, $L:E_6$ \\ \hline 
$E_7$& $A_5+A_2$ & 2 & 1 & 1 & $\emptyset$ & $\Sigma$ &$\leftarrow$&
 \S\ref{sec:actW}\\ \hline
$E_7$& $D_4+A_1$ & 3 & 1 & 1 & $2A_1$ & $\times$ & 1 &\\ \hline
$E_7$& $D_4+2A_1$ & 6& 1 & 1 & $A_1$ & $D_6$ &$\leftarrow$ & \\ \hline
$E_7$& $D_4+3A_1$ & 6 & 1 & 1 & $\emptyset$ & $\Sigma$ &$\leftarrow$ &
 \S\ref{sec:2D4E8} \\ \hline
$E_7$& $D_5+A_1$ & 1 & 1 & 1 & $\emptyset$ & $\Sigma$ & 1 &
  $\langle\setminus6\rangle$\\ \hline
$E_7$& $D_6+A_1$ & 2 & 1 & 1 & $\emptyset$ & $\Sigma$ &$\leftarrow$&\\ \hline\hline
$\Sigma$ & $\Xi$ & $\#$ & $\!\#_{\Xi}\!$ & $\!\#_{\Xi'}\!$  & $\Xi^\perp$ & $\Xi^{\perp\perp}$ & P & \\ 
\hline
$E_8$& $A_1$  & 1 & 1 & 1 & $E_7$ & $\circ$ & 8 &\\ \hline
$E_8$& $A_2$& 1 & 1 & 1 & $E_6$ & $\circ$ & 7  & \\ \hline
$E_8$& $A_3$& 1 & 1 & 1 & $D_5$ & $\circ$ & 7 &\\ \hline
$E_8$& $A_4$& 1 & 1 & 1 & $A_4$ & $\circ$ & 6 &\\ \hline
$E_8$& $A_5$ & 1 & 1 & 1 & $A_2+A_1$ & $\circ$ & 4 &\\ \hline
$E_8$& $A_6$ & 1 & 1 & 1 & $A_1$ & $E_7$ & 3 &\\ \hline
%\newpage\hline
$E_8$& $\phantom{\ \ ]''}A_7\ \ ]''$ & 1 & 1 & 1 & $A_1$ & $E_7$ &$\leftarrow$ &\\ \cline{3-9}
     & $\phantom{\ \ ]'A_7}\ \ ]'$ & 1 & 1 & 1 & $\emptyset$ & $\Sigma$ & 1 &
  $\langle\setminus2\rangle$\\ \hline
$E_8$& $A_8$ & 1 & 1 & 1 & $\emptyset$ & $\Sigma$ &$\leftarrow$ &\\ \hline
$E_8$& $D_4$ & 1 & 1 & 1 & $D_4$ & $\circ$ & 1 &\\ \hline
$E_8$& $D_5$ & 1 & 1 & 1 & $A_3$ & $\circ$ & 2 &\\ \hline
$E_8$& $D_6$ & 1 & 1 & 1 & $2A_1$ & $\circ$ & 1 &\\ \hline
$E_8$& $D_7$ & 1 & 1 & 1 & $\emptyset$ & $\Sigma$ & 1 &
  $\langle\setminus1\rangle$\\ \hline
$E_8$& $D_8$ & 2 & 1 & 1 & $\emptyset$ & $\Sigma$ &$\leftarrow$ &\\ \hline
$E_8$& $E_6$ & 1 & 1 & 1 & $A_2$ & $\circ$ & 1 &\\ \hline
$E_8$& $E_7$ & 1 & 1 & 1 & $A_1$ & $\circ$ & 1 &
 $\langle\setminus8\rangle$\\ \hline
$E_8$& $E_8$ & 1 & 1 & 1 & $\emptyset$ & $\Sigma$ & 1 &\\ \hline
$E_8$& $2A_1$ & 1 & 1 & 1 & $D_6$ & $\circ$ & 21 &\\ \hline
$E_8$& $3A_1$ & 1 & 1 & 1 & $D_4+A_1$ & $\circ$ & 21 &\\ \hline
$E_8$& $\phantom{\ \ ]''}4A_1\ \ ]''$ &1 & 1 & 1 & $D_4$ & $D_4$ &$\leftarrow$ &\\ \cline{3-9}
     & $\phantom{\ \ ]'4A_1}\ \ ]'$ & 1 & 1 & 1 & $4A_1$ & $\circ$ & 7 &
 $\langle2,3,6,8\rangle$\\ \hline
$E_8$& $5A_1$& 5 & 1 & 5 & $3A_1$ & $D_4+A_1$ &$\leftarrow$ &\S\ref{sec:8A1E8}\\ \hline
$E_8$& $6A_1$& 15 & 1 & 15 & $2A_1$ & $D_6$ &$\leftarrow$&\S\ref{sec:8A1E8}\\ \hline
$E_8$& $7A_1$& 30 & 1 & 30 & $A_1$ & $E_7$ &$\leftarrow$&\S\ref{sec:8A1E8}\\ \hline
$E_8$& $8A_1$& 30 & 1 & 30 & $\emptyset$ & $\Sigma$ &$\leftarrow$&
 \S\ref{sec:8A1E8}\\ \hline
$E_8$& $A_2+A_1$& 1 & 1 & 1 & $A_5$ & $\circ$ & $28$   &\\ \hline
$E_8$& $A_2+2A_1$& 1 & 1 & 1 & $A_3$ & $D_5$ & $28$   &\\ \hline
$E_8$& $A_2+3A_1$& 1 & 1 & 1 & $A_1$ & $E_7$ & $7$  &\\ \hline
$E_8$& $A_2+4A_1$& 1 & 1 & 1 & $\emptyset$ & $\Sigma$ & $\to$ &$L:A_2+D_4$\\ \hline
$E_8$& $2A_2$& 1 & 1 & 1 & $2A_2$ & $\circ$ & $8$  &\\ \hline
$E_8$& $2A_2+A_1$& 2 & 1 & 1 & $A_2$ & $E_6$ & $9$  & \\ \hline
$E_8$& $2A_2+2A_1$& 2 & 1 & 1 & $\emptyset$ & $\Sigma$ & $2$  &
  \text{\tiny$\subset 4A_2$}\\ \hline
$E_8$& $3A_2$& 4 & 1 & 1 & $A_2$ & $E_6$ & $\leftarrow$ &\\ \hline
$E_8$& $3A_2+A_1$& 4 & 1 & 1 & $\emptyset$ & $\Sigma$ &$\rightarrow$&
   \text{\tiny$\subset 4A_2$}\ \text{\small$L:E_6+A_1$}\\ \hline
$E_8$& $4A_2$& 8 & 1 & 1 & $\emptyset$ & $\Sigma$ &$\leftarrow$&
 \S\ref{sec:4A2E8}\\ \hline
%%
%$\Sigma$ & $\Xi$ & $\#$ & $\!\#_{\Xi}\!$ & $\!\#_{\Xi'}\!$
%  & $\Xi^\perp$ & $\Xi^{\perp\perp}$ & P & \\ \hline
$E_8$& $A_3+A_1$& 1 & 1 & 1 & $A_3+A_1$ & $\circ$ & $20$ &\\ \hline
$E_8$& $\phantom{\ ]''}A_3+2A_1\ ]''$& 1 & 1 & 1 & $A_3$ & $D_5$ &$\leftarrow$ &\\ \cline{3-9}
     & $\phantom{\ ]'A_3+2A_1}\ ]'$ & 1 & 1 & 1 & $2A_1$ & $D_6$ & $10$  &
 $\langle2,3,4,6,8\rangle$\\ \hline
$E_8$& $A_3+3A_1$& 3 & 1 & 3 & $A_1$ & $E_7$ &$\to$&$L:D_5+A_1$\\ \hline
$E_8$& $A_3+4A_1$& 3 & 1 & 3 & $\emptyset$ & $\Sigma$ &$\to$ &
   \text{\tiny$\subset A_3+D_5$}$\;L:D_7$\\ \hline
$E_8$& $A_3+A_2$& 1 & 1 & 1 & $2A_1$ & $D_6$ & $10$ & \\ \hline
$E_8$& $A_3+A_2+A_1$& 2 & 1 & 1 & $A_1$ & $E_7$ & $4$ & \\ \hline
$E_8$& $\!A_3+A_2+2A_1\!$& 2 & 1 & 1 & $\emptyset$ & $\Sigma$ &&
 \text{\tiny$\subset D_6+2A_1$}\\
     & & & & & & &$\to$&$L:D_5+A_2$\\ \hline
$E_8$&$\phantom{\ \ ]''}2A_3\ \ ]''$& 1 & 1 & 1 & $2A_1$ & $D_6$ &$\leftarrow$&\\ \cline{3-9}
 &$\phantom{\ \ ]'2A_3}\ \ ]'$  & 1 & 1 & 1 & $\emptyset$ & $\Sigma$ & 2 &
 $\langle2,3,4,6,7,8\rangle$\\ \hline
$E_8$& $2A_3+A_1$& 2 & 1 & 1 & $A_1$ & $E_7$ &$\leftarrow$&\\ \hline
$E_8$& $2A_3+2A_1$& 2 & 1 & 1 & $\emptyset$ & $\Sigma$ &$\leftarrow$&
 \text{\tiny$\subset D_6+2A_1$}\\ \hline
$E_8$& $A_4+A_1$& 1 & 1 & 1 & $A_2$ & $E_6$ & $12$ &\\ \hline
$E_8$& $A_4+2A_1$& 1 & 1 & 1 & $\emptyset$ & $\Sigma$ & $5$ &\\ \hline
$E_8$& $A_4+A_2$& 2 & 1 & 1 & $A_1$ & $E_7$ & 4 &\\ \hline
$E_8$& $A_4+A_2+A_1$ & 2 & 1 & 1 & $\emptyset$ & $\Sigma$ & 1 &
  $\langle\setminus4\rangle$\ \text{\tiny$\subset 2A_4$}\\ \hline
$E_8$& $A_4+A_3$ & 2 & 1 & 1 & $\emptyset$ & $\Sigma$ & 1 &
  $\langle\setminus5\rangle$\ \text{\tiny$\subset 2A_4$}\\ \hline
$E_8$& $2A_4$ & 2 & 1 & 1 & $\emptyset$ & $\Sigma$ &$\leftarrow$&
 \S\ref{sec:actW}\\ \hline
$E_8$& $\phantom{\ ]''}A_5+A_1\ ]''$ & 1 & 1 & 1 & $A_2$ & $E_6$ &$\leftarrow$ &\\ \cline{3-9}
     & $\phantom{\ ]'A_5+A_1}\ ]'$& 1 & 1 & 1 & $A_1$ & $E_7$  & 3 &
 $\langle1,4,5,6,7,8\rangle$\\ \hline
$E_8$& $A_5+2A_1$ & 2 & 1 & 2 & $\emptyset$ & $\Sigma$ &  &
   \text{\tiny$\subset A_5+A_2+A_1$}\\
     &     &   &   &   &   & & $\rightarrow$ &
   $L:E_6+A_1$\\ \hline
$E_8$& $A_5+A_2$ & 2 & 1 & 1 & $A_1$ & $E_7$ &$\leftarrow$&\\ \hline
$E_8$& $A_5+A_2+A_1$ & 2 & 1 & 1 & $\emptyset$ & $\Sigma$ &$\leftarrow$&
 \S\ref{sec:actW}\\ \hline
$E_8$& $A_6+A_1$ & 1 & 1 & 1 & $\emptyset$ & $\Sigma$ & 1 &
  $\langle\setminus3\rangle$\\ \hline
$E_8$& $A_7+A_1$ & 1 & 1 & 1 & $\emptyset$ & $\Sigma$ &$\leftarrow$&\\ \hline
$E_8$& $D_4+A_1$ & 1 & 1 & 1 & $3A_1$ & $\circ$ & 2 &\\ \hline
$E_8$& $D_4+2A_1$ & 3 & 1 & 1 & $2A_1$ & $D_6$ &$\leftarrow$&\\ \hline
$E_8$& $D_4+3A_1$ & 6 & 1 & 1 & $A_1$ & $E_7$ &$\leftarrow$&\\ \hline
$E_8$& $D_4+4A_1$ & 6 & 1 & 1 & $\emptyset$ & $\Sigma$ &$\leftarrow$&
 \S\ref{sec:2D4E8}\\ \hline
$E_8$& $D_4+A_2$ & 1 & 1 & 1 & $\emptyset$ & $\Sigma$ & 1 &\\ \hline
$E_8$& $D_4+A_3$ & 3 & 1 & 1 & $\emptyset$ & $\Sigma$ &$\to$&
    \text{\tiny$\subset2D_4$}\ $L:D_7$\\ \hline
$E_8$& $2D_4$ & 6 & 1 & 1 & $\emptyset$ & $\Sigma$ &$\leftarrow$&\S\ref{sec:2D4E8}\\ \hline
$E_8$& $D_5+A_1$ & 1 & 1 & 1 & $A_1$ & $E_7$ & 3 &\\ \hline
$E_8$& $D_5+2A_1$ & 1 & 1 & 1 & $\emptyset$ & $\Sigma$ &$\to$ &$L:D_7$\\ \hline
$E_8$& $D_5+A_2$ & 2 & 1 & 1 & $\emptyset$ & $\Sigma$ & 1 &
  $\langle\setminus6\rangle$
 \   \text{\tiny$\subset D_5+A_3$} \\ \hline
$E_8$& $D_5+A_3$ & 2 & 1 & 1 & $\emptyset$ & $\Sigma$ &$\leftarrow$ &
  \S\ref{sec:actW}\\ \hline
$E_8$& $D_6+A_1$ & 2 & 1 & 1 & $A_1$ & $E_7$ &$\leftarrow$ &\\ \hline
$E_8$& $D_6+2A_1$ & 2 & 1 & 1 & $\emptyset$ & $\Sigma$ &$\leftarrow$ &
 \S\ref{sec:actW}\\ \hline
$E_8$& $E_6+A_1$ & 1 & 1 & 1 & $\emptyset$ & $\Sigma$ & 1 &
  $\langle\setminus7\rangle$\\ \hline
$E_8$& $E_6+A_2$ & 2 & 1 & 1 & $\emptyset$ & $\Sigma$ &$\leftarrow$ &
 \S\ref{sec:actW}\\ \hline
$E_8$& $E_7+A_1$ & 1 & 1 & 1 & $\emptyset$ & $\Sigma$ &$\leftarrow$&\\ \hline
%\newpage
\hline
$\Sigma$ & $\Xi$ & $\#$ & $\!\#_{\Xi}\!$ & $\!\#_{\Xi'}\!$
  & $\Xi^\perp$ & $\Xi^{\perp\perp}$ & P & \\ \hline
$F_4$ & $A_1^L$ & 1 & 1 & 1 & $C_3$ & $\circ$ & 2 & \\ \cline{2-9}
      & $A_1^S$ & 1 & 1 & 1 & $B_3$ & $\circ$ & 2 & \\ \hline
$F_4$ & $A_2^L$ & 1 & 1 & 1 & $A_2^S$ & $\circ$ & 1 &\\ \cline{2-9}
      & $A_2^S$ & 1 & 1 & 1 & $A_2^L$ & $\circ$ & 1 &\\ \hline
$F_4$ & $A_3^L$ & 1 & 1 & 1 & $\emptyset$ & $\Sigma$ &$\to$ &$L:B_3$\\ \cline{2-9}
      & $A_3^S$ & 1 & 1 & 1 & $\emptyset$ & $\Sigma$ &$\to$&$L,\,S:C_3$ \\ \hline
$F_4$ & $D_4^L$ & 1 & 1 & 1 & $\emptyset$ & $\Sigma$ &$\leftarrow$ &\\ \cline{2-9}
      & $D_4^S$ & 1 & 1 & 1 & $\emptyset$ & $\Sigma$ &$\leftarrow$&$S:F_4$\\ \hline
$F_4$ & $B_2$ & 1 & 1 & 1 & $B_2$ & $\circ$ & 1 &\\ \hline
$F_4$ & $B_3$ & 1 & 1 & 1 & $A_1^S$ & $\circ$  & 1 &
 $\langle\setminus4\rangle$\\ \hline
$F_4$ & $C_3$ & 1 & 1 & 1 & $A_1^L$ & $\circ$ & 1 &
 $\langle\setminus1\rangle$\\ \hline
$F_4$ & $B_4$ & 1 & 1 & 1 & $\emptyset$ & $\Sigma$ &$\leftarrow$ &\\ \hline
$F_4$ & $C_4$ & 1 & 1 & 1 & $\emptyset$ & $\Sigma$ &$\leftarrow$ &$S:F_4$\\ \hline
$F_4$ & $F_4$ & 1 & 1 & 1 & $\emptyset$ & $\Sigma$ & 1 &\\ \hline
$F_4$ & $2A_1^L$ & 1 & 1 & 1 & $B_2$ &$B_2$ & $\leftarrow$ & \\ \cline{2-9}
      & $2A_1^S$ & 1 & 1 & 1 & $B_2$ & $B_2$ &$\leftarrow$& $S:B_2$ \\ \cline{2-9} 
      & $A_1^S+A_1^L$ & 1 & 1 & 1 & $A_1^L+A_1^S$ &$\times$& 4 &\\ \hline 
%%%
%\newpage\hline
$F_4$ & $3A_1^L$ & 1 & 1 & 1 & $A_1^L$ & $C_3$ &$\leftarrow$ &\\ \cline{2-9}
      & $3A_1^S$ & 1 & 1 & 1 & $A_1^S$ & $B_3$ &$\leftarrow$&$S:B_3$\\ \cline{2-9} 
      & $A_1^S+2A_1^L$ & 1 & 1 & 1 & $A_1^S$ & $B_3$ &$\leftarrow$ &\\ \cline{2-9} 
      & $2A_1^S+A_1^L$ & 1 & 1 & 1 & $A_1^L$ & $C_3$ &$\leftarrow$&$S:B_2+A_1^L$ \\ \hline 
$F_4$ & $4A_1^L$ & 1 & 1 & 1 & $\emptyset$ & $\Sigma$ &$\leftarrow$ & \\ \cline{2-9}
      & $4A_1^S$ & 1 & 1 & 1 & $\emptyset$ & $\Sigma$ &$\leftarrow$&$S:F_4$ \\ \cline{2-9}
      & $2A_1^S+2A_1^L$ & 1 & 1 & 1 & $\emptyset$ & $\Sigma$ &$\leftarrow$& $S:B_2+2A_1^L$\\ \hline 
$F_4$ & $A_2^L+A_1^S$ & 1 & 1 & 1 & $\emptyset$ & $\Sigma$ & 1 &
 $\langle\setminus3\rangle$\\ \cline{2-9}
      & $A_2^S+A_1^L$ & 1 & 1 & 1 & $\emptyset$ & $\Sigma$ & 1 &
 $\langle\setminus2\rangle$\\ \hline
$F_4$ & $A_2^S+A_2^L$ & 1 & 1 & 1 & $\emptyset$ & $\Sigma$ &$\leftarrow$ &\\ \hline
$F_4$ & $B_2+A_1^L$ & 1 & 1 & 1 & $A_1^L$ & $C_3$ &$\leftarrow$ & \\ \cline{2-9}
      & $B_2+A_1^S$ & 1 & 1 & 1 & $A_1^S$ & $B_3$ &$\leftarrow$&$S:B_3$  \\ \hline
$F_4$ & $B_2+2A_1^L$ & 1 & 1 & 1 & $\emptyset$ & $\Sigma$ &$\leftarrow$ & \\ \cline{2-9}
      & $B_2+2A_1^S$ & 1 & 1 & 1 & $\emptyset$ & $\Sigma$ &$\leftarrow$&$S:B_4$  \\ \hline
$F_4$ & $2B_2$ & 1 & 1 & 1 & $\emptyset$ & $\Sigma$ &$\leftarrow$  &$S:B_4$\\ \hline
$F_4$ & $A_3^S+A_1^L$ & 1 & 1 & 1 & $\emptyset$ & $\Sigma$ &$\leftarrow$ &$S:C_3+A_1^L$\\ \hline
$F_4$ & $A_3^L+A_1^S$ & 1 & 1 & 1 & $\emptyset$ & $\Sigma$ &$\leftarrow$ &\\ \hline
$F_4$ & $C_3+A_1^L$ & 1 & 1 & 1 & $\emptyset$ & $\Sigma$ &$\leftarrow$ &\\ \hline
$F_4$ & $B_3+A_1^S$ & 1 & 1 & 1 & $\emptyset$ & $\Sigma$ &$\leftarrow$&$S:B_4$ \\ \hline
\hline
%$\Sigma$ & $\Xi$ & $\#$ & $\!\#_{\Xi}\!$ & $\!\#_{\Xi'}\!$ & $\Xi^\perp$ & $\Xi^{\perp\perp}$ & P & \\ 
\hline
$G_2$ & $A_1^L$ & 1 & 1 & 1 & $A_1^S$ & $\circ$ & 1 &
 $\langle\setminus2\rangle$\\ \cline{2-9}
      & $A_1^S$ & 1 & 1 & 1 & $A_1^L$ & $\circ$ & 1 &
 $\langle\setminus1\rangle$\\ \hline
$G_2$ & $A_2^L$ & 1 & 1 & 1 & $\emptyset$ & $\Sigma$ &$\leftarrow$ &\\ \cline{2-9}
      & $A_2^S$ & 1 & 1 & 1 & $\emptyset$ & $\Sigma$ &$\leftarrow$ &$S:G_2$\\ \hline
$G_2$ & $G_2$ & 1 & 1 & 1 & $\emptyset$ & $\Sigma$ & 1 &\\ \hline
$G_2$ & $A_1^S+A_1^L$ & 1 & 1 & 1 & $\emptyset$ & $\Sigma$ &$\leftarrow$ &\\ \hline
\end{longtable}
We explain some symbols used in the above table.
\begin{rem}\label{rem:last}
{\rm i)}
In the table we use following notation.
\begin{gather*}
  \Sigma^L:=\{\alpha\in\Sigma\,;\,
  |\beta|\le|\alpha|\quad(\forall\beta\in\Sigma)\},\allowdisplaybreaks\\
  A_m^S\simeq A_m^L\simeq A_m,\quad
  A_m^L\subset \Sigma^L,\quad A_m^S\cap\Sigma^L=\emptyset,\allowdisplaybreaks\\
  D_m^S\simeq D_m^L\simeq D_m,\quad
  D_m^L\subset \Sigma^L,\quad D_m^S\cap\Sigma^L=\emptyset.
\end{gather*}

{\rm ii)} The symbols \ $]'$ \ and \ $]''$ in the column $\Sigma$.\\
Suppose $\Sigma$ is irreducible and of exceptional type.
Then $\#\overline\Hom(\Xi,\Sigma)/\!\Out(\Xi)\le 2$.
When $\#\overline\Hom(\Xi,\Sigma)/\!\Out(\Xi)=2$,
$\Sigma$ is of type $E_7$ or $E_8$ and then
the symbols $[\Xi]'$ and $[\Xi]''$ are used in \cite{Dy} to distinguish 
the equivalence classes of the imbeddings $\Xi\subset\Sigma$.
Then $[\Xi]'$ means that there is a representative $\Xi$
in the equivalence class such that 
\begin{equation}
 \Xi\subset A_n\subset \Sigma=E_n
\end{equation}
with $n=7$ or $8$.
For example, $\#\overline\Hom(4A_1,E_7)/\!\Out(4A_1)=2$ and the symbols
$[4A_1]'$ and $[4A_1]''$ are used in \cite{Dy}, which are expressed
by $\ ]'\ $ and $\ ]''\ $ respectively in the column $\Sigma$ in our table
(cf.~\eqref{eq:4A14}).

In \cite{Dy} the distinction of the elements of 
$\Out(\Sigma)\backslash \overline\Hom(\Xi,\Sigma)/\!\Out(\Xi)$ such as
$\ ]'\ $ and $\ ]''\ $ is not discussed but it is stated there that 
the distinction is due to actual calculation.

{\rm iii)} The structure of $\Out_\Sigma(\Xi)$.\\
If $(\#)=\#\Out(\Xi)$ or $(\#)=1$ in the table,
it follows from \eqref{eq:NOut} that
$\#\Out_\Sigma(\Xi)=1$ or $\Out_\Sigma(\Xi)\overset{\sim}{\to}
\Out(\Xi)$, respectively.
In the column $P$ in the table, a reference such as 
\S\ref{sec:8A1E8} gives the description of $\Out_\Sigma(\Xi)$ 
for other non-trivial cases.

If $\Xi=\Xi_1+\Xi_2\subset\Xi'=\Xi_1+\Xi_1^\perp\subset\Sigma$
and $\Out(\Xi)\overset{\sim}{\leftarrow}\Out(\Xi_1)\times\Out(\Xi_2)$
and $\Xi^\perp=\emptyset$,
we have
\begin{equation}
  \Out_\Sigma(\Xi)\simeq N_{\Aut_{\Sigma}(\Xi')}(\Xi_2)/W_{\Xi}.
\end{equation}
The symbol ``$\subset\Xi'$" is indicated in the column $P$
if $\Out_\Sigma(\Xi)$ is easily obtained by this relation.
For example, $\Out_{E_8}(D_5+A_2)$ is isomorphic to
$\Out_{E_8}(D_5+A_3)$ through the imbedding 
$D_5+A_2\subset D_5+A_3\subset E_8$.
\end{rem}

%%%%%%%%%%%%%%%%%%%%%%%

%%%%%%%%%%%%%%%%%%%%%%%
\end{document}